\pdfoutput=1
\documentclass[11pt]{article}
\usepackage[letterpaper,margin=1in]{geometry}
\usepackage[T1]{fontenc}
\usepackage{lmodern}
\usepackage{amsmath,amssymb,amsthm,mathtools,bm}
\usepackage{booktabs,array,tabularx,longtable,float}
\usepackage{graphicx}
\usepackage{enumitem}
\usepackage[authoryear,round]{natbib}
\usepackage{xcolor}
\usepackage{microtype}
\usepackage[colorlinks=true,linkcolor=blue!50!black,citecolor=blue!50!black,urlcolor=blue!50!black]{hyperref}
\usepackage[capitalise,noabbrev]{cleveref}
\usepackage{caption}
\newtheorem{theorem}{Theorem}
\newtheorem{proposition}{Proposition}
\newtheorem{corollary}{Corollary}
\newtheorem{lemma}{Lemma}[section]
\theoremstyle{definition}
\newtheorem{definition}{Definition}
\newtheorem{assumption}{Assumption}

\theoremstyle{remark}
\newtheorem{remark}{Remark}
\crefname{assumption}{Assumption}{Assumptions}

\newcommand{\R}{\mathbb{R}}
\newcommand{\Z}{\mathbb{Z}}
\newcommand{\N}{\mathbb{N}}
\newcommand{\E}{\mathbb{E}}
\newcommand{\Prob}{\mathbb{P}}
\newcommand{\Var}{\operatorname{Var}}
\newcommand{\Cov}{\operatorname{Cov}}
\newcommand{\CM}{\mathcal{H}}
\newcommand{\cA}{\mathcal{A}}
\newcommand{\cI}{\mathcal{I}}
\newcommand{\cN}{\mathcal{N}}
\newcommand{\cF}{\mathcal{F}}
\newcommand{\cX}{\mathcal{X}}

\newcommand{\cQ}{\mathcal{Q}}

\newcommand{\ind}{\mathbf{1}}
\newcommand{\aff}{\operatorname{Aff}}
\newcommand{\TV}{\mathrm{TV}}
\newcommand{\equivm}{\sim}
\newcommand{\perpm}{\perp}
\newcommand{\abs}[1]{\lvert #1\rvert}
\newcommand{\norm}[1]{\lVert #1\rVert}
\newcommand{\inner}[2]{\langle #1,#2\rangle}
\newcommand{\proofidea}{\noindent\emph{Proof idea.}\ }

\title{The Spatial Confounding Endgame:\\ Resolution, Roughness, and the Limits of Effect Estimation}
\author{Xiao Huang\\ Department of Environmental Sciences, Emory University\\ Atlanta, GA 30329\\ Email: xiao.huang2@emory.edu}
\date{November 2026}

\begin{document}
\maketitle

\begin{abstract}
Across the spatial sciences, covariate effects are estimated from geo-referenced data, and the central threat is unmeasured, spatially structured confounding. When such data identify an effect is known only for stationary fields in Euclidean space, and which of some thirty proposed adjustments recover it, and how precisely, has remained open. In an operator framework, CAR and SAR models, kriging, eigenvector spatial filtering, spatially varying coefficient models and linear graph neural networks become spectral functions of one Laplace-type operator, the graph models exactly on regular graphs. Under spectral infill and vanishing fine-scale confounding, the effect is consistently estimable if and only if the covariate lies almost surely outside the Cameron--Martin space of the residual confounder; for Whittle--Mat\'ern fields of spectral dimension $d_s$ this reads $\nu_X\le\nu_V+d_s/2$. On a fixed regular lattice with an unknown residual spectrum, every value of the effect, of either sign, is observationally equivalent to the true one. Under stated conditions the frontier extends to nonlinear effects, non-Gaussian residuals and counts, and with a covariate measured at a fixed network of stations not even the sign of the effect is identified. Strictly inside the frontier, GLS without a working nugget is consistent if and only if its working smoothness is at least $\nu_X-d_s/2$; OLS, restricted spatial regression and, on refining regular graphs, Leroux, proper CAR and SAR models with fixed hyperparameters never are, whereas on regular grids the intrinsic CAR model obeys the GLS window. Minimax rates quantify the remaining precision. Resolution, roughness and fine-scale confounding, not the choice among nonparametric spatial models, decide whether an effect can be learned, and the last cannot be tested from spatial data. A census of some thirty methods places each within these results.
\end{abstract}

\noindent\textbf{Keywords:} spatial confounding; identifiability; spatial resolution and scale; infill asymptotics; Cameron--Martin space; equivalence of Gaussian measures; aliasing; spectral adjustment; graph Laplacian; conditional autoregressive models; restricted spatial regression.

\section{Introduction}\label{sec:intro}

Across the spatial sciences, researchers estimate how a spatially varying covariate affects a response from data indexed by location: land cover and urban form on surface temperature, green space on house prices, habitat quality on species abundance, deprivation or air pollution on health. We write $X$ for the covariate of interest (an exposure in epidemiology, a treatment in causal inference, an explanatory variable in regression), $Y$ for the response and $U$ for an unmeasured confounder. The central threat is that $U$ varies smoothly in space and affects both $X$ and $Y$. Adding a spatially correlated random effect to the regression is the standard response, and it can change the coefficient of $X$ dramatically \citep{Clayton1993,Reich2006}. Since \citet{HodgesReich2010} asked whether the spatial random effect ``messes up the fixed effect you love'', a large literature has debated what spatial models do under confounding, whether restricted spatial regression (RSR) should be used \citep{Reich2006,HughesHaran2013,Hanks2015,KhanCalder2022,ZimmermanVerHoef2022}, how the scales of covariate and confounder matter \citep{Paciorek2010,Page2017,KellerSzpiro2020}, and how to adjust \citep{ThadenKneib2018,SchnellPapadogeorgou2020,Dupont2022,Guan2023,Dupont2023,Woodward2024,Burman2025}. Much of this literature is set in spatial epidemiology and environmental health, and we borrow examples from that field where they help, but nothing in the theory is specific to it.

The debate has gradually sharpened into an identification question: \emph{when can a single realization of spatial data identify the effect of a covariate at all?} Two recent results settle important special cases. \citet{Gilbert2025} proved that generalized least squares (GLS) with a Mat\'ern (or squared-exponential) working covariance and a positive working nugget is consistent under infill asymptotics whenever the covariate contains some non-spatial variation. \citet{DattaStein2025} characterized, for stationary Gaussian fields in $\R^d$, when the slope is consistently estimable. Under infill asymptotics and a cross-covariance smoother than the covariate, it is estimable if the covariate is less than $d/2$ orders smoother than the confounder and not if it is more than $d/2$ orders smoother, with range and variance parameters playing no role. They settle the boundary case for $d=1$ and integer smoothness. They also cover measurement error, irregular designs, non-stationary and heavy-tailed scale-mixture fields, and show that GLS with an exponential working covariance can be inconsistent. The same $d/2$ threshold appears without confounding in \citet{Yu2022}. \citet{BolinWallin2026} analysed GLS in an operator framework on compact metric spaces, with eigenbasis observations and the Cameron--Martin space of the noise. They proved that GLS converges to zero when the covariate is rougher than the signal it drives (spatial self-confounding). \citet{Guan2023} formulated the now-standard identifying assumption that confounding ``dissipates at local scales''.

\paragraph{What is at stake.} Regression with a spatial random effect, a spatial filter or a smooth function of location is the default analysis of geo-referenced data: in disease mapping, where the BYM model is fitted routinely \citep{BYM1991,Riebler2016,Rue2009}; in environmental epidemiology and ecology \citep{Dormann2007}; in spatial econometrics and regional science \citep{Anselin1988,LeSagePace2009}; and in GIScience, through geographically weighted regression, eigenvector filtering and, increasingly, GeoAI learners (\cref{sec:literature}). In all these fields the coefficient of a covariate is read as an effect, often a causal one. Since \citet{HodgesReich2010} the literature has proposed some thirty adjustment methods and variants (\cref{tab:census}), and simulation and empirical comparisons rank them differently from one setting to the next \citep{Reich2021review,Urdangarin2023,Tec2024space,Pim2026review}. The results just described say when such data identify the effect for stationary fields in $\R^d$. What has been missing is a theory that also covers fixed resolution, non-Euclidean supports and the estimators in use, and that says how precisely the effect can be recovered. This paper develops such a theory within operator-stationary models. Its answer turns on three properties of the design and of the data-generating process (resolution, roughness and fine-scale confounding) rather than on the choice among nonparametric spatial models.

These results leave four questions open, each of practical consequence.
\begin{enumerate}[label=(\roman*),leftmargin=2em,itemsep=1pt]
\item \emph{Resolution.} Most socioeconomic, health and environmental data are areal and come at a fixed resolution (census tracts, counties, satellite pixels). Does the high-frequency assumption identify the effect when the data cannot resolve high frequencies?
\item \emph{Beyond linear--Gaussian.} Covariate effects are often nonlinear, residual confounding need not be Gaussian, responses are often counts, and supports are often networks or the sphere.
\item \emph{Estimators.} Which of the many methods in use reach the identification frontier, and which never do?
\item \emph{Rates.} How much precision does deconfounding cost, and what limits it?
\end{enumerate}

\paragraph{Contributions.} We address all four within one framework (question (i) on regular lattices, and nonlinear confounding mechanisms in (ii) under Cameron--Martin negligibility).
\begin{enumerate}[leftmargin=2em,itemsep=2pt]
\item \textbf{Spectral normal form} (\cref{prop:normalform}, \cref{tab:zoo}). CAR, Leroux and SAR models, Whittle--Mat\'ern/SPDE fields and kriging, graph Mat\'ern processes, Moran eigenvector filtering, (M)GWR-type spatially varying coefficient models in Moran-eigenvector form and linear graph neural networks are spectral functions of a single Laplace-type operator $L$ (the graph-based models exactly on regular graphs and approximately otherwise). Every estimator built from functions of $L$ is a ratio of spectrally weighted quadratic forms, a unifying restatement of known spectral bias formulas. Each estimator therefore targets a weighted average of a \emph{scale-specific association} $\gamma(\lambda)=\beta+\alpha(\lambda)$, where $\alpha(\lambda)$ is the confounding at frequency $\lambda$.
\item \textbf{The Cameron--Martin frontier} (\cref{thm:frontier}, \cref{cor:matern}). Under spectral infill and vanishing high-frequency confounding, $\beta$ is consistently estimable if and only if $\cI=\sum_j f_j/g_j=\infty$, where $f_j$ and $g_j$ are the spectra of the covariate and of the residual confounder; equivalently, the covariate lies almost surely outside the Cameron--Martin space of the residual confounder. For Whittle--Mat\'ern fields defined through the SPDE on supports of spectral dimension $d_s$ (bounded domains, the sphere, river and road networks) this is $\nu_X\le\nu_V+d_s/2$, boundary included; the Euclidean, $\R^d$-stationary analogue is the threshold of \citet{DattaStein2025}. The operator framework and the spectral-infill scheme follow \citet{BolinWallin2026}. We add a jointly stationary confounder, impossibility for every estimator, and a dichotomy for arbitrary confounding classes.
\item \textbf{A resolution theorem} (\cref{thm:resolution}, \cref{cor:latticetarget}). At a fixed lattice resolution with an unknown residual spectrum, \emph{every} effect $\beta'$, of either sign, is observationally equivalent to the true one. Some continuum model with effect $\beta'$ and vanishing high-frequency confounding induces \emph{identical} laws for lattice or block-averaged data, however large the study region. For point samples the same holds with the residual spectrum known whenever the alias variance of the confounding is bounded below or the nugget variance is unknown. With the residual spectrum and the nugget variance both known, $\beta$ is identified exactly when that alias variance comes arbitrarily close to zero. Identification at fixed resolution thus rests on exact knowledge of second-order nuisance quantities. The graph-native assumption that the highest graph frequency is unconfounded holds if and only if the alias-weighted confounding at the Nyquist frequency vanishes, which continuum assumptions do not imply.
\item \textbf{Which estimators reach the frontier} (\cref{thm:window}, \cref{prop:estimators}). Under spectral infill without a nugget and strictly inside the frontier, GLS with a fixed working spectrum is consistent if and only if its working smoothness is at least $\nu_X-d_s/2$. Undersmoothing below this bound is fatal. Oversmoothing does not affect consistency, although with a nugget it can slow convergence severely; at the frontier only the correct tail works. This identifies, in the operator setting, the exact smoothness boundary behind the GLS inconsistency exhibited by \citet{DattaStein2025}. Ordinary least squares (OLS), RSR and eigenvector filtering with a fixed number of eigenvectors are never consistent. On refining graphs of bounded degree (regular graphs for proper CAR and SAR), Leroux, proper CAR and SAR models with fixed hyperparameters are not consistent either, because their precisions stay bounded away from zero at low frequencies, whereas on regular grids the intrinsic CAR model, whose GLS estimator regresses neighbour differences, obeys the window with working smoothness $1-d_s/2$ (\cref{cor:areal}). GLS with a correctly specified tail estimates the finest-scale association $\beta+\alpha(\infty)$, which recovers, within the confounding framework, the self-confounding result of \citet{BolinWallin2026}.
\item \textbf{Minimax rates and the resolution floor} (\cref{thm:minimax}). With $N$ resolved modes, the minimax risk is of order $N^{-\min(2r,p)}$; its bias floor is of the order of the squared worst-case confounding at the finest resolved scale.
\item \textbf{Beyond the linear--Gaussian model} (\cref{thm:shift,prop:shifts,thm:nonlinear,thm:channels}, \cref{cor:kriging}). A shift-equivalence principle transfers the identification frontier to nonlinear covariate effects under Cameron--Martin-negligible confounding (a strong restriction) and to non-Gaussian residual confounders with independent, identically distributed coefficients from a density that is positive almost everywhere and has finite Fisher information (more generally, under the uniform conditions of \cref{prop:shifts}(b)). Its impossibility half also covers arbitrary observation channels. Consistency extends to Poisson counts with growing expected counts; for GLS it requires finite residual variance, and with Cauchy residuals GLS fails on part of the frontier. When the covariate is too smooth, not even the \emph{shape} of a smooth response function is identified, provided the admissible shifts of the residual contain the Sobolev space of order $\nu_V+d/2$, as for restrictions of $\R^d$-stationary Mat\'ern-type fields. When the covariate is measured only at a fixed network of stations, the coefficient on a kriged covariate surface is not identified for the usual kernel choices, and neither is the sign of the structural effect, whereas classical measurement error is asymptotically harmless.
\end{enumerate}
Simulations in \cref{sec:experiments} illustrate each prediction; for the impossibility results they show that natural estimators fail as predicted.

\paragraph{The endgame in brief.} Within the operator-stationary model, purely spatial data identify a covariate effect when three things hold at once. The data must keep resolving finer scales, the covariate must be rough relative to what remains of the confounder after removing its covariate-predictable part, and confounding must die out at fine scales. Each condition is also necessary in a specific sense. Without roughness, no estimator is consistent, even without confounding (\cref{thm:frontier}(i)). Without refining resolution, vanishing high-frequency confounding loses all identifying power on regular lattices when the residual spectrum is unknown (\cref{thm:resolution}). Without vanishing confounding, any confounding class closed under recentring leaves $\beta$ unidentified (\cref{sec:frontier}). Within the nonparametric class of vanishing high-frequency confounding, the choice among the estimators of \cref{tab:zoo}---CAR or kriging, Bayesian or frequentist, GLS or a graph neural network---cannot change where the frontier lies; it decides whether and how fast an estimator reaches it (\cref{thm:window,prop:estimators,cor:areal}). Parametric or graph-native assumptions can restore identification, for instance by asserting how confounding behaves at some resolved scale (that it vanishes at the highest graph frequency, say, or that the association follows the functional form of SLX and lag models) or by tying confounding below the grain to confounding above it; such assertions go beyond vanishing fine-scale confounding, and the data can check them at most in part. \Cref{tab:endgame} summarizes the verdicts.

\paragraph{A guide for spatial scientists.} Three ideas carry every result, and each has a familiar geographic reading. \emph{Scale}: the linear spatial models in common use act on the data scale by scale (exactly on regular graphs, approximately otherwise), as a set of weights on the eigenfunctions (spatial patterns) of a Laplacian of the study region or graph; the scale-specific association $\gamma_j$ is the covariate--response relation seen at one scale. \emph{Roughness}: an effect can be separated from a broad-scale confounder only through the local variation of the covariate, and the Cameron--Martin energy $\cI$ measures how much local variation there is relative to the residual spatial variation of the response. \emph{Resolution}: local variation is visible only down to the grain of the data, and variation below the grain does not disappear but folds back (aliases) onto the resolved scales, much as in the modifiable areal unit problem the grain of the units shapes what an analysis can see. \Cref{tab:notation} translates the notation into these terms. \Cref{sec:literature} reviews the models, adjustment methods and theory in the literature and records, method by method, what the results imply (\cref{tab:census}). Every result is followed by the idea of its proof. \Cref{app:proofs} gives complete proofs; it opens with a map of the proofs and of the classical results they use (\cref{tab:proofmap}), and states every external result it relies on.

\begin{table}[!tbp]
\centering\small
\caption{The endgame at a glance. ``Identified'' means consistently estimable by some procedure; the last column names procedures that attain it. $\Delta\nu=\nu_X-\nu_V$, $d_s$ = spectral dimension of the support.}\label{tab:endgame}
\begin{tabularx}{\textwidth}{>{\raggedright\arraybackslash}p{4.6cm}>{\raggedright\arraybackslash}p{2.0cm}X}
\toprule
Situation & Verdict & Reference and remedy \\
\midrule
Refining resolution, $\Delta\nu<d_s/2$, confounding vanishes at high frequency & identified & \cref{thm:frontier}; under spectral infill, GLS without nugget and with working smoothness $\ge\nu_X-d_s/2$ (on regular grids, the intrinsic CAR model when $\nu_X\le1$, \cref{cor:areal}), high-pass GLS with the correct tail, eigenvector filtering with $K_N\to\infty$ and $K_N/N\to0$ \\
Refining resolution, $\Delta\nu=d_s/2$ & identified, variance of order $1/\log N$ & \cref{cor:matern}; GLS or high-pass GLS with the correct tail \\
Refining resolution, $\Delta\nu>d_s/2$ & not identified & \cref{thm:frontier}(i); no spatial method; external information needed \\
Confounding not vanishing at high frequency & not identified & \cref{thm:frontier}(iii); GLS inside its window estimates $\beta+\alpha(\infty)$ (\cref{thm:window,prop:whatgls}) \\
Fixed lattice resolution (point values or block means), any domain size, residual spectrum unknown & not identified, not even the sign & \cref{thm:resolution}; identification only through functional-form assumptions on $\alpha$ \\
Covariate measured at a fixed network of stations (e.g.\ kriged surface) & sign not identified & \cref{cor:kriging}; the coefficient on a kriged surface is not identified when $2\nu_\ast>\nu_V$ \\
Nonlinear covariate effect (Cameron--Martin-negligible confounding); non-Gaussian residual with i.i.d.\ coefficients of finite Fisher information; Poisson counts with growing expected counts per unit & same identification frontier & \cref{thm:nonlinear,prop:shifts,thm:channels}; GLS needs finite residual variance \\
Under spectral infill: OLS, RSR, eigenvector filtering with fixed $K$, GLS with working smoothness below $\nu_X-d_s/2$ or with a working nugget, penalized spatial+, and Leroux, proper CAR and SAR models with fixed hyperparameters on refining regular graphs & not consistent, even when $\alpha=0$ & \cref{thm:window,prop:estimators,cor:areal} \\
\bottomrule
\end{tabularx}
\end{table}

\paragraph{What the results settle.} \Cref{tab:settled} lists ten questions that the literature has raised, several of them repeatedly, and the answer each receives here. Three features distinguish these answers from further comparisons. Their impossibility parts are \emph{estimator-free}: they hold for every measurable function of the data, from OLS to graph neural networks, so no method that uses the same data under the same nonparametric assumptions can escape them; escaping them requires additional information (time, instruments, negative controls) or stronger, for example parametric, assumptions. They are \emph{sharp} within the operator-stationary model: each positive result is matched by a negative one at the same boundary, and the minimax rates show that what remains is a question of precision, with an explicit floor; the boundaries that remain open (the $\R^d$-stationary boundary, GLS with a nugget under point sampling, estimators that attain the frontier for heavy-tailed residuals) are listed in \cref{sec:discussion}. And they are \emph{partly checkable}: resolution is a property of the design; the roughness gap can be estimated from the covariate and the response residuals (\cref{sec:practice}), with the imprecision usual for smoothness parameters; and the one remaining assumption, vanishing fine-scale confounding, is exactly the one that spatial data cannot test at any resolution, because $(\beta,\alpha)$ and $(\beta+\alpha(\infty),\alpha-\alpha(\infty))$ induce the same law (\cref{sec:frontier}).

\begin{table}[!tbp]
\centering\small\setlength{\tabcolsep}{4pt}
\caption{Questions raised in the literature and the answers given here.}\label{tab:settled}
\begin{tabularx}{\textwidth}{>{\raggedright\arraybackslash}p{5.6cm}X}
\toprule
Question (where it was raised) & Answer (result) \\
\midrule
Does a spatial random effect remove confounding bias or add to it? \citep{Clayton1993,Reich2006,HodgesReich2010,Paciorek2010} & Under spectral infill and strictly inside the frontier, GLS with a working smoothness of at least $\nu_X-d_s/2$ converges in probability to the finest-scale association $\beta+\alpha(\infty)$: it removes confounding exactly when confounding vanishes at fine scales, and when the confounder and the covariate's unconfounded variation are equally smooth its bias can even exceed that of OLS (\cref{thm:window,prop:whatgls}). On refining regular graphs, Leroux, proper CAR and SAR models with fixed hyperparameters converge to the limit of OLS, and the intrinsic CAR model deconfounds only if $\nu_X\le1$ (\cref{cor:areal}) \\
Should restricted spatial regression be used? \citep{Reich2006,HughesHaran2013,Hanks2015,KhanCalder2022,ZimmermanVerHoef2022} & No, if the target is the structural effect: in Gaussian linear models its point estimate is OLS, which is never consistent under spectral infill (\cref{prop:estimators}(a)); see \citet{Bradley2024} for a defence under a different estimand \\
How much smoother than the confounder may the covariate be? \citep{Paciorek2010,Page2017,DattaStein2025} & By at most half the spectral dimension, $\nu_X\le\nu_V+d_s/2$, boundary included for SPDE-defined (operator-stationary) fields, on domains, the sphere and networks (\cref{cor:matern}); beyond it no estimator is consistent, even without confounding (\cref{thm:frontier}(i)) \\
Does the assumption that confounding dissipates at local scales identify the effect? \citep{Guan2023,KellerSzpiro2020,Dupont2022} & With refining resolution, if and only if $\cI=\infty$ (\cref{thm:frontier}); on a fixed regular grid with an unknown residual spectrum, no: every value of $\beta$, of either sign, is observationally equivalent (\cref{thm:resolution}(a)); with known $g$ and nugget variance see \cref{thm:resolution}(b),(c) \\
Does the continuum assumption imply an unconfounded highest graph frequency? \citep{Guan2023,Prim2025} & No: on regular lattices the latter holds if and only if the alias-weighted confounding at the Nyquist frequency vanishes (\cref{cor:latticetarget}) \\
When is GLS with a misspecified covariance consistent? \citep{Gilbert2025,DattaStein2025,BolinWallin2026} & Under spectral infill without a nugget and strictly inside the frontier, if and only if its working smoothness is at least $\nu_X-d_s/2$, so oversmoothing and mild undersmoothing are safe; at the frontier only the correct tail works (\cref{thm:window}). Point-sampled GLS with a nugget is not covered; \citet{Gilbert2025} give a sufficient condition \\
How precisely can the effect be estimated, and what limits precision? & The minimax risk is of order $N^{-\min(2r,p)}$, with a floor set by the confounding at the finest resolved scale (\cref{thm:minimax}) \\
Do nonlinear effects, non-Gaussian residuals, counts or machine-learning estimators change the answer? \citep{Woodward2024,Tec2024space} & Not for identification, under the stated conditions (nonlinear effects under Cameron--Martin-negligible confounding, non-Gaussian residuals with i.i.d.\ coefficients of finite Fisher information, counts with growing expected counts per unit): the frontier is the same (\cref{thm:nonlinear,prop:shifts,thm:channels}), and every impossibility result holds for every measurable estimator \\
Can a covariate surface interpolated from a fixed network of stations be used? \citep{Gryparis2009,SzpiroPaciorek2013} & Not to identify the effect: the coefficient on a kriged surface is not identified when $2\nu_\ast>\nu_V$, and the sign of the structural effect is never identified (\cref{cor:kriging}) \\
Why do comparisons of adjustment methods disagree? \citep{Urdangarin2023,Tec2024space,Pim2026review} & Our results suggest one reason: designs differ in their position relative to the frontier and in resolution. Inside the frontier the consistent methods differ only in efficiency, while OLS, RSR, fixed-$K$ filtering and bounded-weight estimators remain inconsistent; outside it all fail (\cref{sec:practice}). Reporting $\hat\nu_X-\hat\nu_V$ and the resolution of each design would test this explanation \\
\bottomrule
\end{tabularx}
\end{table}
\section{The literature: models, adjustment methods and theory}\label{sec:literature}

This section reviews the spatial models in current use, the methods proposed to adjust for spatial confounding and the theory that analyses them, and places each method within the framework of this paper. \Cref{tab:census} summarizes the verdicts. Recent reviews and comparisons are \citet{Reich2021review}, \citet{Urdangarin2023}, \citet{KhanBerrett2023}, \citet{AkbariWinterTomko2023} and \citet{Pim2026review}; \citet{Hodges2013} gives a book-length treatment of random effects that compete with fixed effects, and \citet{Tec2024space} provide a benchmark environment.

\subsection{The spatial confounding debate}\label{sec:lit-debate}

The phenomenon was first noticed in disease mapping. \citet{Clayton1993} observed that adding spatially structured random effects changed the estimated covariate effect, which they interpreted as adjustment for ``confounding by location''. \citet{Reich2006} and \citet{HodgesReich2010} showed that the random effect and the fixed effect compete for the same low-frequency variation, and \citet{Wakefield2007} stressed that the data alone cannot determine the form and extent of the spatial dependence that such adjustments rely on. Ecology had its own version of the debate: spatial and non-spatial regressions of species distributions gave systematically different coefficients \citep{Dormann2007,Bini2009}. \citet{Paciorek2010} gave the first bias formulas for geostatistical models, showing that spatial models reduce bias only when the covariate varies at finer scales than the confounder. \citet{Page2017} extended the analysis to prediction, \citet{Nobre2021} to hierarchical models, \citet{HuiBondell2022} to generalized estimating equations, and \citet{Narcisi2024} expressed the bias through quadratic forms. \citet{KhanBerrett2023} separated confounding of the analysis model from confounding of the data-generating process, and \citet{Dupont2023} derived a unified analytical framework for the bias, identifying smoothing as the mechanism. Simulation and empirical comparisons \citep{Reich2021review,Urdangarin2023,Tec2024space,Pim2026review} rank the adjustment methods differently across settings; \cref{sec:practice} suggests why.

\subsection{Spatial models and learners in current use}\label{sec:lit-models}

\emph{Areal models.} The intrinsic CAR model \citep{Besag1974}, its convolution with an unstructured effect \citep{BYM1991} and the reparametrized BYM2 model \citep{Riebler2016}, the Leroux model \citep{Leroux2000} and proper CAR models are the defaults of disease mapping, fitted by Markov chain Monte Carlo or by integrated nested Laplace approximation \citep{Rue2009}; \citet{BesagMondal2005} relate the intrinsic CAR model on refining lattices to the de Wijs process. Spatial econometrics adds the simultaneous autoregressive (SAR) error model, the spatial lag model, the SLX model and the spatial Durbin model \citep{Anselin1988,LeSagePace2009,Elhorst2014,HalleckVegaElhorst2015}, estimated by maximum likelihood or by spatial two-stage least squares \citep{KelejianPrucha1998}; \citet{Manski1993} and \citet{GibbonsOverman2012} stress that such models identify effects only through assumed functional forms.

\emph{Geostatistical models.} Kriging and Gaussian-process regression with Mat\'ern covariances \citep{Stein1999} are fitted at scale through the SPDE representation \citep{Lindgren2011}, nearest-neighbour and Vecchia approximations \citep{Vecchia1988,Datta2016}, low-rank and multi-resolution methods \citep{Banerjee2008,CressieJohannesson2008,Nychka2015,Katzfuss2017}; \citet{Heaton2019} compare them. Penalized splines and generalized additive models with smooth functions of location \citep{Wood2017} are a closely related alternative: a penalized smooth is the posterior mean under a Gaussian prior \citep{KimeldorfWahba1970}.

\emph{Basis and local methods.} Moran eigenvector spatial filtering \citep{Griffith2003,Dray2006}, Getis filtering \citep{GetisGriffith2002}, random-effects and shrinkage versions \citep{MurakamiGriffith2015,Donegan2020}, geographically weighted regression and its multiscale version \citep{Brunsdon1996,Fotheringham2002,Fotheringham2017}, and spatially varying coefficient models \citep{Gelfand2003,Murakami2017,Murakami2019scale} are the workhorses of quantitative geography.

\emph{Machine learning and GeoAI.} Graph neural networks \citep{Defferrard2016,Kipf2017,Hamilton2017,Velickovic2018,Wu2019SGC}, which act on graph signals through spectral filters (GCN, ChebNet, SGC; \citealp{Shuman2013,Ortega2018}) or through spatial-domain aggregation (GAT, GraphSAGE), are increasingly used for areal data, together with positional encoders \citep{Mai2020space2vec,Mai2022review,Klemmer2023}, neural networks for geostatistical data \citep{Chen2024deepkriging,ZhanDatta2025}, spatial random forests and boosting \citep{Georganos2021,Saha2023,Sigrist2022}, neural versions of geographically weighted regression \citep{Du2020gnnwr,HagenauerHelbich2022} and post hoc explanations of fitted learners \citep{Li2022shap,Li2024geoshapley}; \citet{WikleZammitMangion2023} review the area.

\Cref{prop:normalform} writes the members of these families that are built from functions of a single operator (exactly on regular graphs) as spectral functions of that operator. Geographically weighted regression, nearest-neighbour approximations and the nonlinear learners fall outside the normal form, but every impossibility result in this paper applies to them, because those results constrain all measurable estimators.

\subsection{Adjustment methods}\label{sec:lit-methods}

\emph{Orthogonality constraints.} Restricted spatial regression constrains the random effect to the orthogonal complement of the covariates \citep{Reich2006,HughesHaran2013}, with geostatistical \citep{Hanks2015}, hierarchical \citep{Nobre2021} and spatio-temporal \citep{Adin2023} versions and a restricted working correlation for estimating equations \citep{HuiBondell2022}. \citet{Hanks2015} and \citet{KhanCalder2022} showed that in Gaussian linear models the RSR point estimate equals the non-spatial one and that its intervals are anticonservative, \citet{ZimmermanVerHoef2022} argued against its use, and \citet{Bradley2024} defends it under a different interpretation of the estimand. Variants project the areal centroids onto the orthogonal complement of the covariates and rebuild the neighbourhood graph \citep[SPOCK and MSPOCK;][]{Prates2019,Azevedo2021}, also for frailty models \citep{Azevedo2023frailty}, replace the constraint by group-lasso shrinkage \citep{Hefley2017}, or use transformed Gaussian Markov random fields \citep{Prates2015} and fixed-rank kriging representations \citep[reviewed by][]{Pim2026review}.

\emph{Residualizing the covariate.} Spatial+ replaces the covariate by its residual from a spatial smooth and then fits a partial spline \citep{Dupont2022,Rice1986}; the geoadditive structural equation model regresses the spatial residuals of the response on those of the covariate \citep{ThadenKneib2018}, which is the double-residual estimator of \citet{Robinson1988} and \citet{Speckman1988}. Bayesian \citep{MarquesWiemann2023}, nonlinear \citep{DupontAugustin2024}, areal and simplified multivariate \citep{Urdangarin2024} versions of spatial+ exist. \citet{Marques2022} model the dependence between the covariate and the spatial effect explicitly, through correlated Gaussian random fields.

\emph{Scale selection and spectral adjustment.} \citet{KellerSzpiro2020} remove the covariate's projection on the smoothest basis functions and select the number of functions; spike-and-slab priors on principal spline functions automate the choice \citep{Zaccardi2025}. The spectral adjustment of \citet{Guan2023} assumes that the coherence between covariate and confounder vanishes at high frequencies, either parametrically (a Mat\'ern coherence) or semiparametrically (in their CAR version, at the largest graph eigenvalue); \citet{Prim2025} extend it to several covariates and responses under a local unconfoundedness assumption, and \citet{WuBanerjee2025} study multivariate areal models. \citet{Woodward2024} show that many of these methods are instrumental-variable procedures in which small-scale covariate variation is the instrument. Basis voting \citep{Burman2025} assumes only that a plurality of basis functions carry no confounding.

\emph{Design-based and causal machine-learning approaches.} Distance-adjusted propensity score matching \citep{Papadogeorgou2019} and spatial propensity scores with doubly robust estimation \citep{Davis2019} treat spatial proximity as a proxy for the confounder. \citet{Gilbert2021causal} formalize location as a proxy within a causal framework and propose double machine learning \citep{Chernozhukov2018}, \citet{Gilbert2025} prove consistency of common estimators when the covariate has variation not explained by location, and \citet{Wiecha2024} develop two-stage estimators with spatial nuisance fits. In economics, spatial first differences \citep{DruckenmillerHsiang2018} and geographic regression discontinuities \citep{KeeleTitiunik2015} exploit local comparisons. Other approaches bring in information beyond a single spatial snapshot: learned representations of observed non-local confounder fields \citep{Tec2023weather2vec}, time-invariant latent confounding in spatio-temporal data \citep{Christiansen2022}, and spatio-temporal treatment processes \citep{Papadogeorgou2022airstrikes}. Interference is modelled with generalized propensity scores under measured confounding \citep{Giffin2023}, with Bayesian models for confounding and interference \citep{PapadogeorgouSamanta2023}, and as a multi-cause signal from which a substitute confounder is reconstructed \citep{Khot2025}. In geography, spatial stratified heterogeneity measures \citep{Wang2010geodetector,WangZhangFu2016} and geographical convergent cross mapping \citep{Gao2023GCCM} are widely used to screen for causal factors.

\emph{Inference under spatial dependence.} Spatial heteroskedasticity- and autocorrelation-consistent standard errors \citep{Conley1999} and spatial correlation robust inference \citep{MullerWatson2022} address uncertainty rather than bias. \citet{MullerWatson2024} and \citet{ConleyKelly2025} show that strongly persistent spatial data produce spurious significance even with such corrections, and \citet{MullerWatson2024} propose transformations that difference away the persistent components.

\subsection{Identification theory}\label{sec:lit-theory}

Population-level identification under parametric covariance models was studied by \citet{SchnellPapadogeorgou2020} and \citet{Tang2026}. On the asymptotic side, \citet{Gilbert2025} proved consistency of GLS with a Mat\'ern or squared-exponential working covariance and a positive working nugget when the covariate has non-spatial variation. \citet{DattaStein2025} characterized consistent estimability for stationary Gaussian fields in $\R^d$ under infill asymptotics: estimable if the covariate is less than $d/2$ orders smoother than the confounder, not if it is more. \citet{Yu2022} found the same threshold without confounding. \citet{BolinWallin2026} analysed GLS in an operator framework, with eigenbasis observations and the Cameron--Martin space of the noise, and proved that GLS converges to zero for covariates rougher than the signal they drive; \citet{Bolin2026} removes this self-confounding by matching the regression operator to the smoothing. \citet{Pan2026} studies partial identification in additive Poisson risk models, bounding residual confounding with a benchmark covariate and a negative-control outcome, and \citet{Engh2026} analyses GLS with averages over refining partitions. \citet{Guan2023} noted that frequencies above $\pi/\Delta s$ cannot be observed from gridded data \citep[see also][]{FuentesReich2010}, but no nonparametric identification or non-identification result for fixed resolution was available.

\subsection{Scale, resolution and aliasing in geography}\label{sec:lit-scale}

Geography has long recognized that the scale and zoning of areal units shape statistical results: the modifiable areal unit problem \citep{GehlkeBiehl1934,OpenshawTaylor1979,Openshaw1984,FotheringhamWong1991,Wong2009}, recently revisited through local models \citep{FotheringhamSachdeva2022}, the ecological fallacy \citep{Robinson1950}, the uncertain geographic context problem \citep{Kwan2012,Kwan2018} and the change-of-support problem \citep{GotwayYoung2002}. Remote sensing and geostatistics treat the support and resolution of measurements explicitly \citep{WoodcockStrahler1987,AtkinsonTate2000}, and \citet{LamQuattrochi1992} linked scale and resolution to roughness through fractal dimension. Sampling theory says that a regular grid cannot resolve frequencies above its Nyquist frequency, which fold back onto lower ones \citep{Shannon1949}. \citet{Tobler1970} and \citet{Goodchild2004} discuss laws of geography; \citet{KolakAnselin2020} and \citet{AkbariWinterTomko2023} review causal inference from a spatial perspective. These literatures describe how results change with scale. The resolution theorem of \cref{sec:resolution} turns this into an identification limit: on a regular grid with an unknown residual spectrum, no amount of territory and no model choice recovers the effect unless one assumes something about confounding at the scale of the units, or a functional form that links confounding below the grain to confounding above it.

\subsection{Where this paper stands}\label{sec:lit-new}

Our frontier theorem is the operator-theoretic counterpart of the threshold of \citet{DattaStein2025}, and our operator framework, the spectral-infill scheme (O1) and the Cameron--Martin criterion for GLS follow \citet{BolinWallin2026} \citep[see also][]{KirchnerBolin2022,BolinKirchner2023,LiTangBanerjee2023}. \citet{DattaStein2025} already treat measurement error, irregular designs, non-stationary fields and heavy-tailed scale-mixture fields in the $\R^d$-stationary setting; our \cref{prop:shifts}(c) and \cref{thm:channels}(b)(i)--(ii) are operator-setting counterparts of their Corollary~8 and Theorem~6. To our knowledge the following are new: the resolution theorem; the consistency window of GLS and its consequences for areal models, which with fixed hyperparameters do not deconfound on refining regular graphs unless they are intrinsic; the minimax theory; the frontier for nonlinear covariate effects, for residual fields with independent non-Gaussian coefficients and for count channels; and the structural non-identification with fixed networks of stations. Mathematically we rely on the classical theory of equivalence of measures \citep{Kakutani1948,Feldman1958,Hajek1958,Shepp1965,IbragimovRozanov1978,Bogachev1998}; on Grenander's observation that a regression coefficient is estimable without error from a continuous record if and only if the regressor lies outside the reproducing kernel Hilbert space of the noise \citep{Grenander1950,Parzen1961}; and on fixed-domain asymptotics \citep{Stein1999,Zhang2004,BolinKirchner2023}.

{\footnotesize
\begin{longtable}{>{\raggedright\arraybackslash}p{4.7cm}>{\raggedright\arraybackslash}p{4.1cm}>{\raggedright\arraybackslash}p{5.8cm}}
\caption{A census of methods in current use and the verdicts of the theory. The impossibility results, \cref{thm:frontier}(i), \cref{thm:resolution}, \cref{thm:shift}(a), \cref{thm:channels}(a) and \cref{cor:kriging}, hold for every measurable estimator and therefore apply to every row; the verdicts record what is specific to each method. ``(O1)'' is spectral infill and ``fixed grid'' is scheme (O3); verdicts for models with hyperparameters assume that these are fixed or, where stated, rescaled at a prescribed rate.}\label{tab:census}\\
\toprule
Family and methods & Spectral reading & Verdict \\
\midrule
\endfirsthead
\toprule
Family and methods & Spectral reading & Verdict \\
\midrule
\endhead
\midrule\multicolumn{3}{r}{\emph{continued}}\\
\endfoot
\bottomrule
\endlastfoot
\multicolumn{3}{l}{\emph{Spatial models and learners}}\\
Areal random effects: ICAR, BYM, BYM2, Leroux, proper CAR, SAR error \citep{Besag1974,BYM1991,Riebler2016,Leroux2000,Anselin1988} & GLS with a working precision $p(\lambda)$ on the graph spectrum (\cref{tab:zoo}); ICAR-GLS is the regression of neighbour differences (\cref{cor:areal}(c)) & On refining graphs of bounded degree (regular for proper CAR and SAR), Leroux, proper CAR and SAR error precisions are bounded above and away from zero, so GLS is not consistent and, when the eigenvalues of fixed modes tend to zero, converges to the limit of OLS (\cref{cor:areal}(a)); on graphs with Weyl-type eigenvalue bounds, such as regular grids, ICAR and BYM have $p(\lambda)\asymp\lambda$ and obey the window with $\nu_h=1-d_s/2$: consistent if and only if $\nu_X<1$ and $\Delta\nu<d_s/2$, or $\nu_X=1$ and $\Delta\nu\le d_s/2$ (\cref{cor:areal}(b)); with $0\le1-\rho_n\le C'n^{-2/d_s}$, Leroux and proper CAR obey the same window, and with $0\le1-\rho_nk\le C'n^{-2/d_s}$ SAR obeys the window with $\nu_h=2-d_s/2$; BYM2, which rescales with the graph, is covered by neither case \\
Geostatistics: kriging and GPs, SPDE--INLA, NNGP and Vecchia, predictive process, FRK, LatticeKrig, MRA \citep{Stein1999,Lindgren2011,Rue2009,Vecchia1988,Datta2016,Banerjee2008,CressieJohannesson2008,Nychka2015,Katzfuss2017} & GLS with an exact or approximated Mat\'ern spectrum; NNGP and Vecchia approximate the covariance, low-rank versions (predictive process, FRK) have a working spectrum that vanishes beyond their rank, and LatticeKrig and MRA are multi-resolution & Strictly inside the frontier, exact Mat\'ern GLS without a nugget is consistent under (O1) if and only if its working smoothness is at least $\nu_X-d_s/2$, and at the frontier if and only if it equals $\nu_V$ (\cref{thm:window}); approximations that preserve the spectral tail inherit this by analogy only, as they are not functions of $L$; versions with a working nugget, including low-rank versions with their nugget idealized as spectral truncations, have bounded weights under (O1) and are not consistent there; a sieve GLS with a working spectrum tail-equivalent to $g$ is consistent for point samples (\cref{prop:pointsamples}); GLS with a nugget under point sampling is not covered (experiment E9), and \citet{Gilbert2025} prove its consistency when the covariate has non-spatial variation, a case outside our $L^2$ framework \\
Penalized splines and GAMs with smooth functions of location \citep{Wood2017} & a penalized smooth with a residual variance is GLS with a power-law spectral prior and a nugget \citep{KimeldorfWahba1970}; its weights $1-t_j$ are bounded (spatial+ with $s_j=0$) & As GLS with a nugget: not consistent under (O1), not covered under point sampling; with a fixed basis dimension, a fixed-rank filter \\
Eigenvector spatial filtering: Moran ESF, RE-ESF, shrinkage ESF \citep{Griffith2003,Dray2006,MurakamiGriffith2015,Donegan2020} & flat prior on the $K$ smoothest modes; RE-ESF and shrinkage ESF add random effects on selected modes and a residual variance & Fixed $K$: never consistent (\cref{prop:estimators}(b)); $K_N\to\infty$ with $K_N/N\to0$: consistent if and only if $\Delta\nu<d_s/2$, with noise variance of order $K_N^{\Delta'-1}$ (\cref{prop:estimators}(c)); RE-ESF and shrinkage ESF have bounded weights (the nugget case) and are not consistent under (O1); data-dependent, stepwise eigenvector selection is not covered. Getis filtering \citep{GetisGriffith2002}, which transforms each variable with the local $G_i$ statistic, is outside the normal form \\
Spatial econometrics: SEM, spatial lag, SLX, spatial Durbin, GS2SLS \citep{Anselin1988,LeSagePace2009,Elhorst2014,HalleckVegaElhorst2015,KelejianPrucha1998,GibbonsOverman2012} & SEM is GLS with a SAR spectrum. On $k$-regular graphs lag, SLX and Durbin models impose a parametric scale profile: $\gamma_j=\beta+\theta(k-\lambda_j)$ (SLX), $\beta/\{1-\rho(k-\lambda_j)\}$ (lag) and $\{\beta+\theta(k-\lambda_j)\}/\{1-\rho(k-\lambda_j)\}$ (Durbin), with $k-\lambda_j$ replaced by $1-\lambda_j/k$ for row-standardized $W$. In all three, $\beta=\gamma(k)$, the association at the scale with zero Moran eigenvalue & SEM: as the SAR error model in the first row. Lag, SLX and Durbin models identify $\beta$ through this functional form across the resolved band, which in particular equates $\beta$ with $\gamma(k)$, the association at the scale with zero Moran eigenvalue, and not through vanishing fine-scale confounding. If the profile reflects genuine spillovers, which our model excludes, $\beta$ and $\theta$ describe direct and indirect effects \citep{LeSagePace2009}, and our results do not cover that case \\
Local and varying-coefficient regression: GWR, MGWR, eigenvector and GP-based SVC models, GNNWR, GWANN \citep{Brunsdon1996,Fotheringham2002,Fotheringham2017,Gelfand2003,Murakami2017,Murakami2019scale,Du2020gnnwr,HagenauerHelbich2022} & local or spatially varying slopes; the Moran-eigenvector form is a spectral SVC prior (\cref{tab:zoo}) & A constant effect is a special case, so the impossibility results bind every functional of $\beta(s)$ that reduces to $\beta$ for a constant effect, such as its average \\
Graph neural networks and GeoAI: GCN, ChebNet, SGC, GAT, GraphSAGE, positional encoders, DeepKriging, NN-GLS, RF-GLS, GPBoost, geographical random forests, boosting with SHAP or GeoShapley \citep{Kipf2017,Defferrard2016,Wu2019SGC,Velickovic2018,Hamilton2017,Mai2020space2vec,Klemmer2023,Chen2024deepkriging,ZhanDatta2025,Saha2023,Sigrist2022,Georganos2021,Li2022shap,Li2024geoshapley} & GCN and SGC filters $p(L)X$ are low-pass (\cref{tab:zoo}); ChebNet learns general polynomial filters; GAT, GraphSAGE and the other learners are nonlinear or aggregate in the spatial domain & Regression on low-pass features concentrates weight on the confounded scales (\cref{sec:estimators}); NN-GLS and RF-GLS are consistent without confounding \citep{Saha2023,ZhanDatta2025}; consistency under confounding not analysed \\
\multicolumn{3}{l}{\emph{Adjustment methods for spatial confounding}}\\
Restricted spatial regression and its geostatistical, hierarchical and spatio-temporal versions \citep{Reich2006,HughesHaran2013,Hanks2015,Nobre2021,Adin2023} & random effect orthogonal to $X$; in Gaussian linear models the point estimate equals OLS \citep{Hanks2015,KhanCalder2022}, in spatial GLMMs approximately the non-spatial estimate & Never consistent in the Gaussian linear model, even when $\alpha=0$ (\cref{prop:estimators}(a)); \citet{Bradley2024} defends RSR under a different estimand \\
Displacement, shrinkage and restricted-correlation variants: SPOCK, MSPOCK and frailty versions, group lasso, restricted GEE, TGMRF and FRK-based adjustments \citep{Prates2019,Azevedo2021,Azevedo2023frailty,Hefley2017,HuiBondell2022,Prates2015,Pim2026review} & modify the spatial basis, its prior or the working correlation & Outside the normal form in general; consistency not analysed \\
Residualizing the covariate: spatial+ and its Bayesian, nonlinear, areal and multivariate versions; gSEM; double spatial regression \citep{Dupont2022,MarquesWiemann2023,DupontAugustin2024,Urdangarin2024,ThadenKneib2018,Wiecha2024} & spatial+: weights $\big((1-t_j)(1-s_j),(1-t_j)(1-s_j)^2\big)$; gSEM (double residual): $\big((1-s_j)(1-t_j),(1-s_j)^2\big)$ (\cref{tab:zoo}) & With a first-stage smoother whose eigenvalues $s_j$ tend to zero and a GLS second stage, the window of \cref{thm:window} (\cref{prop:estimators}(e)); with a penalized second stage, and for gSEM with fixed smoothers, the weights are bounded and these estimators are not consistent under (O1); with covariate variation not explained by location, outside our framework, two-stage estimators are $\sqrt n$-consistent \citep{Wiecha2024} \\
Parametric joint models: correlated Gaussian random fields, Mat\'ern-coherence spectral adjustment, affine estimators, Bayesian models for confounding with interference \citep{Marques2022,Guan2023,SchnellPapadogeorgou2020,Tang2026,PapadogeorgouSamanta2023} & a parametric form for the coherence $\alpha(\lambda)$ or the cross-covariance & Identification through functional form; at fixed resolution this ties confounding below the grain to confounding above it (\cref{sec:resolution}) \\
Scale selection and high-pass adjustment: Keller--Szpiro, spike-and-slab principal splines, high-pass GLS (an idealization of spectral adjustment), multivariate areal adjustment, the IV framework \citep{KellerSzpiro2020,Zaccardi2025,WuBanerjee2025,Woodward2024} & discard or downweight the low frequencies; weights $\ind\{j>K\}$ (flat) or $\ind\{j>K\}/h_j$ (GLS-weighted) & With refining resolution: flat weights with a deterministic cut-off $K_N\to\infty$, $K_N/N\to0$ (Keller--Szpiro, eigenvector filtering) give consistency if and only if $\Delta\nu<d_s/2$ (\cref{prop:estimators}(c)); GLS-weighted high-pass with the correct tail is consistent if the retained modes carry divergent Cameron--Martin energy and inconsistent if that energy stays bounded (\cref{prop:estimators}(d)). On a fixed grid the high-frequency assumption has no identifying power (\cref{thm:resolution}) \\
Lattice- and graph-native assumptions: unconfounded highest graph frequency (semiparametric spectral adjustment), local unconfoundedness and spatial first differences, basis voting \citep{Guan2023,Prim2025,DruckenmillerHsiang2018,Burman2025} & assumptions on the alias-averaged confounding $\bar\alpha^\delta$ itself & On regular lattices they identify $\beta$, because lattice data identify $\beta+\bar\alpha^\delta$ (\cref{cor:latticetarget}); the highest-frequency and neighbour-difference assumptions concern confounding at the scale of the units, while basis voting needs unconfounded support on a plurality of resolved modes, at any scale; none is implied by continuum assumptions. ICAR-GLS is the all-edge version of spatial first differences (\cref{cor:areal}(c)) \\
Proxies for the confounder: distance-adjusted and spatial propensity scores, doubly robust and DML estimators, location as a proxy, interference-aware deconfounding \citep{Papadogeorgou2019,Davis2019,Gilbert2021causal,Gilbert2025,Wiecha2024,Chernozhukov2018,Khot2025} & spatial proximity, smooth functions of location or neighbours' treatments stand in for the confounder & For purely spatial data in our $L^2$ framework the impossibility results bind. Covariate variation not explained by location \citep{Gilbert2025} is a white-noise component of the covariate ($f_j$ bounded below), outside that framework; there even GLS with a nugget is consistent. The single-snapshot identification of \citet{Khot2025} rests on a multi-cause structure induced by interference, outside our no-interference model \\
Information beyond the purely spatial data: temporal variation with time-invariant latent confounders, observed confounder fields, spatio-temporal treatment processes, measured confounding with interference \citep{Christiansen2022,Tec2023weather2vec,Papadogeorgou2022airstrikes,Giffin2023} & add variables or structure outside the purely spatial model & Compatible with the theory, which shows that such information is necessary outside the frontier and on a fixed grid (\cref{thm:frontier}(i), \cref{thm:resolution}) \\
Geographical screening of causal factors: geographical detectors, geographical convergent cross mapping \citep{Wang2010geodetector,WangZhangFu2016,Gao2023GCCM} & functionals of the spatial data & When read as effect estimates, bound by the impossibility results \\
Inference under spatial dependence: spatial HAC, spatial correlation robust inference, spatial unit-root tests \citep{Conley1999,MullerWatson2022,MullerWatson2024,ConleyKelly2025} & correct standard errors \citep{Conley1999,MullerWatson2022} or difference away persistent components \citep{MullerWatson2024} & Complementary: when $\cI<\infty$, or on a fixed grid with an unknown residual spectrum, the effect is not identified, so no correction of the uncertainty can recover it \\
\end{longtable}
}
\section{Framework}\label{sec:framework}

\subsection{Spatial support and spectral dimension}\label{sec:support}

Let $(M,\mu)$ be a compact metric measure space and $L$ a non-negative self-adjoint operator on $L^2(M,\mu)$ with compact resolvent. Its eigenvalues $0\le\lambda_1\le\lambda_2\le\cdots\to\infty$ come with a real orthonormal basis of eigenfunctions $(\varphi_j)_{j\ge1}$. We say that $L$ has \emph{spectral dimension} $d_s$ if $\lambda_j\asymp j^{2/d_s}$ (two-sided bounds for all large $j$). By Weyl's law this holds with $d_s=d$ for the Neumann or Dirichlet Laplacian on a bounded Lipschitz domain $D\subset\R^d$ and for the Laplace--Beltrami operator on a compact $d$-dimensional manifold (the sphere $S^2$ for global data), and with $d_s=1$ for the Kirchhoff Laplacian on a compact metric graph, which models river and road networks \citep{BolinSimasWallin2024}. On a finite graph with $n$ nodes, $L$ is a graph Laplacian with $n$ eigenpairs, and the asymptotic assumptions above concern a sequence of graphs; finite-sample statements such as the bounds of \cref{thm:minimax} then apply with $N=n$ when the spectral bounds hold uniformly in $n$. Graph models are spectral functions of $D-W$ (some only approximately on irregular graphs; \cref{tab:zoo}), but how they behave as graphs refine depends on how their hyperparameters scale with $n$: with fixed hyperparameters, on graphs of bounded degree (regular graphs for proper CAR and SAR), the precisions of the Leroux, proper CAR, SAR and graph Mat\'ern models stay bounded away from zero and infinity on the graph spectrum, and only intrinsic models downweight the smoothest modes without limit (\cref{cor:areal}).

\subsection{Operator-stationary fields}\label{sec:fields}

\begin{definition}[$L$-stationarity]\label{def:Lstat}
A centred random field $Z$ with $Z\in L^2(M)$ almost surely is \emph{$L$-stationary with spectrum $(s_j)$} if its coefficients $z_j=\inner{Z}{\varphi_j}$ are uncorrelated with $\Var(z_j)=s_j$ and $\sum_j s_j<\infty$. Two fields are \emph{jointly $L$-stationary} if the coefficient pairs $(z_j,w_j)$ are uncorrelated across $j$.
\end{definition}

For Gaussian fields, $L$-stationarity means that the coefficients are independent, that is, that the covariance operator is diagonal in the eigenbasis of $L$. The leading example is the \emph{Whittle--Mat\'ern field} $(\kappa^2+L)^{(\nu+d_s/2)/2}Z=\tau\,\mathcal W$ driven by Gaussian white noise $\mathcal W$, whose spectrum is $s_j=\tau^2(\kappa^2+\lambda_j)^{-(\nu+d_s/2)}$; we call $\nu>0$ its smoothness \citep{Whittle1954,Lindgren2011,BolinKirchner2023}. On $\R^d$ this is the Mat\'ern field with smoothness parameter $\nu$.

\begin{remark}[For Whittle--Mat\'ern fields Tobler's first law holds automatically, but it is not what drives identification]
For $L=-\Delta$ on $\R^d$ the Whittle--Mat\'ern covariance is positive and strictly decreasing in distance: for this model class, correlation decaying with distance follows from the form of the generating operator, which is local when $(\nu+d_s/2)/2$ is an integer and is approximated by sparse Markov fields otherwise, rather than being an extra assumption \citep{Whittle1954,Lindgren2011,Tobler1970}. Other local operators, such as the oscillating SPDE models of \citet{Lindgren2011}, produce covariances that change sign. The results below show that identification under confounding is governed by something else: the high-frequency tails of the spectra, that is, by local roughness.
\end{remark}

\subsection{The spectral normal form of the model zoo}\label{sec:zoo}

\Cref{tab:zoo} lists the spectral form of the models and estimators that dominate applied spatial analysis. Its content is collected in the following proposition, whose proof is elementary algebra (\cref{app:normalform}).

\begin{table}[!tbp]
\centering\small
\caption{Spectral normal form. Upper block: models as spectral priors $h(L)$. ICAR and Leroux CAR are exact on every graph with $L=D-W$; proper CAR, SAR and the Moran-eigenvector models are exact on $k$-regular graphs and only approximate on irregular ones; for the graph Mat\'ern GP the exponent $\nu$ plays the role of $\nu+d_s/2$. How graph models behave as graphs refine depends on how their hyperparameters scale with the number of nodes (\cref{cor:areal}). Lower block: estimators $\hat\beta=\sum_j\omega_jx_jy_j/\sum_j\tilde\omega_jx_j^2$ in the eigenbasis of $L$. Here $s_j,t_j\in[0,1]$ are the eigenvalues of linear smoothers and $p$ is a graph filter, whose weights may change sign. Under (O1) the weights of the nugget ($\sigma^2>0$), spatial+ and gSEM rows are bounded (\cref{sec:estimators}).}\label{tab:zoo}
\begin{tabularx}{\textwidth}{>{\raggedright\arraybackslash}p{4.9cm}>{\raggedright\arraybackslash}p{2.6cm}X}
\toprule
Model & Operator $L$ & Spectral prior $h(\lambda)$ \\
\midrule
ICAR \citep{Besag1974,BYM1991} & $D-W$ & $(\tau\lambda)^{-1}$ on $\mathbf 1^\perp$ \\
Leroux CAR \citep{Leroux2000} & $D-W$ & $\tau^{-1}(1-\rho+\rho\lambda)^{-1}$ \\
Proper CAR, degree $k$ & $kI-W$ & $\tau^{-1}(k(1-\rho)+\rho\lambda)^{-1}$ \\
SAR, symmetric $W$, degree $k$ \citep{VerHoef2018} & $kI-W$ & $\sigma^2(1-\rho k+\rho\lambda)^{-2}$ \\
Whittle--Mat\'ern / SPDE; kriging \citep{Lindgren2011} & $-\Delta$ on $M$ & $\tau^2(\kappa^2+\lambda)^{-(\nu+d_s/2)}$ \\
Graph Mat\'ern GP \citep{Borovitskiy2021} & graph Laplacian & $\sigma^2(2\nu/\kappa^2+\lambda)^{-\nu}$ \\
Moran eigenvector filtering \citep{Griffith2003,Dray2006} & via $MCM$ & flat on the $K$ smoothest modes, $0$ otherwise \\
Random-effects ESF \citep{MurakamiGriffith2015} & via $MCM$ & $\sigma^2\mu_j^{\gamma}$ on selected modes ($\mu_j$: Moran eigenvalues) \\
(M)GWR-type SVC in Moran-eigenvector form \citep{Murakami2017,Fotheringham2017} & any of the above & $\beta_\ell(s)=b_\ell+\sum_j\theta_{\ell j}\varphi_j(s)$, $\theta_{\ell j}\sim N(0,h_\ell(\lambda_j))$ \\
Linear GNN / SGC \citep{Wu2019SGC,Defferrard2016} & normalized Laplacian & filter $p(\lambda)$, e.g.\ $(1-\lambda)^K$ \\
\midrule
Estimator & \multicolumn{2}{l}{Weights $(\omega_j,\tilde\omega_j)$} \\
\midrule
OLS & \multicolumn{2}{l}{$(1,1)$} \\
RSR \citep{Reich2006,HughesHaran2013}; its point estimate equals OLS \citep{Hanks2015,KhanCalder2022} & \multicolumn{2}{l}{$(1,1)$} \\
GLS / spatial mixed model, prior $h$, nugget $\sigma^2$ & \multicolumn{2}{l}{$\omega_j=\tilde\omega_j=1/(h(\lambda_j)+\sigma^2)$} \\
Eigenvector filtering, $K$ fixed effects & \multicolumn{2}{l}{$\omega_j=\tilde\omega_j=\ind\{j>K\}$} \\
Spatial+ \citep{Dupont2022} & \multicolumn{2}{l}{$\omega_j=(1-t_j)(1-s_j)$, $\tilde\omega_j=(1-t_j)(1-s_j)^2$} \\
gSEM (double residual) \citep{ThadenKneib2018} & \multicolumn{2}{l}{$\omega_j=(1-s_j)(1-t_j)$, $\tilde\omega_j=(1-s_j)^2$} \\
High-pass GLS (idealized spectral adjustment) \citep{Guan2023} & \multicolumn{2}{l}{$\omega_j=\tilde\omega_j=\ind\{j>K\}/(h(\lambda_j)+\sigma^2)$} \\
Regression on GNN features $p(L)X$ & \multicolumn{2}{l}{$\omega_j=p(\lambda_j)$, $\tilde\omega_j=p(\lambda_j)^2$} \\
\bottomrule
\end{tabularx}
\end{table}

\subsection{Confounding in the spectral domain}\label{sec:confounding}

Let $(X,U)$ be jointly Gaussian and jointly $L$-stationary, with coefficients $(x_j,u_j)$, $f_j=\Var(x_j)>0$, $c_j=\Cov(x_j,u_j)$ and $s^U_j=\Var(u_j)$. Define
\begin{equation}\label{eq:decomp}
\alpha_j=\frac{c_j}{f_j},\qquad v_j=u_j-\alpha_jx_j,\qquad g_j=\Var(v_j)=s^U_j-\frac{c_j^2}{f_j}>0 .
\end{equation}
Then $V=\sum_jv_j\varphi_j$ is independent of $X$ and $U=\alpha(L)X+V$, where $\alpha(L)$ has eigenvalues $\alpha_j$. We call $V$ the \emph{residual confounder}: the part of the confounder that the covariate does not predict. The response $Y=\beta X+U$ has coefficients
\begin{equation}\label{eq:model}
y_j=\gamma_jx_j+v_j,\qquad \gamma_j=\beta+\alpha_j ,
\end{equation}
so $\gamma_j$ is the association between covariate and response at frequency $\lambda_j$; $\alpha_j$ is the confounding at that frequency, the transfer function of \citet{Guan2023}. The class of \emph{vanishing confounding} is $\cA_0=\{\alpha:\alpha_j\to0\}$.

The Cameron--Martin space of $V$ is $\CM_V=\{h=\sum_jh_j\varphi_j:\norm{h}_{\CM_V}^2=\sum_jh_j^2/g_j<\infty\}$. Its role is that a deterministic shift $h$ of the law of $V$ produces an equivalent (mutually absolutely continuous) law if $h\in\CM_V$ and a mutually singular law otherwise \citep{Bogachev1998}. Set
\begin{equation}\label{eq:energy}
q_j=\frac{f_j}{g_j},\qquad \cI=\cI(X\Vert V)=\sum_{j\ge1}q_j=\E\norm{X}_{\CM_V}^2 .
\end{equation}
We call $\cI$ the \emph{Cameron--Martin energy} of the covariate relative to the residual confounder; it is also the expected Fisher information about $\beta$ carried by a complete record of $(X,Y)$. By \cref{lem:chisq}, $X\in\CM_V$ almost surely if $\cI<\infty$ and $X\notin\CM_V$ almost surely if $\cI=\infty$.

\subsection{Observation schemes and consistency}\label{sec:schemes}

We consider four ways in which data accumulate.
\begin{enumerate}[label=(O\arabic*),leftmargin=3em,itemsep=1pt]
\item \emph{Spectral infill.} At stage $N$ one observes $(x_j,y_j)_{j\le N}$. This idealizes the observation of the fields at increasing resolution in a fixed region, $N$ being the number of resolved modes.
\item \emph{Point-sample infill.} At stage $n$ one observes $X$ and $Y$ (the latter with a nugget) at points that fill $M$ (\cref{prop:pointsamples}).
\item \emph{Fixed resolution, increasing domain.} A stationary field on $\R^d$ is observed on the lattice $\delta\Z^d$, or through block means over cells of side $\delta$, in growing windows (\cref{sec:resolution}).
\item \emph{General channels.} Data are any randomized function of the covariate and of the latent response field (\cref{sec:beyond}).
\end{enumerate}

\begin{definition}[Consistent estimability]\label{def:consistent}
A functional $\psi(\theta)$ is \emph{consistently estimable} over a parameter set $\Theta$ under a scheme if there are estimators $T_n$, measurable functions of the stage-$n$ data, with $T_n\to\psi(\theta)$ in $P_\theta$-probability for every $\theta\in\Theta$.
\end{definition}

The basic obstruction is \cref{lem:obstruction}: if $P_\theta$ and $P_{\theta'}$ are mutually absolutely continuous on the $\sigma$-field generated by all stages and $\psi(\theta)\ne\psi(\theta')$, then $\psi$ is not consistently estimable over any $\Theta$ containing $\theta$ and $\theta'$.

\begin{proposition}[Spectral normal form]\label{prop:normalform}
\begin{enumerate}[label=(\alph*),itemsep=1pt]
\item Each model in the upper block of \cref{tab:zoo} is specified through a spectral function of the stated operator: a centred Gaussian prior with covariance operator $h(L)$ (CAR, Leroux, SAR, Whittle--Mat\'ern, graph Mat\'ern, random-effects ESF), an improper flat prior on a spectral subspace (Moran eigenvector filtering), spectral priors on coefficient functions ((M)GWR as a spectral SVC model) or a spectral filter (linear GNNs).
\item Under \eqref{eq:model}, for real weights $\omega_j$ and weights $\tilde\omega_j\ge0$, the estimator $\hat\beta=\sum_j\omega_jx_jy_j/\sum_j\tilde\omega_jx_j^2$ satisfies
\begin{equation}\label{eq:normalform}
\hat\beta-\beta=\sum_j\pi_j\alpha_j+\beta\Big(\sum_j\pi_j-1\Big)+\frac{\sum_j\omega_jx_jv_j}{\sum_j\tilde\omega_jx_j^2},\qquad \pi_j=\frac{\omega_jx_j^2}{\sum_k\tilde\omega_kx_k^2}.
\end{equation}
Each estimator in the lower block of \cref{tab:zoo} is of this form with the stated weights.
\item If $\omega=\tilde\omega\ge0$, then $\sum_j\pi_j=1$ and $\hat\beta=\sum_j\pi_j\gamma_j+\text{noise}$: every such estimator targets a weighted average of the scale-specific associations, with weights fixed by its spectral profile and by the covariate spectrum.
\end{enumerate}
\end{proposition}

\proofidea Each model in \cref{tab:zoo} is defined through a covariance, precision or smoothing operator that is a function of $L$, so it acts diagonally on the eigencoordinates $x_j,y_j$. An estimator built from such operators is therefore a ratio of sums over scales, and substituting $y_j=\gamma_jx_j+v_j$ gives \eqref{eq:normalform}. The only non-obvious rows are RSR, whose restricted random effect is orthogonal to $X$ and drops out of the estimating equation, and spatial+, whose penalized second stage can be profiled explicitly. The complete proof is in \cref{app:normalform}.

Spectral bias formulas of this kind are not new \citep{Paciorek2010,HughesHaran2013,Dupont2023,Woodward2024}; \cref{prop:normalform} is a unifying restatement for estimators built from functions of $L$. Part (c) is the lens through which the rest of the paper reads the debate: estimators differ only in which scales they weight, and confounding biases them exactly to the extent that it is present at the weighted scales. \Cref{tab:notation} collects the notation and its spatial reading.

\begin{table}[!tbp]
\centering\small
\caption{Notation and its spatial reading.}\label{tab:notation}
\begin{tabularx}{\textwidth}{>{\raggedright\arraybackslash}p{2.7cm}>{\raggedright\arraybackslash}p{5.0cm}X}
\toprule
Symbol & Meaning & Spatial reading \\
\midrule
$L$, $(\lambda_j,\varphi_j)$ & Laplace-type operator on the study region $M$ (domain, sphere, network or graph) and its eigenpairs & $\varphi_j$ is a spatial pattern whose scale is of order $\lambda_j^{-1/2}$; large $j$ means fine scale (Moran eigenvectors on graphs) \\
$d_s$ & spectral dimension, $\lambda_j\asymp j^{2/d_s}$ & $1$ on river or road networks, $2$ on surfaces and the sphere, $3$ in volumes \\
$x_j$, $y_j$ & coefficients of covariate and response on $\varphi_j$ & the covariate and the response seen at scale $j$ \\
$f_j$, $g_j$ & spectra of the covariate $X$ and of the residual confounder $V$ & how strongly the covariate, and the unexplained part of the confounder, vary at scale $j$ \\
$\nu_X$, $\nu_V$ & smoothness indices of Whittle--Mat\'ern spectra & small $\nu$: a rough surface with strong local variation \\
$\alpha_j$, $\gamma_j=\beta+\alpha_j$ & confounding and scale-specific association at scale $j$ & the slope carried by variation at scale $j$ alone \\
$\cA_0$ & confounding with $\alpha_j\to0$ & confounders act at broad scales; none at the finest scales \\
$\cI=\sum_jf_j/g_j$ & Cameron--Martin energy & total local signal of the covariate relative to residual spatial variation \\
$\delta$, $T_\delta$, $\pi/\delta$ & lattice spacing, resolved frequency band, Nyquist frequency & the grain (cell size) of the data and the finest resolvable scale \\
$w_m$, $\bar\alpha^\delta$, $v^\delta$ & alias weights, alias mean and alias variance of the confounding & how variation below the grain folds back onto the resolved scales \\
\bottomrule
\end{tabularx}
\end{table}
\section{The Cameron--Martin frontier}\label{sec:frontier}

\subsection{Main theorem}

\begin{theorem}[Cameron--Martin frontier]\label{thm:frontier}
Work under scheme (O1) with spectra $(f_j)$ and $(g_j)$ fixed, and let $\cI=\sum_jf_j/g_j$.
\begin{enumerate}[label=(\roman*),itemsep=2pt]
\item \textbf{Impossibility.} If $\cI<\infty$, then for every confounding sequence $\alpha$ and all $\beta\ne\beta'$ the laws $P_{\beta,\alpha}$ and $P_{\beta',\alpha}$ of the complete record are mutually absolutely continuous. Hence $\beta$ is not consistently estimable over any parameter set containing two such points, even when $\alpha=0$ is known. Quantitatively, for every estimator $T$ of the complete record and $\Delta=\abs{\beta_1-\beta_0}$,
\begin{equation}\label{eq:lecam}
\max_{i\in\{0,1\}}P_{\beta_i,\alpha}\big(\abs{T-\beta_i}\ge\Delta/2\big)\;\ge\;\tfrac14\exp\!\big(-\Delta^2\cI/4\big),
\end{equation}
and for every prior on $\beta$ with an absolutely continuous density $\pi$ of finite Fisher information $J(\pi)=\int(\pi')^2/\pi$, $\int\E_{\beta,\alpha}(T-\beta)^2\,\pi(\beta)\,d\beta\ge(\cI+J(\pi))^{-1}$.
\item \textbf{Possibility.} If $\cI=\infty$, let $h$ be any positive working spectrum with $0<\inf_jg_j/h_j\le\sup_jg_j/h_j<\infty$. The working GLS estimator
\begin{equation}\label{eq:gls}
\hat\beta_N=\frac{\sum_{j\le N}x_jy_j/h_j}{\sum_{j\le N}x_j^2/h_j}
\end{equation}
satisfies $\hat\beta_N\to\beta+\alpha_\infty$ almost surely whenever $\alpha_j\to\alpha_\infty\in\R$. In particular, $\beta$ is consistently estimable over $\R\times\cA_0$, by one and the same estimator for every nuisance spectrum $g$ in the tail-equivalence class of $h$.
\item \textbf{Dichotomy.} Let $\cA$ be any class of confounding sequences and $\cN_q=\{\delta:\sum_j\delta_j^2q_j<\infty\}$. If there exist $c\ne0$ and $\alpha,\alpha'\in\cA$ with $c+\alpha-\alpha'\in\cN_q$, then $\beta$ is not consistently estimable over $\R\times\cA$. If no such triple exists, any two parameters $(\beta,\alpha)$, $(\beta',\alpha')$ in $\R\times\cA$ with $\beta\ne\beta'$ induce mutually singular laws. For $\cA=\cA_0$ the first alternative holds if and only if $\cI<\infty$, and otherwise (ii) supplies a consistent estimator.
\end{enumerate}
\end{theorem}

\proofidea Under (O1) the record is a product over scales: independently across $j$, $x_j\sim N(0,f_j)$ and $y_j\mid x_j\sim N(\gamma_jx_j,g_j)$. Changing $\beta$ by $\Delta$ moves each conditional mean by $\Delta x_j$, and the Hellinger affinity of the two laws at scale $j$ is $(1+\Delta^2q_j/4)^{-1/2}$. By Kakutani's theorem on infinite products, the complete records are equivalent when these scale-wise discrepancies are summable, that is when $\cI<\infty$, and singular otherwise; equivalent laws cannot be told apart by any estimator (\cref{lem:obstruction}). The two quantitative bounds follow from the product formula for affinities and from a Bayesian Cram\'er--Rao (van Trees) bound. For (ii), GLS is a weighted average of the $\gamma_j$, with weights $x_j^2/h_j$ whose total diverges when $\cI=\infty$, plus a noise term whose conditional variance grows no faster than that total; the Toeplitz lemma and a Kolmogorov--Kronecker strong law dispose of the two parts. Part (iii) applies the same affinity computation to a general difference $\gamma'-\gamma$. The complete proof is in \cref{app:frontier}.

The condition $\cI<\infty$ says that the covariate lies in the Cameron--Martin space of the residual confounder: every covariate pattern the data will ever show is a plausible fluctuation of residual confounding, so a change in $\beta$ can be absorbed by a change in a quantity no data set can pin down. Three features deserve emphasis.

First, part (i) holds with no confounding at all: when the covariate is too smooth relative to residual spatial variation, the problem is lack of information, not bias. With confounding the two problems compound: the bias $\sum_j\pi_j\alpha_j$ of any normal-form estimator retains non-vanishing weight on the confounded low frequencies.

Second, part (i) holds at every parameter value, confounded or not, because it is formulated through the residual confounder $V$; the $\R^d$-stationary non-estimability result of \citet[Theorem~5]{DattaStein2025} is stated for independent covariate and confounder. For a deterministic regressor, (i)--(ii) reduce to the classical fact that a regression coefficient is estimable without error from a continuous record exactly when the regressor lies outside the reproducing kernel Hilbert space of the noise \citep{Grenander1950,Parzen1961}.

Third, the abstract criterion (iii) makes precise what any identifying assumption must accomplish: it must prevent a scale-invariant shift $c$ of the association from being mimicked, up to a Cameron--Martin-negligible error, by a difference of two admissible confounding patterns. Vanishing high-frequency confounding is a natural nonparametric assumption that achieves this when $\cI=\infty$, but not the weakest: vanishing along a fixed, known subsequence $(j_k)$ of modes with $\sum_kq_{j_k}=\infty$ also suffices (apply \cref{thm:frontier}(ii) to the modes $(j_k)$), and plurality-type assumptions \citep{Burman2025} play the same role. Every such assumption is useless when $\cI<\infty$. For classes other than $\cA_0$, the second alternative of (iii) gives only pairwise singularity, which does not by itself produce a consistent estimator \citep{BreimanLeCamSchwartz1964}. Conversely, any confounding class that contains a sequence $\alpha$ with $\alpha_j\to\alpha_\infty\ne0$ together with its recentred version $\alpha-\alpha_\infty$ leaves $\beta$ unidentified, since $(\beta,\alpha)$ and $(\beta+\alpha_\infty,\alpha-\alpha_\infty)$ induce the same law; this is the sense in which vanishing confounding is necessary.

\subsection{Smoothness, dimension and the irrelevance of range}

\begin{corollary}[The smoothness frontier on general supports]\label{cor:matern}
Let $L$ have spectral dimension $d_s$, and suppose $f_j\asymp(1+\lambda_j)^{-(\nu_X+d_s/2)}$ and $g_j\asymp(1+\lambda_j)^{-(\nu_V+d_s/2)}$ with $\nu_X,\nu_V>0$. Then $\cI=\infty$ if and only if
\begin{equation}\label{eq:threshold}
\nu_X\le\nu_V+\tfrac{d_s}{2}.
\end{equation}
Consequently, under (O1) and vanishing confounding, $\beta$ is consistently estimable if and only if \eqref{eq:threshold} holds; the verdict depends on no parameter other than the two smoothness indices (in particular not on ranges or variances). At the boundary $\nu_X=\nu_V+d_s/2$, GLS with a working spectrum tail-equivalent to $g$ is consistent, and the variance of its noise term is of order $1/\log N$.
\end{corollary}

\proofidea The spectral dimension turns the decay of $q_j=f_j/g_j$ in $\lambda_j$ into the power law $q_j\asymp j^{-2(\nu_X-\nu_V)/d_s}$, and the series $\sum_jq_j$ diverges exactly when the exponent is at least $-1$. Ranges and variances only change the constants. At the boundary $\sum_{j\le N}q_j\asymp\log N$, and the variance of GLS is the reciprocal of this sum. The complete proof is in \cref{app:matern}.

Thus the covariate may be smoother than the residual confounder, but by at most half the spectral dimension: by $1/2$ on a river or road network, by $1$ on a surface or the sphere, by $3/2$ in three dimensions. Two facts underlie the half-dimension. The sample paths of a Mat\'ern field of smoothness $\nu$ have Sobolev regularity just below $\nu$, whereas its Cameron--Martin space is the Sobolev space of order $\nu+d_s/2$ \citep{Driscoll1973,Steinwart2019,Kanagawa2018}; and Weyl's law converts spectral decay into summability in dimension $d_s$. For $\R^d$-stationary Gaussian fields observed on a bounded domain, the same threshold was obtained by \citet{DattaStein2025}. They settled the boundary $\nu_X=\nu_V+d/2$ for $d=1$ and integer $\nu_X$ (their Corollary~2) and conjectured the remaining cases. In the operator setting the boundary reduces to the divergence of $\sum_jj^{-1}$ and is settled by \cref{cor:matern}, which also extends the threshold to networks and the sphere; this does not resolve the $\R^d$-stationary conjecture.

\subsection{What spatial random-effects models estimate}

\begin{proposition}[The finest-scale association]\label{prop:whatgls}
\begin{enumerate}[label=(\alph*),itemsep=2pt]
\item Under (O1), if $\cI=\infty$ and $\alpha_j\to\alpha_\infty$, every GLS estimator with a working spectrum tail-equivalent to $g$ converges almost surely to $\beta+\alpha_\infty=\lim_j\gamma_j$, the association at the finest scale.
\item \emph{(A three-regime phase diagram.)} Let $X=Z+H$ and $U=cZ$ with $c\ne0$ and $Z,H$ independent, Gaussian and $L$-stationary, with spectra $f^Z_j,f^H_j>0$ for all $j$, $f^Z_j\sim m_Z\lambda_j^{-a_Z}$ and $f^H_j\sim m_H\lambda_j^{-a_H}$. Then $\cI=\infty$ and
\begin{equation}\label{eq:threeregimes}
\alpha_\infty=\begin{cases}0,& a_Z>a_H\ \ (\text{confounder smoother than the covariate's own variation}),\\[2pt] c\,\dfrac{m_Z}{m_Z+m_H},& a_Z=a_H,\\[6pt] c,& a_Z<a_H .\end{cases}
\end{equation}
For Whittle--Mat\'ern spectra $\tau^2(\kappa^2+\lambda)^{-a}$ the constants are $m=\tau^2$; in the classical Mat\'ern parametrization with marginal variance $\sigma^2$, $m\propto\sigma^2\kappa^{2\nu}$, the microergodic parameter \citep{Zhang2004}.
\item \emph{(Self-confounding.)} If the response contains $\beta SX$ for a smoothing operator $S=s(L)$ with $s_j\to0$, instead of $\beta X$, and the confounder is independent of $X$, then $\alpha_j=\beta(s_j-1)\to-\beta$, and GLS with a working spectrum tail-equivalent to $g$ converges to $0$ whenever $\cI=\infty$.
\end{enumerate}
\end{proposition}

\proofidea Part (a) is \cref{thm:frontier}(ii) read for a non-zero limit $\alpha_\infty$. In (b) all spectral quantities are explicit: $\alpha_j=cf^Z_j/(f^Z_j+f^H_j)$, so $\alpha_\infty$ is decided by the ratio of two power laws at high frequency. In (c) the smoothing operator itself creates the confounding $\alpha_j=\beta(s_j-1)$. The complete proof is in \cref{app:whatgls}.

Part (b) recovers and sharpens the scale story of \citet{Paciorek2010}. When the confounder and the covariate's unconfounded variation have the same smoothness, range matters only through the microergodic constants: a shorter range of the unconfounded component (larger $\kappa_H$, hence larger $m_H$) reduces the asymptotic bias of GLS, and a longer one increases it, possibly beyond that of OLS (\cref{fig:glslimit}a). When smoothness differs, the smoothness ordering alone decides between no bias and full bias, whatever the ranges. Part (c) is the self-confounding result of \citet[Theorem~1 and Proposition~1]{BolinWallin2026}, recovered as the special case $\alpha_j=\beta(s_j-1)$ of part (a): for covariates rougher than the signal they drive, the finest-scale association is zero. \citet{Bolin2026} removes self-confounding by matching the regression operator to the smoothing. \Cref{fig:glslimit} confirms \eqref{eq:threeregimes} numerically.

\subsection{Point samples with a nugget}

Scheme (O1) is an idealization. The next result shows that both halves of \cref{thm:frontier} survive the passage to point samples observed with measurement error.

\begin{proposition}[Point samples with a nugget]\label{prop:pointsamples}
Let $\mu(M)<\infty$ and $\operatorname{supp}\mu=M$, let $X$ and $U$ have continuous sample paths, and let the response be observed as $Y_i=Y(s_i)+\varepsilon_i$ with independent and identically distributed errors of mean zero and finite variance (of any distribution), independent of $(X,U)$.
\begin{enumerate}[label=(\alph*),itemsep=2pt]
\item If $\cI<\infty$, then for every design $(s_i)_{i\ge1}$, dense or not, the laws of the complete data sequence under $(\beta,\alpha)$ and $(\beta',\alpha)$ with $\beta\ne\beta'$ are mutually absolutely continuous; no consistent estimator exists.
\item Suppose $\cI=\infty$, $\alpha\in\cA_0$ and $h$ is tail-equivalent to $g$ as in \cref{thm:frontier}(ii); the eigenfunctions are Lipschitz with constants $\Lambda_j$ and sup-norms $B_j$; the paths of $X$ and $U$ are H\"older continuous of some unknown order $\vartheta\in(0,1]$; and $\{A^{(n)}_i\}_{i\le n}$ are partitions of $M$ with $s_i\in A_i$, mesh $\ell_n\to0$ and $w_i=\mu(A_i)$ with $\max_iw_i\to0$. Put $\hat x_j=\sum_iw_iX(s_i)\varphi_j(s_i)$, $\hat y_j=\sum_iw_iY_i\varphi_j(s_i)$, $a_J=\sum_{j\le J}(B_j+\Lambda_j)/h_j$, $\varpi_n=\max(\ell_n,\max_iw_i)$ and $J_n=\max\{J:a_J\le\log(1/\varpi_n)\}$. Then the sieve estimator $\tilde\beta_n=\sum_{j\le J_n}\hat x_j\hat y_j/h_j\big/\sum_{j\le J_n}\hat x_j^2/h_j$ is consistent.
\end{enumerate}
\end{proposition}

\proofidea (a) The point data are a measurable function of the complete fields plus independent noise, so the equivalence of \cref{thm:frontier}(i) passes to them. (b) Quadrature turns point data into estimates $\hat x_j,\hat y_j$ of the first $J$ coefficients, with errors controlled by the mesh, the H\"older regularity of the paths and the averaged nugget. The truncation $J_n$ grows slowly enough that the accumulated error, weighted by $1/h_j$, vanishes, and the sieve estimator is then a small perturbation of the complete-record GLS estimator of \cref{thm:frontier}(ii). The complete proof is in \cref{app:pointsamples}.

The sieve truncation $J_n$ depends on the design, on the eigenfunctions and on the working spectrum $h$, which must be tail-equivalent to $g$, so the tail class of the residual spectrum must be known. It grows only polylogarithmically and is deliberately conservative. Whether GLS with a nugget, which in effect truncates softly at a much higher frequency, is consistent under (O2) is not covered by this result (see experiment E9).
\section{The resolution theorem}\label{sec:resolution}

Scheme (O1) lets the number of resolved modes grow. Most socioeconomic, health and environmental data do not: census tracts, counties, satellite pixels and regular grids have a fixed size, and studies grow by covering more territory. We now show that this difference is decisive.

Let $(X,U)$ be jointly stationary and Gaussian on $\R^d$, written as $U=\alpha\ast X+V$, where $\alpha:\R^d\to\mathbb C$ is a Hermitian transfer function ($\alpha(-\omega)=\overline{\alpha(\omega)}$), $V$ is a stationary Gaussian field independent of $X$ with spectral density $g$, and $X$ has spectral density $f_X$. The response is $Y=\beta X+U$, so that $\gamma(\omega)=\beta+\alpha(\omega)$ is the association at frequency $\omega$. At resolution $\delta>0$ we observe either the \emph{point samples} $D_\delta=(X(\delta k),\,Y(\delta k)+\varepsilon_k)_{k\in\Z^d}$ or the \emph{block means} $B_\delta=(\bar X_k,\,\bar Y_k+\varepsilon_k)_{k\in\Z^d}$ over the cells $Q_k=\delta(k+[0,1)^d)$, with optional nugget errors $\varepsilon_k$, independent and $N(0,\sigma^2)$. Observing the whole lattice is the most favourable increasing-domain limit, so non-identification from $D_\delta$ or $B_\delta$ implies non-identification from every finite window.

\begin{assumption}\label{ass:resolution}
(R1) $f_X$ is continuous, even and strictly positive; $g$ is even, measurable and locally bounded away from zero and infinity; and the alias sums $\sum_{m\in\Z^d}f(\omega+2\pi m/\delta)$ of $f=f_X$ and $f=g$ converge uniformly on compact sets (essentially uniformly for $g$), for instance when $f(\omega)\le C(1+\abs{\omega})^{-d-\epsilon}$. (R2) $\alpha$ is bounded, measurable and Hermitian, with $\operatorname{ess\,sup}_{\abs{\omega}>R}\abs{\alpha(\omega)}\to0$ as $R\to\infty$ (vanishing confounding in the continuum).
\end{assumption}

Write $T_\delta=[-\pi/\delta,\pi/\delta)^d$ and $\omega_m=\omega+2\pi m/\delta$ for the aliases of $\omega\in T_\delta$, and let
\begin{equation}\label{eq:aliasweights}
w_m(\omega)=\frac{f_X(\omega_m)}{\sum_kf_X(\omega_k)},\qquad\bar\alpha^\delta(\omega)=\sum_mw_m(\omega)\,\alpha(\omega_m),\qquad v^\delta(\omega)=\sum_mw_m(\omega)\,\abs{\alpha(\omega_m)-\bar\alpha^\delta(\omega)}^2
\end{equation}
be the alias weights and the alias mean and alias variance of the confounding (for block means replace $f_X$ by $f_X\abs{\hat\psi_\delta}^2$, where $\hat\psi_\delta(\omega)=\prod_i\operatorname{sinc}(\delta\omega_i/2)$).

\begin{theorem}[Resolution theorem]\label{thm:resolution}
Fix $\delta>0$, let \cref{ass:resolution} hold, and let $\beta'\in\R$ be arbitrary.
\begin{enumerate}[label=(\alph*),itemsep=1pt]
\item \textbf{Unknown residual spectrum.} There exist $(\alpha',g')$ satisfying \cref{ass:resolution}, with $\alpha'-\alpha$ and $g'-g$ bounded and compactly supported, such that the law of $D_\delta$, respectively of $B_\delta$, is the same under $(\beta',\alpha',f_X,g')$ as under $(\beta,\alpha,f_X,g)$. If $\alpha$ and $g$ are continuous, so are $\alpha'$ and $g'$.
\item \textbf{Known $g$, point samples.} Let $g$ be known. If $\operatorname{ess\,inf}_{\omega\in T_\delta}v^\delta(\omega)>0$, the conclusion of (a) holds for $D_\delta$ with $g'=g$ and the same nugget variance. If the data carry a nugget of positive, unknown variance $\sigma^2$, it holds with $g'=g$ and any nugget variance $\sigma'^2\in[0,\sigma^2)$ in place of $\sigma^2$. If instead the nugget variance is known (or there is no nugget) and $\operatorname{ess\,inf}v^\delta=0$, then $\beta$ is identified from $D_\delta$: no $(\beta',\alpha',f_X,g)$ with $\beta'\ne\beta$ and $\alpha'$ satisfying (R2) induces the same law.
\item \textbf{Known $g$, block means.} If $g$ is known and $\alpha$ vanishes on all aliases of a set of frequencies of positive measure (for instance $\alpha\equiv0$), then $\beta$ is identified from $B_\delta$, whether or not the nugget variance is known.
\end{enumerate}
Consequently, at a fixed lattice resolution with an unknown residual spectrum, nothing about $\beta$, not even its sign, is identified, however large the observation window and with or without a nugget.
\end{theorem}

\proofidea The proof is a construction. Lattice data are stationary on $\Z^d$, and their law is determined by the aliased spectral matrix. Its cross term is $F^\delta_{XY}=F^\delta_X(\beta+\bar\alpha^\delta)$, and its response term is $F^\delta_Y=F^\delta_X(\abs{\beta+\bar\alpha^\delta}^2+v^\delta)+\sum_mg(\omega_m)$, where $F^\delta_X=\sum_mf_X(\omega_m)$. To move $\beta$ to $\beta'=\beta+\Delta$, we replace the confounding below a cutoff frequency $R$ by its alias average shifted by $-\Delta$, with a correction that keeps the aliased cross term exact. The cutoff can be as high as we like, so the leftover alias variance can be made arbitrarily small. Under (a) the residual spectrum absorbs it; under (b) a zero-mean perturbation restores the original alias variance exactly. The alternative confounding still vanishes beyond the cutoff: both models satisfy the identifying assumption, and the data cannot tell them apart. Parts (b) and (c) show that the only escape is second-moment information. If $g$ and the nugget variance were known exactly, the response spectrum would reveal the alias average of $\abs{\gamma}^2$. The Cauchy--Schwarz inequality then pins down $\beta$ wherever the alias variance of the confounding comes arbitrarily close to zero. For point samples this happens if and only if $\operatorname{ess\,inf}v^\delta=0$; an unknown nugget variance already destroys it. This identification rests entirely on exact knowledge of second-order nuisance quantities and has no counterpart in practice. The complete proof, in twenty-one steps, is in \cref{app:resolution}.

\begin{corollary}[What lattice data identify]\label{cor:latticetarget}
Under \cref{ass:resolution}, the law of lattice data depends on $(\beta,\alpha,g)$ only through the alias-averaged association
\begin{equation}\label{eq:aliasavg}
\gamma^\delta(\omega)=\sum_{m\in\Z^d}w_m(\omega)\,\gamma(\omega_m)=\beta+\bar\alpha^\delta(\omega),\qquad\omega\in T_\delta,
\end{equation}
and the aliased response spectrum $F^\delta_X(\abs{\gamma^\delta}^2+v^\delta)+\sum_mg(\omega_m)$. In particular, an analysis that assumes the coherence to vanish at the highest resolved frequency $\omega^\ast$, as the semiparametric CAR model of \citet{Guan2023} does at the largest graph eigenvalue, targets
\begin{equation}\label{eq:latticebias}
\beta+\sum_mw_m(\omega^\ast)\,\alpha(\omega^\ast_m)
\end{equation}
(for continuous $\alpha$), which equals $\beta$ if and only if the alias-weighted confounding at $\omega^\ast$ vanishes. If $\alpha\ge0$ everywhere and $\alpha(\omega^\ast_m)>0$ for some alias, the target exceeds $\beta$.
\end{corollary}

\proofidea The first statement is the lattice-law computation that opens the proof of \cref{thm:resolution}; the second evaluates the alias mean at $\omega^\ast$, where all alias weights are positive (\cref{app:resolution}, Steps~4 and~21).

Three consequences are worth stating plainly. First, the assumption that confounding vanishes at high frequencies has \emph{no identifying power} at fixed resolution. Any identification of $\beta$ from gridded data must come from assumptions that tie $\alpha$ below the Nyquist frequency to $\alpha$ above it, such as a parametric coherence function \citep{SchnellPapadogeorgou2020,Guan2023,Tang2026}, monotonicity or analyticity, or that bound the frequency beyond which confounding is known to vanish. Such assumptions extrapolate beyond the resolution of the data. Second, the graph-native version of the assumption, an unconfounded highest graph frequency, is not implied by its continuum version. By \eqref{eq:latticebias} it holds if and only if the alias-weighted confounding at $\omega^\ast=(\pi/\delta,\dots,\pi/\delta)$ vanishes, essentially when confounding is absent at the scale of the cells. Because the alias weights concentrate near $\omega^\ast$, the bias in \eqref{eq:latticebias} is of the order of the confounding at the resolution scale: in the configuration of \cref{fig:resolution}c it is about $0.8\,\alpha(\pi/\delta)$ for point samples and $0.99\,\alpha(\pi/\delta)$ for block means. Third, the resolution floor is quantitative. \Cref{thm:minimax} shows that with $N$ resolved modes the squared bias floor is of the order of the squared worst-case confounding at the finest resolved scale.

How far apart can observationally equivalent effects be? By \cref{thm:resolution}(a) arbitrarily far: the identified set is the whole real line. The price of a distant alternative is the frequency up to which its confounding must differ from $\alpha$. The leftover alias variance is controlled by the alias mass beyond the cutoff $R$, so a shift $\Delta$ needs $R$ large enough that this mass is small compared with $\Delta^{-2}$. Experiment E4 (\cref{sec:experiments}) constructs models with $\beta=1$ and with $\beta'=-1$ that induce the same lattice law. For point samples it does so with the residual spectrum held fixed (part b); for point samples and block means it also does so with an adjusted residual spectrum (part a).

The theorem is proved for the cubic lattice $\delta\Z^d$ with point values or cubic block means, where stationarity makes the aliasing structure explicit. Irregular partitions such as counties are not stationary, and we do not claim the theorem for them. The mechanism, however, does not depend on regularity: a fixed partition supplies a bounded number of observations per unit area, and such data cannot constrain the behaviour of $\alpha$ at arbitrarily high frequencies. We expect, but do not prove, the same conclusion for irregular areal units.
\section{Which estimators reach the frontier}\label{sec:estimators}

\Cref{thm:frontier} asserts that \emph{some} estimator is consistent inside the frontier. We now ask which of the estimators in \cref{tab:zoo} are. Throughout this section we work under (O1), with fixed working spectra and no working nugget, in the power-law regime
\begin{equation}\label{eq:powerlaw}
f_j\asymp j^{-a},\qquad g_j\asymp j^{-b},\qquad h_j\asymp j^{-\varkappa},\qquad a>1,
\end{equation}
where $h$ is a working spectrum. For Whittle--Mat\'ern spectra on a support of spectral dimension $d_s$, $a=1+2\nu_X/d_s$, $b=1+2\nu_V/d_s$ and $\varkappa=1+2\nu_h/d_s$, with $\nu_h$ the working smoothness, or more generally the working spectral exponent (a valid working covariance requires $\nu_h>0$). A working nugget changes the picture. Under (O1) the weights $1/(h_j+\sigma^2)$ of \cref{tab:zoo} are bounded, so $\sum_jx_j^2/(h_j+\sigma^2)<\infty$, and GLS with any $\sigma^2>0$ converges to a non-degenerate random limit and is not consistent, even when $\alpha=0$ (\cref{app:estimators}, Step~7). A nugget is natural only under (O2), where it describes measurement error; that case is not covered by our theorems and is examined numerically in experiment E9.

\begin{theorem}[The GLS consistency window]\label{thm:window}
Under \eqref{eq:powerlaw}, let $u=\varkappa-a$ and $\Delta'=a-b$, and consider the working GLS estimator \eqref{eq:gls} with weights $1/h_j$.
\begin{enumerate}[label=(\roman*),itemsep=1pt]
\item If $u<-1$, then $\hat\beta_N$ is not consistent: even when $\alpha=0$ it does not converge in probability to $\beta$.
\item If $u>-1$, then $\hat\beta_N\to\beta$ in probability for every $\alpha\in\cA_0$ if and only if $\Delta'<1$.
\item If $u=-1$, then $\hat\beta_N\to\beta$ in probability for every $\alpha\in\cA_0$ if and only if $\Delta'\le1$.
\end{enumerate}
In Whittle--Mat\'ern terms, GLS with working spectral exponent $\nu_h$, any working range and any working variance, is consistent if and only if either $\nu_h>\nu_X-d_s/2$ and $\nu_X<\nu_V+d_s/2$, or $\nu_h=\nu_X-d_s/2$ and $\nu_X\le\nu_V+d_s/2$. If instead $\alpha_j\to\alpha_\infty\in\R$, the same conditions give $\hat\beta_N\to\beta+\alpha_\infty$ in probability, the finest-scale association.
\end{theorem}

\proofidea GLS is again signal plus noise. Under the power laws the signal weights $x_j^2/h_j$ behave like $j^u$; their total diverges if and only if $u\ge-1$, and only then does the signal forget the confounded low frequencies. Given $X$, the noise is Gaussian with variance $V_N/D_N^2$, where $D_N\asymp\sum_{j\le N}j^u$ and $V_N\asymp\sum_{j\le N}j^{2u+\Delta'}$, and comparing these power sums gives the two conditions. The complete proof is in \cref{app:estimators}.

The window is asymmetric. Under (O1) without a nugget, oversmoothing the working model, that is, choosing $\nu_h$ larger than the truth, does not affect consistency inside the open frontier: it upweights high frequencies, which is what deconfounding needs. Undersmoothing below $\nu_X-d_s/2$ is fatal: the weights then decay too fast for the high frequencies to dominate, and the estimator retains a non-vanishing share of the confounded low frequencies. This locates the exact smoothness boundary behind the GLS inconsistency that \citet[Supplement S1.2]{DattaStein2025} exhibit for an exponential working model. In the operator setting, undersmoothing is fatal if and only if $\nu_h<\nu_X-d_s/2$, so mildly undersmoothed working models remain consistent. Whether the same boundary holds for $\R^d$-stationary fields, where they conjecture inconsistency for other undersmoothed working classes, remains open. The theorem also explains their observation that GLS can fail under arbitrarily small misspecification. When $\nu_V$ is close to $\nu_X-d_s/2$, that is near the frontier, a working smoothness slightly below the true $\nu_V$ already leaves the window. At the frontier itself only a correctly specified tail works, for GLS or for high-pass GLS (\cref{prop:estimators}(d)). The consistency verdict does not depend on working ranges or variances, but finite-sample accuracy does. With a nugget under (O2), convergence rates depend on them as well. An oversmoothed working model then treats high frequencies as noise and effectively uses only the modes where its spectrum exceeds the nugget. Their number grows like a small power of $n$, so oversmoothing can make GLS nearly as poor as an undersmoothed fit at realistic sample sizes (experiment E9).

\begin{proposition}[Other estimators]\label{prop:estimators}
Under \eqref{eq:powerlaw} and $\alpha\in\cA_0$:
\begin{enumerate}[label=(\alph*),itemsep=1pt]
\item OLS and restricted spatial regression converge almost surely to the non-degenerate limit $\beta+\sum_j\alpha_jx_j^2/\sum_jx_j^2+\sum_jx_jv_j/\sum_jx_j^2$ and are therefore inconsistent for every $\alpha$, including $\alpha=0$; conditionally on $X$ the limit is biased unless $\sum_j\alpha_jx_j^2=0$.
\item Eigenvector filtering with a fixed number $K$ of eigenvectors converges to the analogous non-degenerate limit over the modes $j>K$; it is inconsistent for every $\alpha$, and biased conditionally on $X$ unless $\sum_{j>K}\alpha_jx_j^2=0$.
\item Eigenvector filtering with $K_N\to\infty$ and $K_N/N\to0$ is consistent if and only if $\Delta'<1$; the variance of its noise term, conditionally on $X$ in probability and unconditionally, is of order $K_N^{\Delta'-1}$.
\item High-pass GLS with weights $\ind\{K_N<j\le N\}/h_j$, $h$ tail-equivalent to $g$ and $K_N\to\infty$, is consistent if $\sum_{K_N<j\le N}q_j\to\infty$ and inconsistent if this sum stays bounded; in particular it is consistent whenever $\cI=\infty$ and $K_N\to\infty$ slowly enough.
\item Spatial+ whose first-stage smoother satisfies $s_j\to0$ and whose second stage is GLS with working spectrum $h$ obeys the window of \cref{thm:window}. The spatial+ of \cref{tab:zoo}, whose second stage is a penalized smooth, and gSEM with fixed smoothers have bounded weights; like OLS, they converge to a non-degenerate limit under (O1).
\end{enumerate}
\end{proposition}

\proofidea Each estimator is a normal-form ratio with explicit weights, and the signal-plus-noise comparison used for \cref{thm:window} applies. Bounded total weight (OLS, RSR, fixed $K$, penalized spatial+, gSEM) leaves a non-degenerate random limit. Discarding the low frequencies gives consistency when the noise variance of the retained modes vanishes: for unweighted filtering, part (c), this variance is of order $K_N^{\Delta'-1}$; for high-pass GLS with the correct tail, part (d), it is the reciprocal of the Cameron--Martin energy of the retained modes. The complete proof is in \cref{app:estimators}.

The same computation settles the areal models of \cref{tab:zoo}, once one asks how their working precisions behave as the graph is refined.

\begin{corollary}[Areal models under refinement]\label{cor:areal}
Work under (O1) and \eqref{eq:powerlaw}, and read a sequence of graphs $G_n$ with $n$ nodes as in \cref{sec:support}, identifying the $j$th graph mode with the $j$th mode of $L$ and ignoring aliasing: at stage $n$ the data are $(x_j,y_j)_{j\le n}$, and GLS uses the weights $w^{(n)}_j=p_n(\lambda^{(n)}_j)$, where $\lambda^{(n)}_1\le\dots\le\lambda^{(n)}_n$ are the eigenvalues of the graph Laplacian $L_n=D-W$ of $G_n$ and $p_n$ is the working precision as a function of the eigenvalue. The sums run over $2\le j\le n$; the constant eigenvector is absorbed by an intercept.
\begin{enumerate}[label=(\alph*),itemsep=1pt]
\item If $p_-\le w^{(n)}_j\le p_+$ for all sufficiently large $n$ and $2\le j\le n$, with constants $0<p_-\le p_+<\infty$, then GLS is not consistent, even when $\alpha=0$. If moreover $p_n=p$ does not depend on $n$, $p$ is continuous at $0$ with $p(0)>0$, and $\lambda^{(n)}_j\to0$ for each fixed $j$, then GLS converges almost surely to $\sum_{j\ge2}x_jy_j/\sum_{j\ge2}x_j^2$, the limit of OLS with an intercept. The bounds hold, with hyperparameters that do not depend on $n$, for the Leroux model with $0\le\rho<1$ and for graph Mat\'ern models on graphs of bounded degree, and for the proper CAR model with $\abs\rho<1$ and the SAR error model with $\abs\rho k<1$ on $k$-regular graphs; these precisions are continuous with $p(0)>0$, so they converge to the limit of OLS whenever the eigenvalues of fixed modes tend to zero.
\item If $c\,(j/n)^{2m/d_s}\le w^{(n)}_j\le C\,(j/n)^{2m/d_s}$ for all sufficiently large $n$ and $2\le j\le n$, for some $m>0$ and constants $0<c\le C$, then (i)--(iii) of \cref{thm:window} hold with $u=2m/d_s-a$. In Whittle--Mat\'ern terms GLS then has working smoothness $\nu_h=m-d_s/2$: for every $\alpha\in\cA_0$ it is consistent if and only if either $\nu_X<m$ and $\nu_X<\nu_V+d_s/2$, or $\nu_X=m$ and $\nu_X\le\nu_V+d_s/2$, and if $\nu_X>m$ it is not consistent even when $\alpha=0$. When $c\,(j/n)^{2/d_s}\le\lambda^{(n)}_j\le C\,(j/n)^{2/d_s}$ for $2\le j\le n$, as on cycles ($d_s=1$) and on discrete tori ($d_s=d$), the hypothesis holds with $m=1$ for the intrinsic CAR model, $p(\lambda)=\tau\lambda$, and for the BYM model with fixed variances, $p(\lambda)=\lambda/(\sigma_u^2+\sigma_v^2\lambda)$; with $m=2$ for second-order intrinsic models, $p(\lambda)=\tau\lambda^2$; with $m=1$ for the Leroux model and, on $k$-regular graphs, the proper CAR model when $0\le1-\rho_n\le C'n^{-2/d_s}$; and with $m=2$ for the SAR error model on $k$-regular graphs when $0\le1-\rho_nk\le C'n^{-2/d_s}$.
\item The intrinsic CAR estimator is the least-squares regression of neighbour differences: for symmetric weights $W_{ii'}$, $\inner X{L_nY}=\sum_{i<i'}W_{ii'}(X_i-X_{i'})(Y_i-Y_{i'})$, so that $\hat\beta=\sum_{i<i'}W_{ii'}(X_i-X_{i'})(Y_i-Y_{i'})\big/\sum_{i<i'}W_{ii'}(X_i-X_{i'})^2$.
\end{enumerate}
\end{corollary}

\proofidea (a) With weights bounded above and below, the denominator of GLS stays bounded while, given $X$, the variance of its noise term stays bounded away from zero, exactly as for OLS; when the weights of fixed modes converge to $p(0)$, dominated convergence gives the limit of OLS. (b) GLS does not change when all weights are multiplied by $n^{2m/d_s}$; after this rescaling they lie between multiples of $j^{2m/d_s}$, uniformly in $n$. The orders of the sums in the proof of \cref{thm:window} depend on the weights only through such bounds, and the steps of that proof that use weights fixed in $N$ are redone in probability. (c) is the summation-by-parts identity for the graph Laplacian. The complete proof is in \cref{app:areal}.

The corollary sorts the areal models of disease mapping and spatial econometrics by what their precision does at the lowest graph frequencies. A precision that stays bounded away from zero there, as in the Leroux, proper CAR and SAR models with fixed $\rho$ and in graph Mat\'ern models with fixed range, downweights the confounded broad scales by at most a bounded factor, so on refining graphs of bounded degree (regular graphs for proper CAR and SAR) these models are not consistent, even without confounding, and when the eigenvalues of fixed modes tend to zero they converge to the limit of OLS. A precision that vanishes at frequency zero, as in the intrinsic CAR and BYM models, downweights the broad scales without limit and, on graphs with Weyl-type eigenvalue bounds such as regular grids, places the estimator in the GLS window with working smoothness $1-d_s/2$. The intrinsic CAR estimator is then an all-edge version of the spatial first differences of \citet{DruckenmillerHsiang2018}, and it is consistent only for covariates rough enough to have infinite Dirichlet energy ($\nu_X\le1$). Hyperparameters that move towards the boundary at the rate $n^{-2/d_s}$, as when CAR and SAR precisions are read as discretizations of the continuum operator $\kappa^2-\Delta$ \citep{Lindgren2011,BesagMondal2005}, give the windows of $\nu_h=1-d_s/2$ (CAR) and $\nu_h=2-d_s/2$ (SAR). Two features of (b) differ from the setting of \cref{thm:window}. The BYM model adds a fixed variance $\sigma_v^2$ to the intrinsic CAR model, much like a nugget, yet obeys the window: as the graph refines, the intrinsic CAR variance of a fixed mode grows like $n^{2/d_s}$, so $\sigma_v^2$ acts as a vanishing nugget, unlike the fixed nugget of (O1). And the working smoothness $1-d_s/2$ of the intrinsic CAR model is not positive when $d_s\ge2$; such improper working spectra are allowed, because \cref{thm:window} does not require $\varkappa>1$. Experiment E10 confirms these predictions for the intrinsic CAR, BYM, Leroux and SAR models on cycles. In practice $\rho$ and the variances are estimated, which our results do not cover.

Regressions on features produced by a low-pass graph filter, such as those computed by linear graph neural networks, target a different quantity: by \eqref{eq:normalform} their weights vanish at high frequency, so they concentrate exactly where confounding lives. Nonlinear graph neural networks fall outside the normal form, but they remain subject to the impossibility statements of \cref{thm:frontier,thm:resolution,thm:shift}, which constrain every measurable estimator.

\subsection{Minimax rates and the resolution floor}

\begin{theorem}[Minimax rate]\label{thm:minimax}
Work under (O1) with $f$ and $g$ known and $q_j=f_j/g_j\asymp j^{p-1}$ for some $p>0$. For $r,C_\alpha>0$ let $\Theta(r,C_\alpha)=\{(\beta,\alpha):\beta\in\R,\ \abs{\alpha_j}\le C_\alpha j^{-r}\text{ for all }j\}$. There are constants $0<c\le C$ such that for all $N\ge N_0$
\begin{equation}\label{eq:minimax}
c\,\big(N^{-2r}+N^{-p}\big)\;\le\;\inf_{T_N}\sup_{\Theta(r,C_\alpha)}\E\,(T_N-\beta)^2\;\le\;C\,\big(N^{-2r}+N^{-p}\big),
\end{equation}
where the infimum runs over estimators based on $(x_j,y_j)_{j\le N}$. The oracle GLS estimator attains the upper bound. The lower bound $N^{-2r}$ holds for every $N$ and is witnessed by two parameters whose laws on the first $N$ modes coincide.
\end{theorem}

\proofidea The upper bound computes the mean squared error of oracle GLS: the bias is a weighted average of the $\alpha_j$, controlled by their decay, and the variance is an inverse moment of a weighted $\chi^2$ sum. Both lower bounds are two-point arguments. In the first, two parameters share the same association on every resolved mode, so their data laws coincide; in the second, the parameters are as far apart as the Cameron--Martin energy of the resolved modes allows while their Hellinger affinity stays bounded below. The complete proof is in \cref{app:minimax}.

The two terms have transparent meanings. The variance term $N^{-p}$ is the reciprocal of the Cameron--Martin energy of the resolved modes, and $p=1+2(\nu_V-\nu_X)/d_s$ in Whittle--Mat\'ern terms, so it degrades continuously as the roughness gap approaches the frontier. The bias term $N^{-2r}$ is the squared worst-case confounding at the finest resolved scale; if $\alpha(\lambda)\asymp\lambda^{-\varrho}$, then $r=2\varrho/d_s$. The least-favourable pair in the proof has no confounding at unresolved scales and constant confounding $-\Delta_N$ on the resolved ones: confounding of the size allowed at the finest resolved scale cannot be distinguished from a shift in $\beta$. On a sequence of graphs with $n$ nodes, the same bounds hold with $N=n$ when $q_j\asymp j^{p-1}$ holds uniformly in $n$, and the floor is then set by the confounding at the scale of the units. This is the quantitative face of \cref{thm:resolution}. \Cref{fig:frontier}c confirms the rate $N^{-\min(r,p/2)}$ for the root mean squared error.
\section{Beyond the linear--Gaussian model}\label{sec:beyond}

The frontier of \cref{sec:frontier} is not an artefact of linearity or Gaussianity. Its engine is a single principle about shifts.

\subsection{The shift-equivalence principle}

Let the covariate $X$ be a random element of a Polish space $\cX$ of functions on $M$, with arbitrary law. Let the residual confounder $V$ be independent of $X$ and take values in a separable Fr\'echet space $E$ of functions (or coefficient sequences) on $M$, with law $P_V$. The latent response is
\begin{equation}\label{eq:latent}
\eta=\Psi(X)+V,\qquad \Psi=\mathsf M+\mathsf A,
\end{equation}
where $\mathsf M:\cX\to E$ is an effect map ($x\mapsto\beta x$, or a nonlinear response function $x\mapsto m\circ x$, a dose--response curve in epidemiology) and $\mathsf A:\cX\to E$ is a confounding mechanism, any measurable and possibly nonlinear and non-local map; in the Gaussian case of \cref{sec:confounding}, $\mathsf A(x)=\alpha(L)x=\E[U\mid X=x]$. The \emph{admissible shifts} of $V$ are $\cQ(P_V)=\{h\in E:P_V(\cdot-h)\text{ and }P_V\text{ are mutually absolutely continuous}\}$, and $P_V$ has the \emph{dichotomy property} if every $h\notin\cQ(P_V)$ produces a mutually singular law.

\begin{theorem}[Shift equivalence]\label{thm:shift}
Let $\Psi_1,\Psi_2:\cX\to E$ be measurable and let $P_i$ be the law of $(X,\eta)$ under $\Psi_i$.
\begin{enumerate}[label=(\alph*),itemsep=1pt]
\item If $\Psi_1(X)-\Psi_2(X)\in\cQ(P_V)$ almost surely, then $P_1$ and $P_2$ are mutually absolutely continuous, and so are the laws of any data obtained from $(X,\eta)$ through a common Markov kernel (sampling design, aggregation, nugget, counts, measurement error). Any functional that differs between the two models is then not consistently estimable under any observation scheme.
\item If $P_V$ has the dichotomy property, then $\TV(P_1,P_2)\ge\Prob\big(\Psi_1(X)-\Psi_2(X)\notin\cQ(P_V)\big)$, where the event is assumed measurable (it is in all classes of \cref{prop:shifts}; otherwise read $\Prob$ as an inner probability); in particular $P_1\perp P_2$ when the shift is almost surely inadmissible.
\end{enumerate}
\end{theorem}

\proofidea Given $X$, the latent response is the law of $V$ translated by $\Psi_i(X)$, so the two conditional laws are equivalent exactly when the difference of the translations is an admissible shift; integrating over $X$ gives (a). The affinity of the joint laws is the average of the conditional affinities, which vanish on inadmissible shifts under the dichotomy property; this gives (b). The complete proof is in \cref{app:shift}.

\begin{proposition}[Admissible shifts beyond the Gaussian case]\label{prop:shifts}
The following classes of residual fields have explicit admissible shifts.
\begin{enumerate}[label=(\alph*),itemsep=1pt]
\item If $P_V$ is Gaussian, then $\cQ(P_V)=\CM_V$ and the dichotomy property holds.
\item Let $E$ be continuously embedded in $L^2(M)$, or a sequence space with continuous coordinates, and let $V=\sum_j\sigma_je_j\varphi_j$, the series converging almost surely in $E$, with independent $e_j$ whose densities $p_j$ are absolutely continuous and positive almost everywhere, and suppose that
\[
\sup_jJ(p_j)<\infty\qquad\text{and}\qquad\inf_j\inf_{t\ne0}\frac{H^2\big(p_j,p_j(\cdot-t)\big)}{\min(t^2,1)}>0,
\]
where $J$ is the Fisher information and $H^2$ the squared Hellinger distance. Then $\cQ(P_V)=\{h:\sum_jh_j^2/\sigma_j^2<\infty\}$, the Cameron--Martin space of the Gaussian field with the same scales, and the dichotomy property holds. Both conditions hold for independent and identically distributed $e_j$ from any fixed density with finite Fisher information that is positive almost everywhere: Laplace (as in sparsity-promoting Besov fields), logistic, hyperbolic secant, finite Gaussian mixtures, and Student-$t_\nu$ for every $\nu>0$.
\item If $V=\sqrt W\,G$ with $W>0$ independent of a centred Gaussian field $G$, then $\cQ(P_V)=\CM_G$ and the dichotomy property holds.
\end{enumerate}
In each case (in case (c) with $G$ $L$-stationary) the equivalence and singularity statements of \cref{thm:frontier}(i),(iii) hold verbatim, so the identification frontier $\cI=\infty$ of \cref{cor:matern} is the same in all these classes. The consistency statement (ii) and the variance statement of \cref{cor:matern} require additional moment conditions. In case (b) they hold when, in addition, $\E e_j=0$ and $\sup_j\E e_j^2<\infty$. Then the Cram\'er--Rao inequality gives $\Var e_j\ge1/J(p_j)$, so the coefficient variances are comparable to $\sigma_j^2$, and $\cI$ is computed with $g_j=\sigma_j^2$. In case (c) consistency holds as stated, by conditioning on $W$; the variance statement holds conditionally on $W$, and unconditionally when $\E W<\infty$.
\end{proposition}

\proofidea (a) is the Cameron--Martin theorem. In (b) the coefficients are independent, so Kakutani's theorem reduces the question to the Hellinger distance between a density and its translate, which the two conditions make comparable to $\min(t^2,1)$. In (c) one conditions on the mixing variable $W$. The complete proof, with the transfer to \cref{thm:frontier} and the Cauchy example below, is in \cref{app:shift}.

Without finite variance, GLS can miss part of the frontier even though identification is unaffected. For independent Cauchy coefficients (with $V$ viewed in the sequence space $E=\R^{\N}$; an $L^2(M)$-valued $V$ would also require $\sum_j\sigma_j<\infty$) and oracle weights $w_j=\sigma_j^{-2}$, the noise part of the GLS error is, conditionally on $X$, exactly Cauchy with scale $\sum_{j\le N}\sqrt{q_j}\,\abs{\zeta_j}\big/\sum_{j\le N}q_j\zeta_j^2$, where $\zeta_j=x_j/\sqrt{f_j}$. When $q_j\asymp j^{p-1}$ with $p>0$ this scale is of order $N^{(1-p)/2}$, so GLS is consistent only for $p>1$, that is for $\nu_X<\nu_V$ in Whittle--Mat\'ern terms, whereas the laws are mutually singular for all $p\ge0$ (experiment E7 confirms this numerically). Whether a score-based estimator attains the whole frontier for such residuals is open.

\subsection{Nonlinear covariate effects}

Let $M\subset\R^d$ be a bounded connected Lipschitz domain and $\eta=m\circ X+\mathsf A(X)+V$. We use the Sobolev scale $H^s$ and the following conditions.
\begin{enumerate}[label=(N\arabic*),leftmargin=3em,itemsep=1pt]
\item $P_V$ has the dichotomy property, and for some $t>d/2$, $\cQ(P_V)\subseteq H^t_{\mathrm{loc}}(M)$.
\item $X$ has continuous paths with $X(M)$ almost surely a non-degenerate interval, and \emph{exact local regularity} $s_X$: almost surely, (ER$-$) $X|_B\notin H^s(B)$ for every open ball $B\subseteq M$ and every $s\ge s_X$ with $s>d/2$, and (ER$+$) $X\in H^s(M)$ for every $s\in(d/2,s_X)$.
\item Response functions are smooth, $m\in C^\infty(\R)$, and confounding mechanisms are Cameron--Martin negligible: $\mathsf A\in\cA_{\mathrm{CM}}=\{\mathsf A:\mathsf A(X)\in\cQ(P_V)\ \text{a.s.}\}$.
\end{enumerate}
Condition (N1) holds with $t=\nu_V+d/2$ for restrictions to $M$ of $\R^d$-stationary Gaussian fields whose spectral density is comparable to $(1+\abs\omega^2)^{-(\nu_V+d/2)}$: their Cameron--Martin space is the reproducing kernel Hilbert space of the covariance \citep{Kanagawa2018}, which on $\R^d$ is $H^{\nu_V+d/2}(\R^d)$ with an equivalent norm \citep[Corollary~10.13]{Wendland2005}, and restricting a kernel to $M$ restricts its space \citep{Aronszajn1950}, so that $\CM_V=H^{\nu_V+d/2}(M)$. It holds for their scale mixtures (\cref{prop:shifts}(c)), whose admissible shifts form the same space. For Whittle--Mat\'ern fields defined through the SPDE with the Neumann or Dirichlet Laplacian, $\CM_V$ is the domain of $(\kappa^2+L)^{(\nu_V+d/2)/2}$; when $(\nu_V+d/2)/2$ is an integer, iterating interior elliptic regularity \citep[Section~6.3, Theorem~2]{Evans2010} shows that this domain lies in $H^{\nu_V+d/2}_{\mathrm{loc}}(M)$, so (N1) holds, and the same applies to the fields of \cref{prop:shifts}(b) built on these scales. By \cref{lem:regularity}, (N2) holds with $s_X=\nu_X$ for stationary Gaussian covariates with Mat\'ern-type spectral density and for their smooth monotone transforms, such as log-normal concentrations, prices or abundances $X=e^{G}$. The class $\cA_{\mathrm{CM}}$ contains every mechanism that acts through a sufficiently strong smoothing of arbitrary, possibly nonlinear, functionals of the covariate. It contains the constants whenever they lie in $\cQ(P_V)$ (as for periodic or Neumann operators and for restrictions of $\R^d$-stationary Mat\'ern-type fields), and then an intercept is absorbed into $\mathsf A$. The restriction is strong: in the linear case it requires $\sum_j\alpha_j^2f_j/g_j<\infty$, far more than vanishing confounding. Such confounding is undetectable and asymptotically harmless, and \cref{thm:nonlinear} is silent about nonlinear confounding that vanishes at high frequency without being Cameron--Martin negligible.

\begin{theorem}[Nonlinear covariate effects]\label{thm:nonlinear}
Assume (N1)--(N3) and write $t$ for the order in (N1).
\begin{enumerate}[label=(\alph*),itemsep=1pt]
\item \textbf{Shape is identified when the covariate is rough enough.} If $s_X\le t$, then for all $m_1,m_2$ and $\mathsf A_1,\mathsf A_2\in\cA_{\mathrm{CM}}$ the laws of $(X,\eta)$ satisfy $\TV\ge\Prob(m_1-m_2\text{ is non-constant on }X(M))$. They are mutually singular whenever $m_1-m_2$ is real-analytic and non-constant, for instance a non-constant polynomial.
\item \textbf{Nothing is identified when the covariate is too smooth.} If $s_X>t$ and $H^t(M)\subseteq\cQ(P_V)$, then for all $m_1,m_2\in C^\infty(\R)$ and every $\mathsf A\in\cA_{\mathrm{CM}}$ the laws under $(m_1,\mathsf A)$ and $(m_2,\mathsf A)$ are mutually absolutely continuous: no feature of the response function (slope, curvature, a smoothed threshold, monotonicity) is consistently estimable under any observation scheme.
\item \textbf{Parametric consistency.} Suppose $s_X\le t$, $V$ is Gaussian with spectrum $g$ in the eigenbasis of $L$ and $\CM_V\subseteq H^t_{\mathrm{loc}}(M)$, and $m_\theta=\sum_{k\le K}\theta_kb_k$ with $b_k\in C^\infty(\R)$ linearly independent modulo constants on every non-degenerate interval (for example $b_k(x)=x^k$). If at stage $N$ one observes $X$ and $(y_j)_{j\le N}$, then with $z_j=(\inner{b_k(X)}{\varphi_j})_{k\le K}$ and $I_N=\sum_{j\le N}z_jz_j^\top/g_j$, the GLS estimator $\hat\theta_N=I_N^{-1}\sum_{j\le N}z_jy_j/g_j$ converges to $\theta$ in probability for every $\mathsf A\in\cA_{\mathrm{CM}}$.
\end{enumerate}
\end{theorem}

\proofidea The difference $\delta=m_1-m_2$, composed with the covariate, is a shift of the latent response. Where $\delta'\ne0$, $\delta$ is locally invertible, so on some small ball the covariate is a smooth function of $\delta\circ X$. If $\delta\circ X$ were an admissible shift, and hence $H^t$-smooth, the covariate would be $H^t$-smooth on that ball, which exact local regularity excludes when $s_X\le t$; this gives (a). When the covariate is smoother than order $t$, every smooth function of it is an admissible shift, which gives (b). In (c) the same local inversion shows that the information matrix diverges in every direction. The complete proof is in \cref{app:nonlinear}, and the regularity lemmas are proved in \cref{app:regularity}.

For Mat\'ern-type fields $s_X=\nu_X$ and $t=\nu_V+d/2$, so the threshold $s_X\le t$ is exactly the linear frontier $\nu_X\le\nu_V+d/2$: nonlinearity costs nothing in identifiability. Part (b) is the strongest negative statement of the paper. For residual fields whose admissible shifts contain $H^{\nu_V+d/2}(M)$, such as the restrictions of $\R^d$-stationary Mat\'ern-type fields just discussed, a covariate more than $d/2$ orders smoother than the residual confounder leaves purely spatial data unable even to reveal whether the response function is increasing. For SPDE residual fields with Dirichlet or Neumann boundary conditions, $H^t(M)\not\subseteq\cQ(P_V)$ in general, because the domains of powers of $L$ impose boundary conditions, and (b) does not apply as stated.

\subsection{Observation channels and predicted covariates}

\begin{theorem}[Observation channels]\label{thm:channels}
\begin{enumerate}[label=(\alph*),itemsep=1pt]
\item The non-identification statements \cref{thm:frontier}(i),(iii), \cref{thm:shift}(a) and \cref{thm:nonlinear}(b) hold for any data generated from $(X,\eta)$ by a Markov kernel, with any design and sample size; \cref{cor:kriging}(a) holds for any data generated from $(\hat X,Y)$, and \cref{cor:kriging}(b) for any data generated from $(Z,Y)$, by a Markov kernel; and the non-identification statements of \cref{thm:resolution}(a),(b) hold for any data generated from the lattice data $D_\delta$ or $B_\delta$ by a Markov kernel.
\item In the linear model of \cref{sec:confounding} with $\cI=\infty$ and $\alpha\in\cA_0$, suppose the stage-$n$ data yield estimators $\hat x_j,\hat\eta_j$ of the coefficients of $X$ and $\eta=Y$ such that, for deterministic $a_J\uparrow\infty$ and $\epsilon_n\to0$, and with $h$ tail-equivalent to $g$, the quantities $\sum_{j\le J}(\abs{\hat x_j-x_j}+\abs{\hat\eta_j-\eta_j})/h_j\big/(\epsilon_na_J)$ are tight uniformly in $n$ and $J$. Then the sieve GLS estimator with $J_n=\max\{J:a_J\le\epsilon_n^{-1/2}\}$ is consistent. The condition holds for (i) a nugget of any distribution with finite variance, (ii) classical measurement error in $X$, (iii) Poisson counts with intensity $\Lambda_n\rho\,e^{\eta}$, for a known $\rho$ bounded away from zero and infinity, aggregated over cells $A_i$ with $\Lambda_n\min_i\mu(A_i)\to\infty$, and (iv) binomial counts with success probabilities $1/(1+e^{-\eta(s_i)})$ whose numbers of trials grow uniformly over cells ($\min_im_i\to\infty$); in (iii) and (iv) $\eta$ has H\"older continuous paths. In (i)--(iv) the design, the eigenfunctions and the covariate are as in \cref{prop:pointsamples}(b): partitions with $\ell_n\to0$ and $\max_iw_i\to0$, Lipschitz eigenfunctions, H\"older paths of $X$ and $U$, and $X$ observed at the points $s_i$ (with classical error in (ii)); $a_J$ is as in \cref{prop:pointsamples}(b).
\end{enumerate}
\end{theorem}

\proofidea (a) Data obtained through a Markov kernel from equivalent laws have equivalent laws. (b) If the coefficient estimates are accurate enough relative to the working spectrum, the sieve GLS estimator is a small perturbation of the complete-record estimator, as in \cref{prop:pointsamples}(b); the examples verify the required accuracy through moment bounds for averages of noise and of log-transformed counts. The complete proof is in \cref{app:channels}.

Classical measurement error is therefore asymptotically harmless for identification: it adds independent noise that averages out at every fixed scale. Smoothing is different, because it removes the high-frequency content of the covariate for good.

\begin{corollary}[Predicted covariates]\label{cor:kriging}
\begin{enumerate}[label=(\alph*),itemsep=1pt]
\item \textbf{Regression on a predicted surface.} Let $\hat X$ be the covariate field entering the response analysis, and model $Y=\beta\hat X+\mathsf A(\hat X)+V$ as in \eqref{eq:latent}, with the law $P_V$ of $V$ not depending on $\beta$ and $\cQ(P_V)$ a linear space. If $\hat X\in\cQ(P_V)$ almost surely, then the coefficient $\beta$ is not consistently estimable under any observation scheme, whatever $\mathsf A$. In particular, let $\hat X=\sum_{i\le K}c_ik_\ast(\cdot,m_i)$ be a kriging predictor from a fixed set of $K$ stations with a Whittle--Mat\'ern kernel $k_\ast$ of smoothness $\nu_\ast$, on a support whose heat kernel satisfies the on-diagonal bound $p_\tau(x,x)=\sum_je^{-\tau\lambda_j}\varphi_j(x)^2\le C\tau^{-d_s/2}$ for $\tau\le1$, and let $V$ be Gaussian with $g_j\asymp(1+\lambda_j)^{-(\nu_V+d_s/2)}$. Then $\hat X\in\CM_V$, and $\beta$ is not identified, whenever $2\nu_\ast>\nu_V$.
\item \textbf{The structural effect.} Assume the model of \cref{sec:confounding} with continuous eigenfunctions, with bounded $\alpha$ (as every $\alpha\in\cA_0$ is), and with $\sum_jf_j\varphi_j(m)^2<\infty$ at the stations. Suppose the covariate is observed only at $K$ fixed stations, $Z=(X(m_1),\dots,X(m_K))$, where $X(m)$ is the $L^2(\Prob)$-limit of $\sum_jx_j\varphi_j(m)$ (the point value when $X$ has continuous paths, \cref{lem:conditioning}), with the covariance matrix of $Z$ nonsingular. Suppose the response, an $L^2(M)$-valued Gaussian element, is observed through any Markov kernel. If $\beta\ne0$, the laws of $(Z,Y)$ under $(\beta,\alpha)$ and under $(-\beta,-\alpha)$ are mutually absolutely continuous, so the sign of $\beta$ is not consistently estimable for any confounding class closed under sign change, such as $\cA_0$. If, in addition, the variance of $X$ is unknown, then for every $\beta'\ne0$ the laws under $(\beta,\alpha,f)$ and under $(\beta',\alpha\beta'/\beta,(\beta/\beta')^2f)$ are mutually absolutely continuous: no two non-zero values of $\beta$ are distinguished. Both statements hold whatever the residual spectrum and the response sampling scheme.
\end{enumerate}
\end{corollary}

\proofidea (a) A kriging predictor from a fixed network is a finite combination of kernel sections; a smoothness count with the heat-kernel bound shows that each section lies in the Cameron--Martin space of $V$ when $2\nu_\ast>\nu_V$, so changing $\beta$ is an admissible shift. (b) Given the station values $Z$, the response is a Gaussian element with mean $\gamma(L)\E[X\mid Z]$, a finite combination of fixed functions, and with a conditional covariance that is unchanged under $(\beta,\alpha)\mapsto(-\beta,-\alpha)$. That mean is an admissible shift of the conditional noise, so the two conditional laws are equivalent for every $Z$. The complete proof, including a lemma on conditioning a Gaussian field on point values, is in \cref{app:channels}.

Part (a) concerns the working regression on $\hat X$, with the residual law left unrestricted. That GLS cannot be consistent when the regressor lies in the Cameron--Martin space of the noise is Theorem~1(iii) of \citet{BolinWallin2026} (see also \citet{Wang2020kriging} for fixed regressors); part (a) extends this to every estimator. The condition $2\nu_\ast>\nu_V$ holds for essentially every practical choice; with exponential kernels on both sides, for instance, $\nu_\ast=\nu_V=1/2$. The reason is a smoothness count: $\hat X$ is a finite combination of kernel sections, which lie in $H^s$ for every $s<2\nu_\ast+d_s/2$, whereas $\CM_V$ is the Sobolev space of order $\nu_V+d_s/2$. In the structural model, however, the residual relative to $\hat X$, $(\beta+\alpha(L))(X-\hat X)+V$, depends on $\beta$, so (a) alone does not settle the structural question. Part (b) does. With a fixed network of stations, the covariate information is $K$-dimensional, and every shift it can produce in the response lies in the Cameron--Martin space of the response noise. The sign of the effect is therefore never identified. If the covariate variance is unknown, as it is when only $K$ station values are available, no two non-zero effects are distinguished, whatever the residual spectrum. If all spectra were known exactly, $\abs{\beta}$ might still be identified from the second-order structure of $Y$ alone, through the Feldman--H\'ajek theorem; that identification owes nothing to the covariate data. Two-stage analyses that regress a response on a covariate surface interpolated from a sparse fixed network of stations, for instance health outcomes on kriged pollution surfaces \citep{Gryparis2009,SzpiroPaciorek2013}, or ecological and economic responses on interpolated climate or soil surfaces, therefore cannot be rescued by any spatial adjustment of the second stage, however many response observations are added.
\section{Numerical experiments}\label{sec:experiments}

We illustrate every prediction numerically. For the impossibility results, simulations can only show that natural estimators fail as predicted; the constructions of E4 verify the resolution theorem directly. Experiments E1, E3, E4 and E7 work in the spectral model \eqref{eq:model}, where exact computations are possible for up to $2^{18}$ modes. Experiments E2, E5, E6, E8 and E9 work in physical space on the torus, with point samples, nuggets, nonlinear features and counts; there infill is emulated by subsampling a single fine realization. Experiment E10 checks the areal models of \cref{cor:areal} in both settings. Unless they are labelled as misspecified (E1, E9) or use a prescribed working model (the working spectrum of E5 and the graph-model weights of E10), GLS fits use the true residual spectrum, and all physical-space experiments use periodic domains. Settings are given in \cref{app:sims}, and all code is available as supplementary material.

\begin{figure}[t]
\centering
\includegraphics[width=\textwidth]{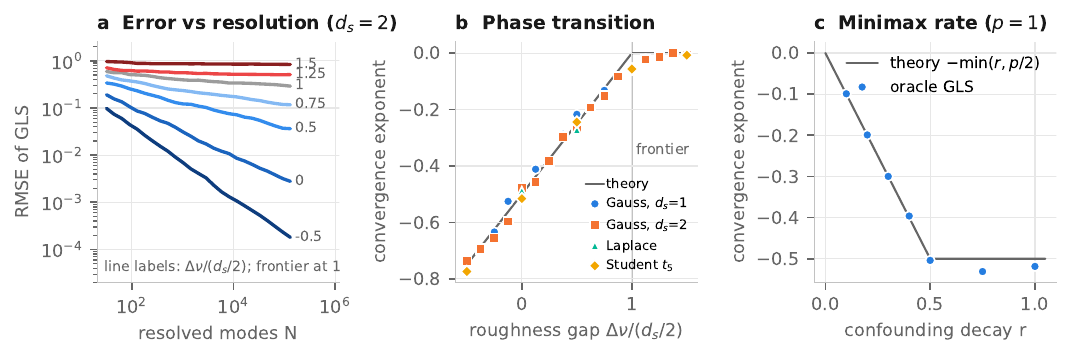}
\caption{The Cameron--Martin frontier (E1, E7). (a) RMSE of oracle GLS versus resolved modes $N$, $d_s=2$; line labels give $\Delta\nu/(d_s/2)$. (b) Fitted convergence exponents versus the roughness gap for Gaussian residuals with $d_s=1,2$ and for Laplace and Student-$t_5$ residuals ($d_s=1$); the solid line is the theory $-\max\{0,(1-\Delta\nu/(d_s/2))/2\}$. (c) Exponent versus confounding decay $r$ at $p=1$, with the minimax prediction $-\min(r,p/2)$.}\label{fig:frontier}
\end{figure}

\paragraph{E1: the frontier, its location and its rates.} \Cref{fig:frontier}a shows the root mean squared error (RMSE) of GLS against the number of resolved modes for $d_s=2$. Inside the frontier the error decays polynomially; at the boundary it decays logarithmically; outside it does not decay. Panel (b) plots the fitted convergence exponent against the roughness gap $\Delta\nu/(d_s/2)$ for $d_s=1$ and $d_s=2$: both collapse onto the theoretical curve $-\max\{0,(1-\Delta\nu/(d_s/2))/2\}$, with the kink exactly at the frontier. Panel (c) varies the decay $r$ of the confounding at $\Delta\nu=0$ ($p=1$); the exponents match $-\min(r,1/2)$ of \cref{thm:minimax} to within $0.04$. \Cref{tab:estimators} compares estimators inside and outside the frontier for $d_s=1$, where a valid undersmoothed working model ($\nu_h=0.05$) lies below the window in both configurations. Every verdict of \cref{thm:window,prop:estimators} is reproduced, including the exponent $-0.117$ of eigenvector filtering with $K=\sqrt N$ (theory $-0.125$), the stalled undersmoothed GLS, and the divergence of oversmoothed and high-pass estimators outside the frontier (exponents $+0.139$, theory $+0.125$).

\begin{table}[!tbp]
\centering\small
\caption{Estimators inside ($\Delta\nu/(d_s/2)=0.5$, $\nu_X=0.75$) and outside ($1.25$, $\nu_X=1.125$) the frontier, $d_s=1$, $\nu_V=0.5$, $N=2^{17}$ modes, 200 replicates. ``Exponent'' is the fitted slope of $\log$ RMSE on $\log N$; the theory column gives the prediction of \cref{thm:window,prop:estimators} (0 marks a non-vanishing limit).}\label{tab:estimators}
\begin{tabular}{lcccccc}
\toprule
 & \multicolumn{3}{c}{Inside} & \multicolumn{3}{c}{Outside} \\
\cmidrule(lr){2-4}\cmidrule(lr){5-7}
Estimator & RMSE & exponent & theory & RMSE & exponent & theory \\
\midrule
OLS (= RSR) & 0.736 & $0.000$ & 0 & 1.253 & $0.000$ & 0 \\
Eigenvector filtering, $K=20$ & 0.356 & $0.000$ & 0 & 1.648 & $0.000$ & 0 \\
Eigenvector filtering, $K=\sqrt N$ & 0.186 & $-0.117$ & $-0.125$ & 2.270 & $+0.058$ & $+0.063$ \\
GLS, undersmoothed ($\nu_h=0.05$) & 0.300 & $-0.010$ & 0 & 0.866 & $0.000$ & 0 \\
GLS, oversmoothed ($\nu_h=\nu_V+1$) & 0.060 & $-0.266$ & $-0.25$ & 3.729 & $+0.139$ & $+0.125$ \\
GLS, oracle & 0.034 & $-0.259$ & $-0.25$ & 0.546 & $-0.010$ & 0 \\
High-pass GLS, $K=N/4$ & 0.047 & $-0.274$ & $-0.25$ & 3.425 & $+0.139$ & $+0.125$ \\
\bottomrule
\end{tabular}
\end{table}

\begin{figure}[t]
\centering
\includegraphics[width=0.86\textwidth]{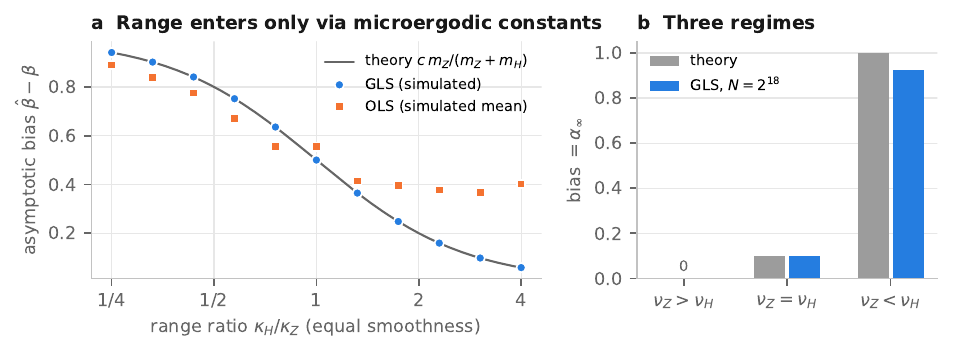}
\caption{What spatial GLS estimates (E3). (a) Asymptotic bias versus the range ratio at equal smoothness: simulated GLS (dots) against the microergodic formula of \cref{prop:whatgls}(b), and the simulated mean of OLS. (b) The three smoothness regimes of \eqref{eq:threeregimes}.}\label{fig:glslimit}
\end{figure}

\paragraph{E3: what GLS estimates.} In the model of \cref{prop:whatgls}(b) with equal smoothness, the simulated asymptotic bias of GLS matches $c\,m_Z/(m_Z+m_H)$ to within $3\times10^{-4}$ over a sixteen-fold range of $\kappa_H/\kappa_Z$ (\cref{fig:glslimit}a); OLS behaves quite differently. Across the three smoothness regimes the bias is $0.0003$, $0.1001$ and $0.925$, against theoretical limits $0$, $0.1$ and $1$; the last approaches its limit slowly because $\alpha_j\to c$ at rate $j^{-1/2}$.

\begin{figure}[t]
\centering
\includegraphics[width=\textwidth]{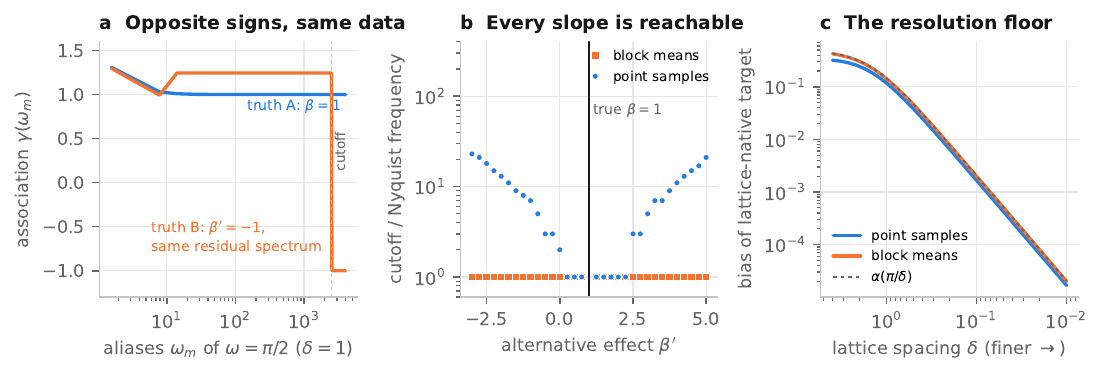}
\caption{The resolution theorem (E4). (a) Associations $\gamma(\omega_m)$ at the aliases of $\omega=\pi/2$ for two continuum models with effects $\beta=1$ and $\beta'=-1$ and the same residual spectrum; their lattice laws coincide. (b) Cutoff, in multiples of the Nyquist frequency, needed by the construction of \cref{thm:resolution}(a) to reach each $\beta'$. (c) Bias of the lattice-native target at the Nyquist frequency versus the spacing $\delta$, against $\alpha(\pi/\delta)$.}\label{fig:resolution}
\end{figure}

\paragraph{E4: the resolution theorem.} We implemented the constructions in the proof of \cref{thm:resolution} at $\delta=1$ (configuration in \cref{app:sims}, $\beta=1$). With an adjusted residual spectrum (part a), $\beta'=-1$ and $\beta'=3$ are reached with cutoffs of $7$ and $5$ times the Nyquist frequency for point samples, and of $1$ times for block means. With the residual spectrum held fixed (part b, point samples), $\beta'=-1$ needs a cutoff of $800$ times the Nyquist frequency, because the alias variance of this $\alpha$ is smallest at the Nyquist frequency (\cref{fig:resolution}a). In every case the aliased spectral matrices of the two models agree to a relative error of about $10^{-15}$. This check uses the same alias sums as the construction. As an independent check, we computed the lattice covariances of both models at lags $0,\dots,10$ by direct quadrature of the continuum spectral differences over $\R$, without alias sums. They agree to within about $10^{-8}$ at all lags for block means and at lags $1,\dots,10$ for point samples. At lag $0$ for point samples the residuals, up to $1.2\times10^{-4}$, come from truncating at $\abs m\le3000$ the alias sums used to build the alternative. Panel (b) shows the cutoff needed in part (a) for $\beta'\in[-3,5]$. It is finite throughout and grows roughly quadratically in $\abs{\beta'-\beta}$ for point samples, because for this $f_X$ the alias mass beyond $R$ decays like $1/R$. Panel (c) shows that the target of a lattice-native analysis exceeds $\beta$ by $0.74$--$0.82\,\alpha(\pi/\delta)$ for point samples and $0.98$--$0.99\,\alpha(\pi/\delta)$ for block means, and vanishes only as $\delta\to0$.

\begin{figure}[t]
\centering
\includegraphics[width=0.92\textwidth]{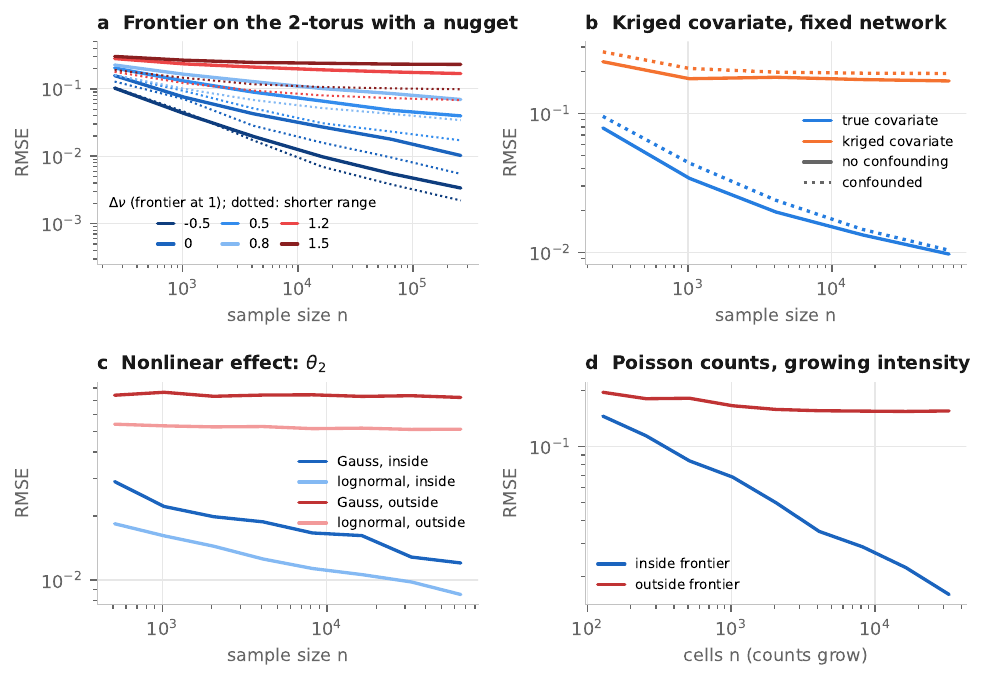}
\caption{Physical-space experiments. (a) GLS on the 2-torus with a nugget; solid and dotted lines use different ranges. (b) Regression on a kriged covariate from a fixed network versus the true covariate. (c) Quadratic coefficient of a nonlinear response function under a nonlinear, non-local confounding mechanism. (d) Poisson counts with growing intensity.}\label{fig:physical}
\end{figure}

\paragraph{E2, E5, E6, E8: physical space.} \Cref{fig:physical} reports the four physical-space experiments. On the 2-torus with a nugget (panel a), the RMSE of GLS keeps falling inside the frontier and flattens outside it. Fitted exponents over the last sixteen-fold increase of $n$ are $-0.39,-0.35,-0.18,-0.14,-0.05,-0.01$ for $\Delta\nu=-0.5,0,0.5,0.8,1.2,1.5$ in one range setting and $-0.42,-0.38,-0.21,-0.15,-0.06,-0.03$ in a setting with ranges $8/3$ times shorter: ranges change levels, not the verdict. Outside the frontier the RMSE is still declining slowly at the largest $n$; the theory predicts a positive floor, which we have not computed. With a kriged covariate from 36 fixed stations (panel b), the RMSE stalls at $0.17$ without confounding and at $0.19$ with confounding (bias $0.10$), while regression on the true covariate converges. For a quadratic response function with a nonlinear, non-local confounding mechanism (panel c), both coefficients converge inside the frontier, for Gaussian and log-normal covariates, and neither does outside it. With Poisson counts whose expected values grow with resolution (panel d), the RMSE falls from $0.146$ to $0.016$ inside the frontier and stalls at $0.156$ (bias $0.095$) outside. This is the favourable regime of \cref{thm:channels}(b); in real areal count data (disease cases, crimes, species records), counts per unit shrink as units are refined. The confounding mechanism of panel (c) is Cameron--Martin negligible by construction, so E6 cannot reveal failures under nonlinear confounding that vanishes at high frequency without being negligible.

\paragraph{E7: non-Gaussian residuals.} With Laplace and Student-$t_5$ residual coefficients the convergence exponents of GLS track the Gaussian ones (\cref{fig:frontier}b). With Cauchy coefficients, which have finite Fisher information but infinite variance, GLS behaves as predicted after \cref{prop:shifts} (oracle weights, no confounding, $d_s=1$). For $p=1.5$ it is consistent, with exponent $-0.21$ (theory $-0.25$). For $p=1$ its median absolute error stays near $0.9$ (from $0.82$ at $N=32$ to $0.94$ at $N=2^{17}$). For $p=0.5$ it diverges, from $1.74$ to $10.0$ (exponent $+0.27$, theory $+0.25$). Over the same range the Gaussian and $t_5$ errors fall to $0.002$ ($p=1$) and $0.025$ ($p=0.5$). Identification is the same in all three cases; only GLS fails.

\paragraph{E9: GLS with a nugget and misspecified working models.} This regime, point sampling with a working nugget and a misspecified working covariance, is not covered by our theorems. On the 1-torus with $\nu_V=0.5$, $\nu_X=0.9$ (inside the frontier) and nugget variance $0.01$, the RMSE of GLS for $n=2^9,2^{11},\dots,2^{17}$ is as follows (\cref{tab:nugget}). With a working smoothness below the window ($\nu_h=0.05$) it stalls near $0.18$ with bias $0.11$, as \cref{thm:window} predicts. With the correct smoothness it falls from $0.12$ to $0.077$. With oversmoothed working models it falls only from $0.165$ to $0.143$ ($\nu_h=2$), or not at all at about $0.20$ ($\nu_h=6$). The correct smoothness with an eight times shorter working range gives $0.26$ falling to $0.22$, with bias $0.14$. With a nugget, an oversmoothed working model treats the high frequencies as noise. The number of modes it effectively uses then grows like a small power of $n$, so oversmoothing that is harmless under (O1) slows convergence drastically, and the working range matters at realistic sample sizes.

\begin{table}[!tbp]
\centering\small
\caption{GLS with a working nugget under point sampling (E9): RMSE over 300 replicates, $\nu_V=0.5$, $\nu_X=0.9$, nugget variance $0.01$; the window of \cref{thm:window} (without nugget) is $\nu_h\ge0.4$.}\label{tab:nugget}
\begin{tabular}{lcccccc}
\toprule
Working model & $n=2^9$ & $2^{11}$ & $2^{13}$ & $2^{15}$ & $2^{17}$ & bias at $2^{17}$ \\
\midrule
$\nu_h=0.05$ (below window) & 0.189 & 0.185 & 0.183 & 0.182 & 0.182 & $+0.109$ \\
$\nu_h=0.5$ (correct) & 0.122 & 0.103 & 0.090 & 0.081 & 0.077 & $+0.022$ \\
$\nu_h=2$ (oversmoothed) & 0.165 & 0.158 & 0.153 & 0.149 & 0.143 & $0.000$ \\
$\nu_h=6$ (oversmoothed) & 0.201 & 0.199 & 0.201 & 0.204 & 0.204 & $+0.005$ \\
$\nu_h=0.5$, range$/8$ & 0.260 & 0.247 & 0.236 & 0.226 & 0.217 & $+0.138$ \\
\bottomrule
\end{tabular}
\end{table}

\paragraph{E10: areal models under refinement.} \Cref{tab:areal} checks \cref{cor:areal} on cycle graphs ($d_s=1$) whose nodes are point samples of one fine realization on the 1-torus, with confounding that vanishes at high frequency. Inside the frontier with $\nu_X=0.75<1$, the intrinsic CAR and BYM models and the Leroux model with $\rho_n$ approaching $1$ at the rate $n^{-2}$ converge at the rate of oracle GLS. The Leroux model with $\rho=0.9$ and the SAR error model with $\rho=0.4$ instead drift towards the error of OLS as $n$ grows: refining the data makes these models worse. With $\nu_X=1.5>1$, still inside the frontier, oracle GLS converges, but the intrinsic CAR and BYM models stall, as \cref{cor:areal}(b) predicts for covariates smoother than their window allows. The exponents of the converging estimators are slightly steeper than the asymptotic value $-0.25$, for oracle GLS as for the others. This point-sample design lies outside (O1), because the graph coefficients are then aliased sums of continuum coefficients; the spectral design of \cref{app:sims}, which matches the reading of graph refinement in \cref{cor:areal}, gives the same verdicts, and so does the design without confounding.

\begin{table}[!tbp]
\centering\small
\caption{Areal models under refinement (E10): cycle graphs with $n$ nodes sampled from one realization on the 1-torus, $\nu_V=\nu_X-0.25$ (inside the frontier, $d_s=1$), confounding that vanishes at high frequency, 100 replicates. RMSE at $n=2^7$ and $n=2^{15}$, fitted exponent, and the prediction of \cref{cor:areal}, which is stated for the spectral design of \cref{app:sims} (0 marks a non-vanishing limit).}\label{tab:areal}
\begin{tabular}{lcccccccc}
\toprule
 & \multicolumn{4}{c}{$\nu_X=0.75$} & \multicolumn{4}{c}{$\nu_X=1.5$} \\
\cmidrule(lr){2-5}\cmidrule(lr){6-9}
Estimator & $2^7$ & $2^{15}$ & exponent & theory & $2^7$ & $2^{15}$ & exponent & theory \\
\midrule
OLS & 0.301 & 0.301 & $0.000$ & 0 & 0.414 & 0.414 & $0.000$ & 0 \\
ICAR (neighbour differences) & 0.144 & 0.024 & $-0.293$ & $-0.25$ & 0.203 & 0.182 & $-0.003$ & 0 \\
BYM & 0.152 & 0.027 & $-0.288$ & $-0.25$ & 0.232 & 0.182 & $-0.008$ & 0 \\
Leroux, $\rho_n=1/(1+\kappa^2/n^2)$ & 0.134 & 0.025 & $-0.297$ & $-0.25$ & 0.279 & 0.259 & $-0.002$ & 0 \\
Leroux, $\rho=0.9$ & 0.151 & 0.300 & $+0.012$ & 0 & 0.338 & 0.414 & $0.000$ & 0 \\
SAR error, $\rho=0.4$ & 0.260 & 0.301 & $+0.001$ & 0 & 0.405 & 0.414 & $0.000$ & 0 \\
Oracle GLS & 0.135 & 0.026 & $-0.311$ & $-0.25$ & 0.164 & 0.028 & $-0.290$ & $-0.25$ \\
\bottomrule
\end{tabular}
\end{table}
\section{Implications for practice}\label{sec:practice}

The results turn several long-running questions into checkable statements.

\paragraph{Does the spatial random effect ``eat'' the covariate effect?} Inside the frontier, no. Under spectral infill, GLS with a fixed working spectrum whose smoothness is not too low converges to the finest-scale association $\beta+\alpha(\infty)$, which equals $\beta$ under vanishing confounding (\cref{thm:window,prop:whatgls}). On areal data the answer depends on the model: on refining regular graphs, Leroux, proper CAR and SAR models with fixed hyperparameters do not deconfound at all and converge to the limit of OLS, and on regular grids the intrinsic CAR and BYM models deconfound only when $\nu_X\le1$ (\cref{cor:areal}, experiment E10). Mixed models with estimated covariance parameters and nuggets are not covered by our theorems; experiment E9, which uses fixed working models with a nugget under point sampling, shows that such fits can converge slowly. Outside the frontier every method fails, so the instability of estimates across spatial specifications can be a symptom of a covariate that lies (nearly) in the Cameron--Martin space of the residual confounder rather than a defect of any particular model; fixed resolution, confounding that does not vanish at fine scales and small samples can produce it too.

\paragraph{Should restricted spatial regression be used?} No. Its point estimate is OLS, which retains a non-vanishing weight on the confounded low frequencies (\cref{prop:estimators}(a)), in agreement with \citet{Hanks2015}, \citet{KhanCalder2022} and \citet{Gilbert2025}.

\paragraph{Why do comparisons of adjustment methods disagree?} By \cref{prop:normalform}, every linear adjustment is a ratio of weighted sums over the scale-specific associations, a weighted average of them when $\omega=\tilde\omega\ge0$ as for OLS, GLS and eigenvector filtering (\cref{prop:normalform}(c)), so its bias in a simulation is fixed by where the simulated confounding sits relative to the method's weights. Inside the frontier with refining resolution, the consistent methods differ only in efficiency and finite-sample bias, and a comparison measures how quickly each moves its weight to the finest scales. Outside the frontier, or on a fixed grid, no method is consistent: some converge to a confounded target and others diverge (\cref{tab:estimators}), and the ranking reflects which estimate happens to lie closest to $\beta$ for the confounding pattern and sample size chosen. Designs that differ in the roughness gap, in the resolution or in the scale of the confounding should therefore rank the methods differently, which may explain the disagreements among the simulation and empirical comparisons of \citet{Urdangarin2023}, \citet{Tec2024space} and \citet{Pim2026review}. Reporting the roughness gap and the resolution of each design would make such comparisons interpretable and would test this explanation.

\paragraph{A four-step protocol.}
\begin{enumerate}[leftmargin=2em,itemsep=2pt]
\item \emph{Check the resolution.} If the response is available only on a fixed regular grid and the residual spectrum is unknown, the high-frequency assumption cannot identify the effect, not even its sign (\cref{thm:resolution}); for irregular areal units the same is expected but not proved. Report spectral-adjustment estimates as conditional on an assumption about confounding at the scale of the areal units, or on a parametric form for the coherence, and seek finer-resolution data where possible.
\item \emph{Estimate the roughness gap.} The smoothness $\nu_X$ of the covariate is microergodic and can be estimated from the covariate alone, for example from increments or the behaviour of the variogram near the origin \citep{KentWood1997,Loh2021}. Because confounding vanishes at high frequencies under the identifying assumption, the smoothness of the residual confounder can be read off the high-frequency behaviour of the response residuals. If $\hat\nu_X-\hat\nu_V$ clearly exceeds $d_s/2$ (one half on networks, one on surfaces), no purely spatial analysis will identify the effect: temporal variation, negative controls \citep{Lipsitch2010}, instruments or natural experiments are needed.
\item \emph{Choose an estimator inside the window.} Avoid a working smoothness below $\nu_X-d_s/2$, which is fatal (\cref{thm:window}). Do not oversmooth heavily either. With a nugget, a strongly oversmoothed working model can make GLS nearly as poor as an undersmoothed fit at realistic sample sizes (experiment E9), so estimate the working smoothness or use high-pass GLS. Eigenvector filtering with a growing number of eigenvectors is another valid alternative, with different finite-sample trade-offs (\cref{prop:estimators}). On areal data, prefer intrinsic models, such as the intrinsic CAR model or the regression of neighbour differences, to Leroux, proper CAR or SAR models with fixed hyperparameters, which do not deconfound under refinement; the intrinsic CAR model in turn requires $\nu_X\le1$, and second-order intrinsic models, which suit smoother covariates, require $\nu_X\le2$ (\cref{cor:areal}).
\item \emph{Use covariate data that refine with the response data.} With a fixed network of stations, neither the coefficient on a kriged surface nor the sign of the effect is identified (\cref{cor:kriging}); classical measurement error in densely observed covariates is harmless (\cref{thm:channels}).
\end{enumerate}
The same protocol applies to nonlinear covariate effects under Cameron--Martin-negligible confounding and to count responses with growing expected counts (\cref{thm:nonlinear,thm:channels}).

\section{Discussion}\label{sec:discussion}

We have argued that the identification of covariate effects from spatial data has an endgame; \cref{tab:settled} lists the questions it settles. Three conditions---resolution, roughness and vanishing confounding---are jointly sufficient, and each is necessary in the senses made precise in \cref{thm:frontier,thm:resolution}. As conditions for identification they do not depend on the choice among the estimators of \cref{tab:zoo}, which decides only whether and how fast an estimator reaches the frontier (\cref{thm:window,cor:areal}), and their impossibility parts are robust to nonlinearity, non-Gaussianity and the observation channel. The organizing device is old: equivalence and singularity of measures, and Grenander's insight that what a single continuous record reveals is exactly what distinguishes non-equivalent measures. What is new is how far the device reaches once the model zoo is written in its spectral normal form.

\paragraph{Limitations.} Our positive results assume operator stationarity, which is exact for Whittle--Mat\'ern fields on manifolds and networks, for stationary fields on tori and for graph-based models on regular graphs, and approximate otherwise. For $\R^d$-stationary fields observed on bounded domains, the frontier was established by \citet{DattaStein2025}. Transferring our impossibility half through \cref{thm:shift} requires the admissible shifts of the conditional residual $V=U-\E[U\mid X|_M]$, which is not of Mat\'ern type in general; we have not carried this out. \Cref{thm:nonlinear}(a) with linear $m$ does give mutual singularity at the boundary $\nu_X=\nu_V+d/2$ for restrictions of $\R^d$-stationary Mat\'ern-type fields, when the residual is independent of the covariate and the confounding mechanism is Cameron--Martin negligible. Consistency at that boundary remains open. Accordingly, ``boundary included'' refers to operator-stationary fields. The resolution theorem is proved for the cubic lattice with point values or cubic block means; irregular partitions call for a non-stationary version. The nonlinear theory uses Cameron--Martin-negligible confounding mechanisms. In the linear case vanishing confounding suffices; closing this gap for nonlinear mechanisms is open, and our experiment E6 does not probe it. The GLS theory assumes fixed working spectra without a nugget under spectral infill, the areal results assume fixed hyperparameters or a prescribed rescaling, and the experiments use prescribed working spectra, mostly the true residual spectrum, on periodic domains. GLS with a nugget and estimated covariance parameters under point sampling, and estimators that attain the frontier for heavy-tailed residuals, are open. Other natural next steps are rates for nonparametric response functions, spatio-temporal designs, interference, and confounders that are deterministic functions of location. For the latter, identification requires covariate variation not explained by location \citep[Proposition~1]{Gilbert2025}.

\paragraph{The broader point.} Much of the spatial confounding literature has compared models, and the census in \cref{tab:census} shows how many there now are. The theory here suggests comparing \emph{data}: their resolution, the roughness of the covariate, and what can plausibly be assumed about confounding at the finest scales the data resolve. These three properties, rather than the choice among nonparametric spatial models, decide whether a covariate effect can be learned; parametric assumptions can restore identification, but they do so by asserting more about confounding than its disappearance at fine scales, and the data can check such assertions at most in part. For geography, the results give the scale problem a precise form on regular grids: below the grain of the data, confounding is not a nuisance that a better model can remove, but a limit on what any model can know without assumptions about confounding at the scale of the units.

\section*{Use of generative AI}
The author directed the research, guided the development of the manuscript and its proofs, and reviewed the completed work. Under the author's direction, Claude (Anthropic), a large language model, provided substantial assistance in drafting the text and proofs, producing the figures, writing the simulation code that reproduces the numerical experiments, and checking the proofs and references. The author takes full responsibility for the final manuscript and all claims it contains.

\section*{Code availability}
Code that reproduces every experiment in \cref{sec:experiments} and \cref{app:sims} is provided as ancillary files with this article.

\bibliographystyle{plainnat}
\bibliography{refs}

\appendix
\section{Complete proofs}\label{app:proofs}

\subsection{Guide to the proofs}\label{app:guide}

This appendix proves every result stated in the paper. Each proof of a result stated in \crefrange{sec:framework}{sec:beyond} opens with a short paragraph giving its idea, which matches the proof idea stated after the result in the main text, and then proceeds in numbered steps; the auxiliary lemmas are proved directly. Classical results that we use without proof are stated precisely, with references, in \cref{app:classical}; every other auxiliary fact is proved in \cref{app:prelim}, \cref{app:regularity} or \cref{app:channels}. \Cref{tab:proofmap} shows where each result is proved, which tools it uses and on which other results it rests.

\emph{Conventions.} All random elements are defined on one probability space. The \emph{law} of a random element is its distribution on the relevant Polish space with its Borel $\sigma$-field; all spaces used below ($\R^{\N}$, $L^2(M)$, $C(M)$, separable Fr\'echet spaces and their products) are Polish. Two probability measures $P,Q$ are \emph{equivalent}, $P\equivm Q$, if they are mutually absolutely continuous, and \emph{singular}, $P\perpm Q$, if they are mutually singular. For deterministic sequences, $a_j\asymp b_j$ means $c\le a_j/b_j\le C$ for constants $0<c\le C<\infty$ and all large $j$ (similarly with $N$ as index). For random quantities, $A_N\asymp_PB_N$ means that $A_N/B_N$ and $B_N/A_N$ are both bounded in probability. Unless they are introduced as specific quantities (such as the shift $c$ of \cref{thm:frontier}(iii)), the letters $c,C$ denote positive constants that may change from line to line but never depend on the index ($j$, $N$ or $n$) under consideration. In the Gaussian model of \cref{sec:confounding} we write $\xi_j=x_j/\sqrt{f_j}$; the $\xi_j$ are independent standard normal and independent of $V$.

\subsection{Classical results used without proof}\label{app:classical}

\begin{enumerate}[label=(C\arabic*),leftmargin=3em,itemsep=2pt]
\item \emph{Kakutani's theorem} \citep{Kakutani1948}. Let $P=\bigotimes_jP_j$ and $Q=\bigotimes_jQ_j$ be product probability measures on a countable product of measurable spaces, with $P_j\equivm Q_j$ for every $j$. Then $\aff(P,Q)=\prod_j\aff(P_j,Q_j)$. If $\sum_jH^2(P_j,Q_j)<\infty$, then $P\equivm Q$, and $dP/dQ$ is the $Q$-almost sure and $L^1(Q)$ limit of the products $\prod_{j\le N}(dP_j/dQ_j)$, read as functions of the first $N$ coordinates. Otherwise $P\perpm Q$.
\item \emph{Kolmogorov's convergence criterion} \citep[Section~2.5]{Durrett2019}. If $\varepsilon_1,\varepsilon_2,\dots$ are independent real random variables with $\E\varepsilon_j=0$ and $\sum_j\E\varepsilon_j^2<\infty$, then $\sum_j\varepsilon_j$ converges almost surely and in $L^2$.
\item \emph{Lusin--Souslin theorem} \citep[Theorem~15.1 and Corollary~15.2]{Kechris1995}. If $S$ and $T$ are Polish spaces and $\Phi:S\to T$ is Borel measurable and injective, then $\Phi(B)$ is a Borel subset of $T$ for every Borel set $B\subseteq S$.
\item \emph{Measurable densities} \citep[Theorem~V.58]{DellacherieMeyer1982}. If $\cX,\mathcal Y$ are standard Borel spaces and $K,K'$ are Markov kernels from $\cX$ to $\mathcal Y$, there are jointly measurable functions $k,k'\ge0$ on $\cX\times\mathcal Y$ with $K(x,dy)=k(x,y)\lambda_x(dy)$ and $K'(x,dy)=k'(x,y)\lambda_x(dy)$, where $\lambda_x=\tfrac12\{K(x,\cdot)+K'(x,\cdot)\}$.
\item \emph{Gaussian measures} \citep[Sections~2.4, 2.5 and~2.8]{Bogachev1998}. Let $\gamma$ be a centred Gaussian measure on a separable Fr\'echet space $E$, with Cameron--Martin space $\CM(\gamma)=\{h\in E:\abs{h}_{\CM(\gamma)}<\infty\}$, where $\abs{h}_{\CM(\gamma)}=\sup\{\ell(h):\ell\in E^\ast,\ \int\ell^2d\gamma\le1\}$. (a) \emph{Cameron--Martin theorem.} The shifted measure $\gamma(\cdot-h)$ is equivalent to $\gamma$ if $h\in\CM(\gamma)$ and singular to $\gamma$ otherwise. (b) \emph{Zero--one law} \citep[see also][]{Kallianpur1970}. Every $\gamma$-measurable linear subspace of $E$ has $\gamma$-measure $0$ or $1$. (c) \emph{Fernique's theorem} \citep{Fernique1970}. If $E$ is a separable Banach space, then $\int\exp(\epsilon\norm x^2)\,\gamma(dx)<\infty$ for some $\epsilon>0$; in particular $\int\norm x^2\gamma(dx)<\infty$. For an $L$-stationary Gaussian field with spectrum $(s_j)$, viewed in $E=L^2(M)$, the formula for $\abs{\cdot}_{\CM(\gamma)}$ with $\ell=\inner{\cdot}{a}$ gives $\sup\{\sum_jh_ja_j:\sum_ja_j^2s_j\le1\}=(\sum_jh_j^2/s_j)^{1/2}$, which is the Cameron--Martin norm of \cref{sec:confounding}.
\item \emph{Sobolev spaces} \citep[Chapter~4 and Section~5.3]{RunstSickel1996}. For $s\ge0$, $H^s(\R^d)=\{u\in L^2(\R^d):\int(1+\abs\omega^2)^s\abs{\hat u(\omega)}^2d\omega<\infty\}$. For open $\Omega\subseteq\R^d$, $H^s(\Omega)$ is the space of restrictions to $\Omega$ of elements of $H^s(\R^d)$, so that $H^s(\Omega)\subseteq H^t(\Omega)$ for $s\ge t$; and $H^s_{\mathrm{loc}}(M)$ is the space of functions $u$ on $M$ with $u|_B\in H^s(B)$ for every ball $B$ with $\bar B\subset M$. Multiplication by a function in $C_c^\infty(\R^d)$ maps $H^s(\R^d)$ into itself. If $s>d/2$ and $F\in C^\infty(\R)$ with $F(0)=0$, then $u\mapsto F\circ u$ maps $H^s(\R^d)$ into itself.
\item \emph{Mercer's theorem} \citep[Theorem~4.49]{SteinwartChristmann2008}. Let $M$ be a compact metric space, $\mu$ a finite Borel measure with $\operatorname{supp}\mu=M$, and $k:M\times M\to\R$ a continuous, symmetric, positive semidefinite kernel. Then the integral operator $f\mapsto\int k(\cdot,t)f(t)\,\mu(dt)$ on $L^2(M,\mu)$ has an orthonormal system of eigenfunctions $e_j$, continuous, with eigenvalues $s_j>0$ spanning the closure of its range, and $k(s,t)=\sum_js_je_j(s)e_j(t)$, the series converging absolutely and uniformly on $M\times M$.
\item \emph{Spectral representation of stationary fields} \citep[Section~2.5]{Stein1999}. A centred, mean-square continuous, stationary random field $X$ on $\R^d$ with spectral density $f$ satisfies $X(s)=\int e^{is\cdot\omega}Z_X(d\omega)$ almost surely for each $s$, with an orthogonal random measure $Z_X$ such that $\E\abs{Z_X(d\omega)}^2=f(\omega)\,d\omega$ and $Z_X(-d\omega)=\overline{Z_X(d\omega)}$, and $\phi\mapsto\int\phi\,dZ_X$ is a linear isometry from $L^2(f\,d\omega)$ into $L^2(\Prob)$. Consequently, for bounded, measurable and Hermitian $\alpha$, $\int e^{is\cdot\omega}\alpha(\omega)Z_X(d\omega)$ is a real stationary field with spectral density $\abs\alpha^2f$, jointly stationary with $X$ with cross-spectral density $\alpha f$; and for an integrable kernel $\psi$, $(X\ast\psi)(s)=\int e^{is\cdot\omega}\hat\psi(\omega)Z_X(d\omega)$ almost surely, with $\hat\psi(\omega)=\int e^{-iu\cdot\omega}\psi(u)\,du$, by the isometry and Fubini's theorem.
\end{enumerate}

\begin{table}[H]
\centering\small
\caption{Map of the proofs. ``(C$k$)'' refers to the classical results of \cref{app:classical}.}\label{tab:proofmap}
\begin{tabularx}{\textwidth}{>{\raggedright\arraybackslash}p{2.4cm}>{\raggedright\arraybackslash}p{2.35cm}X>{\raggedright\arraybackslash}p{3.0cm}}
\toprule
Result & Proof & Main tools & Rests on \\
\midrule
\cref{prop:normalform} & \cref{app:normalform} & functional calculus of $L$; Frisch--Waugh--Lovell; profiling a penalized fit & --- \\
\cref{thm:frontier} & \cref{app:frontier} & Kakutani's theorem (C1); Gaussian affinities; Le Cam's two-point bound; a Bayesian Cram\'er--Rao bound; Toeplitz and Kronecker lemmas; Kolmogorov's criterion (C2) & \cref{lem:distances,lem:obstruction,lem:gaussaff,lem:chisq,lem:toeplitz,lem:kolmogorov,lem:vantrees} \\
\cref{cor:matern} & \cref{app:matern} & spectral dimension; power sums; concentration of weighted $\chi^2$ sums & \cref{thm:frontier}, \cref{lem:chisq,lem:powersums} \\
\cref{prop:whatgls} & \cref{app:whatgls} & explicit spectra; AM--GM & \cref{thm:frontier} \\
\cref{prop:pointsamples} & \cref{app:pointsamples} & quadrature errors under H\"older paths; perturbation of ratios & \cref{thm:frontier}, \cref{lem:kernels,lem:obstruction} \\
\cref{thm:resolution}, \cref{cor:latticetarget} & \cref{app:resolution} & spectral representation (C8); periodization of spectra; explicit construction; Cauchy--Schwarz & \cref{lem:kernels} \\
\cref{thm:window}, \cref{prop:estimators} & \cref{app:estimators} & power sums; conditionally Gaussian ratios; Kolmogorov's criterion (C2) & \cref{prop:normalform}, \cref{lem:chisq,lem:powersums,lem:toeplitz,lem:condgauss,lem:inverse} \\
\cref{cor:areal} & \cref{app:areal} & rescaling of weights; eigenvalues of cycles and tori; summation by parts on graphs & \cref{thm:window}, \cref{lem:chisq,lem:powersums,lem:condgauss} \\
\cref{thm:minimax} & \cref{app:minimax} & inverse moments; two-point lower bounds; the affinity formula \eqref{eq:affinity} & \cref{lem:distances,lem:inverse,lem:powersums} \\
\cref{thm:shift}, \cref{prop:shifts} & \cref{app:shift} & kernels and affinities; Kakutani (C1); Cameron--Martin theorem (C5); Lusin--Souslin (C3) & \cref{lem:distances,lem:kernels,lem:obstruction,lem:chisq,lem:powersums,lem:location} \\
\cref{lem:regularity,lem:composition} & \cref{app:regularity} & spectral representation (C8); Peetre's inequality; zero--one law and Fernique's theorem (C5); Lusin--Souslin (C3); Sobolev composition (C6) & --- \\
\cref{thm:nonlinear} & \cref{app:nonlinear} & local inversion of smooth maps; Sobolev composition & \cref{thm:shift}, \cref{prop:shifts}(a), \cref{lem:distances,lem:regularity,lem:composition} \\
\cref{thm:channels}, \cref{cor:kriging} & \cref{app:channels} & kernels; moment bounds for counts; heat-kernel bounds; conditioning a Gaussian field on point values; Kakutani (C1); Kolmogorov (C2); Mercer (C7) & \cref{prop:pointsamples,thm:shift}, \cref{lem:gaussaff,lem:kernels,lem:obstruction,lem:counts,lem:conditioning} \\
\bottomrule
\end{tabularx}
\end{table}

\subsection{Auxiliary lemmas}\label{app:prelim}

For probability measures $P,Q$ with densities $p,q$ with respect to a dominating measure $\lambda$, we use the total variation distance $\TV(P,Q)=\sup_A\abs{P(A)-Q(A)}=\tfrac12\int\abs{p-q}\,d\lambda$, the Hellinger affinity $\aff(P,Q)=\int\sqrt{pq}\,d\lambda$ and the squared Hellinger distance $H^2(P,Q)=1-\aff(P,Q)=\tfrac12\int(\sqrt p-\sqrt q)^2d\lambda$; none depends on the choice of $\lambda$.

\begin{lemma}[Distances and two-point bounds]\label{lem:distances}
\begin{enumerate}[label=(\alph*),itemsep=1pt]
\item $H^2(P,Q)\le\TV(P,Q)$ and $1-\TV(P,Q)\ge\tfrac12\aff(P,Q)^2$.
\item $P\perpm Q$ if and only if $\aff(P,Q)=0$, if and only if $\TV(P,Q)=1$.
\item If $T$ is a real statistic and $b_0,b_1\in\R$ with $\abs{b_1-b_0}=\Delta>0$, then $P(\abs{T-b_0}\ge\Delta/2)+Q(\abs{T-b_1}\ge\Delta/2)\ge1-\TV(P,Q)$.
\end{enumerate}
\end{lemma}
\begin{proof}
(a) Since $\abs{\sqrt p-\sqrt q}\le\sqrt p+\sqrt q$, we have $H^2=\tfrac12\int(\sqrt p-\sqrt q)^2\le\tfrac12\int\abs{\sqrt p-\sqrt q}(\sqrt p+\sqrt q)=\tfrac12\int\abs{p-q}=\TV$. Next, $\min(p,q)=\tfrac12(p+q-\abs{p-q})$ gives $1-\TV=\int\min(p,q)$. By the Cauchy--Schwarz inequality, $\aff^2=\big(\int\sqrt{\min(p,q)\max(p,q)}\big)^2\le\int\min(p,q)\int\max(p,q)$, and $\int\max(p,q)\le\int(p+q)=2$.
(b) If $P(A)=1$ and $Q(A)=0$, then $p=0$ $\lambda$-a.e.\ off $A$ and $q=0$ $\lambda$-a.e.\ on $A$, so $pq=0$ a.e.\ and $\aff=0$. If $\aff=0$, then $pq=0$ a.e., and $A=\{p>0\}$ has $P(A)=1$, $Q(A)=0$. Finally $\TV=1$ if and only if $\int\min(p,q)=0$, that is, if and only if $pq=0$ a.e.
(c) The events $E_i=\{\abs{T-b_i}<\Delta/2\}$, $i=0,1$, are disjoint, so $E_1\subseteq E_0^c$ and $P(E_0^c)+Q(E_1^c)\ge P(E_1)+1-Q(E_1)\ge1-\TV(P,Q)$.
\end{proof}

\begin{lemma}[Markov kernels and equivalence]\label{lem:kernels}
Let $\cX,\mathcal Y,\mathcal W,\mathcal Z$ be standard Borel spaces.
\begin{enumerate}[label=(\alph*),itemsep=1pt]
\item Let $\mu$ be a probability measure on $\cX$ and $K,K'$ Markov kernels from $\cX$ to $\mathcal Y$ with $K(x,\cdot)\ll K'(x,\cdot)$ for $\mu$-a.e.\ $x$. Then $\mu\otimes K\ll\mu\otimes K'$. If $\mu'\equivm\mu$, then $\mu\otimes K'\equivm\mu'\otimes K'$.
\item If $P\ll P'$ on $\mathcal W$ and $R$ is a Markov kernel from $\mathcal W$ to $\mathcal Z$, then $PR\ll P'R$, where $PR(B)=\int R(w,B)\,P(dw)$. This covers image measures under a common measurable map ($R(w,\cdot)=\delta_{\Phi(w)}$) and data obtained by adding independent noise with a common law.
\item $\aff(\mu\otimes K,\mu\otimes K')=\int\aff(K(x,\cdot),K'(x,\cdot))\,\mu(dx)$, the integrand being measurable in $x$.
\item If $P_i\equivm Q_i$ on $\mathcal W_i$ for $i=1,2$, then $P_1\otimes P_2\equivm Q_1\otimes Q_2$.
\item If $\Phi:\mathcal W\to\mathcal Z$ is Borel and injective, then $P\equivm P'$ if and only if $P\circ\Phi^{-1}\equivm P'\circ\Phi^{-1}$, and $P\perpm P'$ if and only if $P\circ\Phi^{-1}\perpm P'\circ\Phi^{-1}$.
\end{enumerate}
\end{lemma}
\begin{proof}
(a) Let $A\subseteq\cX\times\mathcal Y$ be measurable with $(\mu\otimes K')(A)=\int K'(x,A_x)\,\mu(dx)=0$, where $A_x=\{y:(x,y)\in A\}$. Then $K'(x,A_x)=0$ for $\mu$-a.e.\ $x$, hence $K(x,A_x)=0$ for $\mu$-a.e.\ $x$, and $(\mu\otimes K)(A)=0$. If $\mu'\equivm\mu$, the set $\{x:K'(x,A_x)>0\}$ is $\mu$-null if and only if it is $\mu'$-null, so $(\mu\otimes K')(A)=0$ if and only if $(\mu'\otimes K')(A)=0$.
(b) If $P'R(B)=0$, then $R(\cdot,B)=0$ $P'$-a.e., hence $P$-a.e., and $PR(B)=0$. Adding noise $\varepsilon$ with law $\rho$ through a measurable map $G(w,e)$ is the kernel $R(w,B)=\rho\{e:G(w,e)\in B\}$.
(c) By (C4) there are jointly measurable densities $k,k'$ of $K(x,\cdot),K'(x,\cdot)$ with respect to $\lambda_x$. Then $\mu\otimes K$ and $\mu\otimes K'$ have densities $k$ and $k'$ with respect to the measure $\mu(dx)\lambda_x(dy)$, and Tonelli's theorem gives $\aff(\mu\otimes K,\mu\otimes K')=\int\big(\int\sqrt{kk'}\,d\lambda_x\big)\mu(dx)$, with a measurable inner integral.
(d) With $p_i=dP_i/dQ_i$, which is positive $Q_i$-a.e., the product $P_1\otimes P_2$ has the density $p_1(w_1)p_2(w_2)$ with respect to $Q_1\otimes Q_2$, and this density is positive $Q_1\otimes Q_2$-a.e. A measure with an a.e.\ positive density is equivalent to the reference measure.
(e) If $P\equivm P'$, the images are equivalent by (b). Conversely, let $A\subseteq\mathcal W$ be Borel with $P'(A)=0$. By (C3), $\Phi(A)$ is Borel, and $\Phi^{-1}(\Phi(A))=A$ by injectivity, so $P'\circ\Phi^{-1}(\Phi(A))=0$, hence $P\circ\Phi^{-1}(\Phi(A))=P(A)=0$; by symmetry $P\equivm P'$. If $P(A)=1$ and $P'(A)=0$, the Borel set $\Phi(A)$ separates the images; if $B$ separates the images, $\Phi^{-1}(B)$ separates $P$ and $P'$.
\end{proof}

\begin{lemma}[Equivalence obstructs consistency]\label{lem:obstruction}
Let $\cF_n$ be the $\sigma$-fields generated by the stage-$n$ data and $\cF_\infty=\bigvee_n\cF_n$. If $P_2\ll P_1$ on $\cF_\infty$ (in particular, if $P_1\equivm P_2$) and $\psi_1\ne\psi_2$, there are no $\cF_n$-measurable $T_n$ with $T_n\to\psi_i$ in $P_i$-probability for both $i=1,2$.
\end{lemma}
\begin{proof}
Suppose such $T_n$ exist. Convergence in probability implies almost sure convergence along a subsequence, so there is a subsequence $(n_k)$ with $T_{n_k}\to\psi_1$ $P_1$-a.s.; passing to a further subsequence, still denoted $(n_k)$, also $T_{n_k}\to\psi_2$ $P_2$-a.s. The event $E=\{T_{n_k}\to\psi_1\}$ belongs to $\cF_\infty$ and $P_1(E^c)=0$, so $P_2(E^c)=0$ and $P_2(E)=1$. On $E$ the sequence $T_{n_k}$ converges to $\psi_1\ne\psi_2$, which contradicts $P_2(T_{n_k}\to\psi_2)=1$.
\end{proof}

\begin{lemma}[Gaussian integrals]\label{lem:gaussaff}
Let $\phi_g$ denote the $N(0,g)$ density.
\begin{enumerate}[label=(\alph*),itemsep=1pt]
\item For $a,b\in\R$ and $g>0$, $\int\sqrt{\phi_g(y-a)\phi_g(y-b)}\,dy=\exp\{-(a-b)^2/(8g)\}$.
\item If $Z\sim N(0,f)$ and $s\ge0$, then $\E e^{-sZ^2}=(1+2sf)^{-1/2}$.
\item For $t\ge0$, $\tfrac12\min(t,1)\le\log(1+t)\le t$.
\end{enumerate}
\end{lemma}
\begin{proof}
(a) $\sqrt{\phi_g(y-a)\phi_g(y-b)}=(2\pi g)^{-1/2}\exp\{-[(y-a)^2+(y-b)^2]/(4g)\}$, and $(y-a)^2+(y-b)^2=2\big(y-\tfrac{a+b}2\big)^2+\tfrac12(a-b)^2$. The integral over $y$ of $\exp\{-(y-\frac{a+b}2)^2/(2g)\}$ equals $(2\pi g)^{1/2}$, which leaves $\exp\{-(a-b)^2/(8g)\}$.
(b) $\E e^{-sZ^2}=(2\pi f)^{-1/2}\int\exp\{-z^2(1+2sf)/(2f)\}\,dz=(1+2sf)^{-1/2}$.
(c) The function $\log(1+t)-t/2$ vanishes at $0$ and has derivative $1/(1+t)-1/2\ge0$ on $[0,1]$, so $\log(1+t)\ge t/2$ there; for $t\ge1$, $\log(1+t)\ge\log2>1/2$. The upper bound is the elementary inequality $\log(1+t)\le t$.
\end{proof}

\begin{lemma}[Weighted $\chi^2$ sums]\label{lem:chisq}
Let $\xi_1,\xi_2,\dots$ be independent standard normal.
\begin{enumerate}[label=(\alph*),itemsep=1pt]
\item Let $a_j\ge0$. If $\sum_ja_j<\infty$, then $\sum_ja_j\xi_j^2<\infty$ a.s.; if $\sum_ja_j=\infty$, then $\sum_ja_j\xi_j^2=\infty$ a.s.
\item For each $N$ let $a_{N,j}\ge0$ for $j$ in a finite set $J_N$, $A_N=\sum_ja_{N,j}$ and $S_N=\sum_ja_{N,j}\xi_j^2$. If $A_N\to\infty$, then $S_N\to\infty$ in probability. If $\sum_ja_{N,j}^2=o(A_N^2)$, then $S_N/A_N\to1$ in probability. If $\sup_NA_N<\infty$, then $(S_N)$ is bounded in probability.
\item In particular, $X\in\CM_V$ a.s.\ if $\cI<\infty$ and $X\notin\CM_V$ a.s.\ if $\cI=\infty$, because $\norm{X}^2_{\CM_V}=\sum_jq_j\xi_j^2$.
\end{enumerate}
\end{lemma}
\begin{proof}
(b) For $t>0$ and $s\ge0$, Markov's inequality and \cref{lem:gaussaff}(b) give
\[
\Prob(S_N\le t)=\Prob\big(e^{-sS_N}\ge e^{-st}\big)\le e^{st}\prod_j(1+2sa_{N,j})^{-1/2}\le e^{st}(1+2sA_N)^{-1/2},
\]
because $\prod_j(1+y_j)\ge1+\sum_jy_j$ for $y_j\ge0$. With $s=1/t$ the bound is $e\,(1+2A_N/t)^{-1/2}$, which tends to zero when $A_N\to\infty$. Since $\E S_N=A_N$ and $\Var S_N=2\sum_ja_{N,j}^2$, Chebyshev's inequality gives $\Prob(\abs{S_N/A_N-1}>\epsilon)\le2\sum_ja_{N,j}^2/(\epsilon^2A_N^2)$. Markov's inequality gives $\Prob(S_N>K)\le A_N/K$.
(a) The first claim follows from $\E\sum_ja_j\xi_j^2=\sum_ja_j<\infty$. For the second, apply (b) with $J_N=\{1,\dots,N\}$ and $a_{N,j}=a_j$: $S_N\to\infty$ in probability, and since $S_N$ is non-decreasing in $N$, its limit $S_\infty\ge S_N$ satisfies $\Prob(S_\infty\le K)\le\Prob(S_N\le K)\to0$ for every $K$.
(c) $\norm{X}^2_{\CM_V}=\sum_jx_j^2/g_j=\sum_jq_j\xi_j^2$; apply (a).
\end{proof}

\begin{lemma}[Power sums]\label{lem:powersums}
Let $e\in\R$. Then $\sum_{j\le N}j^e\asymp N^{e+1}$ if $e>-1$, $\sum_{j\le N}j^e\asymp\log N$ if $e=-1$, and $\sum_jj^e<\infty$ if $e<-1$. If $e\ge-1$, then $\sum_{j\le N}j^{2e}=o\big((\sum_{j\le N}j^e)^2\big)$. If $e<-1$ and $N\ge2K\ge2$, then $\sum_{K<j\le N}j^e\asymp K^{e+1}$. The same statements hold with $j^e$ replaced by any $a_j>0$ with $a_j\asymp j^e$.
\end{lemma}
\begin{proof}
Compare with integrals: for $e\ge0$, $\int_0^Nx^edx\le\sum_{j\le N}j^e\le\int_1^{N+1}x^edx$; for $e<0$, $\int_1^{N+1}x^edx\le\sum_{j\le N}j^e\le1+\int_1^Nx^edx$. This gives the three cases. For the second claim: if $e>-1/2$, $\sum_{j\le N}j^{2e}\asymp N^{2e+1}=o(N^{2e+2})$; if $-1<e\le-1/2$, $\sum_{j\le N}j^{2e}=O(\log N)=o(N^{2e+2})$; if $e=-1$, $\sum_jj^{-2}<\infty$ while $(\log N)^2\to\infty$. For the window, $\int_{K+1}^{2K+1}x^edx\le\sum_{K<j\le N}j^e\le\int_K^\infty x^edx=K^{e+1}/\abs{e+1}$, and the left side is at least $(2K+1)^e K\ge3^eK^{e+1}$. If $a_j\asymp j^e$, the sums of $a_j$ and of $j^e$ are within constant factors of each other up to finitely many terms, which do not affect divergent sums or the window bound for large $K$.
\end{proof}

\begin{lemma}[Toeplitz and Kronecker]\label{lem:toeplitz}
\begin{enumerate}[label=(\alph*),itemsep=1pt]
\item If $w_j\ge0$, $W_N=\sum_{j\le N}w_j\to\infty$ and $\alpha_j\to\alpha_\infty$, then $\sum_{j\le N}w_j\alpha_j/W_N\to\alpha_\infty$.
\item If $0<b_1\le b_2\le\cdots\to\infty$ and $\sum_jz_j/b_j$ converges, then $\sum_{j\le N}z_j/b_N\to0$.
\end{enumerate}
\end{lemma}
\begin{proof}
(a) Given $\epsilon>0$ choose $j_0$ with $\abs{\alpha_j-\alpha_\infty}\le\epsilon$ for $j>j_0$. Then $\abs{\sum_{j\le N}w_j\alpha_j/W_N-\alpha_\infty}\le\sum_{j\le j_0}w_j\abs{\alpha_j-\alpha_\infty}/W_N+\epsilon$, and the first term tends to zero.
(b) Let $r_j=\sum_{k>j}z_k/b_k$ for $j\ge0$, so that $r_j\to0$ and $z_j=b_j(r_{j-1}-r_j)$. Summation by parts gives $\sum_{j\le N}z_j=b_1r_0-b_Nr_N+\sum_{j=1}^{N-1}(b_{j+1}-b_j)r_j$. After division by $b_N$, the first two terms tend to zero, and so does the third by (a), applied with the weights $b_{j+1}-b_j\ge0$, whose partial sums $b_N-b_1$ diverge, and with $r_j\to0$.
\end{proof}

\begin{lemma}[A strong law for weighted sums]\label{lem:kolmogorov}
Let $\varepsilon_1,\varepsilon_2,\dots$ be independent with mean zero and variances $\sigma_j^2$, and let $c_j\ge0$ with $c_1>0$, $b_N=\sum_{j\le N}c_j\to\infty$ and $\sigma_j^2\le C\,c_j$ for all $j$. Then $\sum_{j\le N}\varepsilon_j/b_N\to0$ almost surely.
\end{lemma}
\begin{proof}
Since $c_j=b_j-b_{j-1}$ for $j\ge2$,
\[
\sum_j\Var(\varepsilon_j/b_j)\le C\sum_j\frac{c_j}{b_j^2}\le C\Big(\frac1{b_1}+\sum_{j\ge2}\frac{b_j-b_{j-1}}{b_{j-1}b_j}\Big)=C\Big(\frac1{b_1}+\sum_{j\ge2}\Big(\frac1{b_{j-1}}-\frac1{b_j}\Big)\Big)\le\frac{2C}{b_1}.
\]
By (C2) the series $\sum_j\varepsilon_j/b_j$ converges almost surely, and \cref{lem:toeplitz}(b) gives the claim.
\end{proof}

\begin{lemma}[Conditionally Gaussian ratios]\label{lem:condgauss}
Let $\mathcal G$ be a $\sigma$-field, let $D_N>0$ and $V_N\ge0$ be $\mathcal G$-measurable, and let $M_N$ have conditional law $N(0,V_N)$ given $\mathcal G$.
\begin{enumerate}[label=(\alph*),itemsep=1pt]
\item $M_N/D_N\to0$ in probability if and only if $V_N/D_N^2\to0$ in probability. More precisely, if $\Prob(V_N/D_N^2\ge\eta)\ge\eta$ for some $\eta>0$ and infinitely many $N$, then for these $N$ and every $\epsilon>0$, $\Prob(\abs{M_N/D_N}>\epsilon)\ge\eta\,\Prob(\abs Z>\epsilon/\sqrt\eta)$ with $Z\sim N(0,1)$.
\item Let $S_N$ be $\mathcal G$-measurable, and suppose $V_N/D_N^2\ge R_0$ for all $N$, where $R_0$ is $\mathcal G$-measurable, does not depend on $N$ and is positive almost surely. Then $S_N+M_N/D_N$ does not converge in probability to any constant.
\end{enumerate}
\end{lemma}
\begin{proof}
(a) Let $R_N=V_N/D_N^2$ and $\psi(x)=\Prob(\abs Z>x)$ for $x\in[0,\infty]$, with $\psi(\infty)=0$. Given $\mathcal G$, $M_N/D_N\sim N(0,R_N)$, so $\Prob(\abs{M_N/D_N}>\epsilon\mid\mathcal G)=\psi(\epsilon/\sqrt{R_N})$. If $R_N\to0$ in probability, then $\psi(\epsilon/\sqrt{R_N})\to0$ in probability, and since $\psi\le1$, also in $L^1$, so $\Prob(\abs{M_N/D_N}>\epsilon)\to0$. Since $\psi$ is non-increasing, $\Prob(\abs{M_N/D_N}>\epsilon)\ge\E[\psi(\epsilon/\sqrt{R_N});R_N\ge\eta]\ge\eta\,\psi(\epsilon/\sqrt\eta)$ whenever $\Prob(R_N\ge\eta)\ge\eta$.
(b) Given $\mathcal G$, $T_N=S_N+M_N/D_N$ is normal with variance $R_N\ge R_0$, so its conditional density is bounded by $(2\pi R_0)^{-1/2}$. For $b\in\R$ and $\epsilon>0$, $\Prob(\abs{T_N-b}\le\epsilon\mid\mathcal G)\le\min\{1,2\epsilon(2\pi R_0)^{-1/2}\}$, hence $\Prob(\abs{T_N-b}\le\epsilon)\le\pi(\epsilon)=\E\min\{1,2\epsilon(2\pi R_0)^{-1/2}\}$, which does not depend on $N$ or $b$ and tends to zero as $\epsilon\downarrow0$ by dominated convergence. Choosing $\epsilon$ with $\pi(\epsilon)<1$ shows that $T_N$ does not converge in probability to $b$.
\end{proof}

\begin{lemma}[A Bayesian Cram\'er--Rao bound; a version of the van Trees inequality \citep{GillLevit1995}]\label{lem:vantrees}
Let $\{P_\beta:\beta\in\R\}$ be probability measures on a standard Borel space $\mathcal D$, and let $\mathsf x:\mathcal D\to\cX$ be measurable with a law $P^{\mathsf x}$ that does not depend on $\beta$. Suppose there are measurable $S:\mathcal D\to\R$ and $Q:\cX\to(0,\infty)$ such that, for $P^{\mathsf x}$-a.e.\ $x$, the conditional laws $P^x_\beta=P_\beta(\cdot\mid\mathsf x=x)$ satisfy $dP^x_\beta/dP^x_0=\exp\{\beta S-\beta^2Q(x)/2\}$ and, under $P^x_\beta$, $S-\beta Q(x)$ has mean $0$ and variance $Q(x)$. Let $\pi$ be an absolutely continuous probability density on $\R$ with $J(\pi)=\int(\pi')^2/\pi<\infty$ (with $0/0=0$). Then every measurable $T:\mathcal D\to\R$ satisfies
\[
\int\E_\beta(T-\beta)^2\,\pi(\beta)\,d\beta\;\ge\;\int\frac{P^{\mathsf x}(dx)}{Q(x)+J(\pi)}\;\ge\;\frac1{\int Q\,dP^{\mathsf x}+J(\pi)} .
\]
\end{lemma}
\begin{proof}
\emph{Step 1 (properties of $\pi$).} Since $\pi\ge0$ is absolutely continuous, $\pi$ attains its minimum at every point of $\{\pi=0\}$ where it is differentiable, so $\pi'=0$ there; hence $\pi'=0$ a.e.\ on $\{\pi=0\}$, and the Cauchy--Schwarz inequality gives $\int\abs{\pi'}=\int_{\pi>0}(\abs{\pi'}/\sqrt\pi)\sqrt\pi\le J(\pi)^{1/2}$. Hence $\pi$ has limits at $\pm\infty$, which are zero because $\pi$ is integrable, and $\pi(\beta)=\int_{-\infty}^\beta\pi'\le J(\pi)^{1/2}$ is bounded.

\emph{Step 2 (integration by parts in $\beta$).} Fix $x$ with $Q=Q(x)>0$ and write $\E^x_\beta$ for expectation under $P^x_\beta$. Assume $r(x)=\int\E^x_\beta(T-\beta)^2\pi(\beta)\,d\beta<\infty$; otherwise there is nothing to prove for this $x$. For each data point $d$ let $\ell_d(\beta)=e^{\beta S(d)-\beta^2Q/2}\pi(\beta)$. It is absolutely continuous on bounded intervals with derivative $\ell_d'=(S-\beta Q)\ell_d+e^{\beta S-\beta^2Q/2}\pi'$. Since $e^{\beta S-\beta^2Q/2}\le e^{-\beta^2Q/4}$ for $\abs\beta\ge4\abs S/Q$, and since $\pi$ is bounded and $\pi'$ integrable, the functions $\ell_d$ and $\beta\ell_d$ tend to zero at $\pm\infty$, and $\ell_d'$ and $\beta\ell_d'$ are integrable. Integration by parts on $[-A,A]$ and $A\to\infty$ give $\int\ell_d'\,d\beta=0$ and $\int\beta\ell_d'\,d\beta=-\int\ell_d\,d\beta$, so
\begin{equation}\label{eq:vtibp}
\int\big(T(d)-\beta\big)\,\ell_d'(\beta)\,d\beta=\int\ell_d(\beta)\,d\beta .
\end{equation}

\emph{Step 3 (integration over the data).} For non-negative measurable $G$, Tonelli's theorem and the density hypothesis give $\iint G(d,\beta)\,\ell_d(\beta)\,d\beta\,P^x_0(dd)=\int\pi(\beta)\,\E^x_\beta[G(\cdot,\beta)]\,d\beta$. With $G=1$ the right side of \eqref{eq:vtibp} integrates to $\int\pi=1$. Write $\ell_d'=s_d\ell_d$ with $s_d(\beta)=S-\beta Q+\pi'(\beta)/\pi(\beta)$ on $\{\pi>0\}$; on $\{\pi=0\}$ both $\ell_d$ and $\ell'_d$ vanish a.e. Then
\[
\iint s_d^2\,\ell_d\,d\beta\,dP^x_0=\int\pi\,\E^x_\beta(S-\beta Q)^2\,d\beta+2\int\pi'(\beta)\,\E^x_\beta(S-\beta Q)\,d\beta+J(\pi)=Q+J(\pi),
\]
where Fubini's theorem applies to the cross term because $\int\abs{\pi'(\beta)}\,\E^x_\beta\abs{S-\beta Q}\,d\beta\le Q^{1/2}\int\abs{\pi'}<\infty$,
and $\iint(T-\beta)^2\ell_d\,d\beta\,dP^x_0=r(x)$. By the Cauchy--Schwarz inequality $\iint\abs{T-\beta}\,\abs{\ell'_d}\,d\beta\,dP^x_0\le\{r(x)(Q+J(\pi))\}^{1/2}<\infty$, so Fubini's theorem applies to the left side of \eqref{eq:vtibp}, and
\[
1=\iint(T-\beta)\,s_d\,\ell_d\,d\beta\,dP^x_0\le\{r(x)\,(Q+J(\pi))\}^{1/2}.
\]
Hence $r(x)\ge1/(Q(x)+J(\pi))$.

\emph{Step 4 (integration over $\mathsf x$).} Since the law of $\mathsf x$ does not depend on $\beta$, $\int\E_\beta(T-\beta)^2\pi(\beta)d\beta=\int r(x)\,P^{\mathsf x}(dx)$. The first inequality follows, and the second is Jensen's inequality for the convex function $t\mapsto1/(t+J(\pi))$.
\end{proof}

\begin{lemma}[Inverse moments]\label{lem:inverse}
Let $D_N=\sum_{j\le N}q_j\xi_j^2$ with $q_j\asymp j^{p-1}$ and $p>0$. Then $\sum_{j\le N}q_j\asymp N^p$, and for each $k\ge1$ there are $C_k,N_k$ with $\E D_N^{-k}\le C_kN^{-pk}$ for $N\ge N_k$.
\end{lemma}
\begin{proof}
The first claim is \cref{lem:powersums}. Let $J_N=\{\lceil N/2\rceil,\dots,N\}$, with $m=\abs{J_N}\ge N/2$ elements. For $j\in J_N$, $q_j\ge c\,j^{p-1}\ge c\min(1,2^{1-p})N^{p-1}=:c'N^{p-1}$, so $D_N\ge c'N^{p-1}\chi^2_m$, where $\chi^2_m=\sum_{j\in J_N}\xi_j^2$. For $m>2k$, $\E(\chi^2_m)^{-k}=2^{-k}\Gamma(m/2-k)/\Gamma(m/2)$, and $\Gamma(m/2)/\Gamma(m/2-k)=\prod_{i=1}^k(m/2-i)\ge(m/4)^k$ when $m\ge4k$. Hence $\E(\chi^2_m)^{-k}\le(2/m)^k\le(4/N)^k$ and $\E D_N^{-k}\le(c')^{-k}4^kN^{-pk}$ for $N\ge8k$.
\end{proof}

\begin{lemma}[Logarithms of counts]\label{lem:counts}
There is a universal constant $C$ with the following property. Let $N$ take values in $\{0,1,2,\dots\}$, with mean $\kappa\ge1$, $\E(N-\kappa)^2\le\kappa$, $\E(N-\kappa)^4\le\kappa+3\kappa^2$ and $\Prob(N\le\kappa/2)\le e^{-\kappa/7}$. Then $Z=\log\{(N+\tfrac12)/\kappa\}$ satisfies $\abs{\E Z}\le C/\kappa$ and $\E Z^2\le C/\kappa$. The hypotheses hold for the Poisson law with mean $\kappa$ and for the binomial law with $m$ trials, success probability $p$ and $\kappa=mp$.
\end{lemma}
\begin{proof}
Let $D=(N+\frac12-\kappa)/\kappa$, so $Z=\log(1+D)$. Then $\E D=1/(2\kappa)$, $\E D^2=\{\E(N-\kappa)^2+\frac14\}/\kappa^2\le2/\kappa$ and, by $(a+b)^4\le8(a^4+b^4)$, $\E D^4\le8\{\E(N-\kappa)^4+1\}/\kappa^4\le40/\kappa^2$. Let $A=\{N>\kappa/2\}$. On $A$, $D>-1/2$, and Taylor's theorem gives $\abs{\log(1+x)-x}=x^2/\{2(1+\theta x)^2\}\le2x^2$ for $x\ge-1/2$ and some $\theta\in(0,1)$. On $A^c$, $(N+\frac12)/\kappa\in[1/(2\kappa),1]$, so $\abs Z\le\log(2\kappa)$. Writing $\E Z=\E D-\E[D;A^c]+\E[Z-D;A]+\E[Z;A^c]$,
\[
\abs{\E Z}\le\frac1{2\kappa}+\Big(\frac2\kappa\Big)^{1/2}e^{-\kappa/14}+\frac4\kappa+\log(2\kappa)\,e^{-\kappa/7}\le\frac C\kappa,
\]
because $\kappa^{1/2}e^{-\kappa/14}$ and $\kappa\log(2\kappa)e^{-\kappa/7}$ are bounded on $[1,\infty)$. Similarly $\E Z^2\le2\E D^2+2\E[(Z-D)^2;A]+\E[Z^2;A^c]\le4/\kappa+320/\kappa^2+\log^2(2\kappa)e^{-\kappa/7}\le C/\kappa$.
For the Poisson law, $\E(N-\kappa)^2=\kappa$ and $\E(N-\kappa)^4=\kappa+3\kappa^2$. For the binomial law, $\E(N-\kappa)^2=mp(1-p)\le\kappa$ and $\E(N-\kappa)^4=mp(1-p)\{1+3(m-2)p(1-p)\}\le\kappa+3\kappa^2$. In both cases $\E e^{-sN}\le\exp\{\kappa(e^{-s}-1)\}$ for $s>0$ (for the binomial law because $1+y\le e^y$), so Markov's inequality with $s=\log2$ gives $\Prob(N\le\kappa/2)\le\exp\{s\kappa/2+\kappa(e^{-s}-1)\}=\exp\{-\kappa(1-\log2)/2\}\le e^{-\kappa/7}$.
\end{proof}

\begin{lemma}[Location families]\label{lem:location}
Let $p$ be an absolutely continuous probability density on $\R$ with $J(p)=\int(p')^2/p<\infty$.
\begin{enumerate}[label=(\alph*),itemsep=1pt]
\item $\sqrt p$ is absolutely continuous, with derivative $p'/(2\sqrt p)$ on $\{p>0\}$ and $0$ elsewhere, and $H^2(p,p(\cdot-t))\le t^2J(p)/8$ for all $t\in\R$.
\item If $\xi$ has density $p$ and $\E\xi^2<\infty$, then $\Var\xi\ge1/J(p)$.
\item The map $t\mapsto H^2(t)=H^2(p,p(\cdot-t))$ is continuous, positive for $t\ne0$ and tends to $1$ as $\abs t\to\infty$, and $H^2(t)/t^2\to J(p)/8>0$ as $t\to0$. Consequently $\inf_{t\ne0}H^2(t)/\min(t^2,1)>0$.
\end{enumerate}
\end{lemma}
\begin{proof}
(a) Since $p\ge0$ is absolutely continuous, $p$ attains its minimum at every point of $\{p=0\}$ where it is differentiable, so $p'=0$ a.e.\ on $\{p=0\}$, and the Cauchy--Schwarz inequality gives $\int\abs{p'}=\int_{p>0}(\abs{p'}/\sqrt p)\sqrt p\le J(p)^{1/2}$; in particular $p'$ is integrable. For $\epsilon>0$ the function $\sqrt{p+\epsilon}$ is absolutely continuous on bounded intervals, with derivative $p'/(2\sqrt{p+\epsilon})$, because $y\mapsto\sqrt y$ is continuously differentiable and Lipschitz on $[\epsilon,\infty)$. As $\epsilon\downarrow0$ these derivatives converge a.e.\ to $u'=p'/(2\sqrt p)\ind\{p>0\}$ and are dominated by $\abs{u'}$, which is square integrable with $\int(u')^2=J(p)/4$ and hence integrable on bounded intervals. Dominated convergence in $\sqrt{p+\epsilon}(b)-\sqrt{p+\epsilon}(a)=\int_a^bp'/(2\sqrt{p+\epsilon})$ shows that $u=\sqrt p$ is absolutely continuous with derivative $u'$. For $t>0$, $u(x)-u(x-t)=\int_{x-t}^xu'(v)\,dv$, so the Cauchy--Schwarz inequality gives $(u(x)-u(x-t))^2\le t\int_{x-t}^xu'(v)^2dv$; integrating in $x$ and using Tonelli's theorem, $\int(u(x)-u(x-t))^2dx\le t^2\int(u')^2=t^2J(p)/4$, and $H^2=\frac12\int(u-u(\cdot-t))^2$. The case $t<0$ is symmetric.
(b) By the Cauchy--Schwarz inequality $\int\abs x\abs{p'(x)}dx\le(\int x^2p)^{1/2}J(p)^{1/2}<\infty$. Because $p'$ is integrable, $p$ has limits at $\pm\infty$, both zero, and $\int p'=0$. Because $\int\abs xp(x)dx<\infty$, there are $A_k\to\infty$ with $A_k\{p(A_k)+p(-A_k)\}\to0$; otherwise $p(x)+p(-x)\ge c/x$ for all large $x$, which is not integrable. Integrating by parts on $[-A_k,A_k]$ and letting $k\to\infty$ gives $\int xp'(x)dx=-1$. With $m=\E\xi$ and $\int p'=0$, $1=-\int(x-m)p'(x)dx=-\E[(\xi-m)p'(\xi)/p(\xi)]\le(\Var\xi)^{1/2}J(p)^{1/2}$.
(c) Let $u=\sqrt p$ and $\tau_tu=u(\cdot-t)$, so that $H^2(t)=\frac12\norm{u-\tau_tu}_2^2$ and $\norm u_2=1$. Continuity follows from $\abs{\norm{u-\tau_tu}_2-\norm{u-\tau_su}_2}\le\norm{u-\tau_{s-t}u}_2$ and the continuity of translation in $L^2$. If $H^2(t)=0$ with $t\ne0$, then $p=p(\cdot-t)$ a.e., and $\int_\R p=\sum_{k\in\Z}\int_{kt}^{(k+1)t}p$ would be $0$ or infinite, which is impossible. As $\abs t\to\infty$, $\inner u{\tau_tu}\to0$: approximate $u$ within $\epsilon$ in $L^2$ by a compactly supported $v$; then $\inner v{\tau_tv}=0$ for large $\abs t$ and, since $\norm v_2\le1+\epsilon$, $\abs{\inner u{\tau_tu}-\inner v{\tau_tv}}\le\epsilon(2+\epsilon)$. Hence $H^2(t)=1-\inner u{\tau_tu}\to1$. Finally $(u-\tau_tu)/t=\int_0^1\tau_{\theta t}u'\,d\theta$, which converges to $u'$ in $L^2$ by Minkowski's integral inequality and the continuity of translation; hence $H^2(t)/t^2\to\frac12\norm{u'}_2^2=J(p)/8$. Here $J(p)>0$, because $J(p)=0$ would force $p'=0$ a.e.\ and $p$ constant. The function $H^2(t)/\min(t^2,1)$ is therefore continuous and positive on $t\ne0$, with the positive limits $J(p)/8$ at $0$ and $1$ at $\pm\infty$, so its infimum is positive.
\end{proof}

\subsection{Proof of \texorpdfstring{\cref{prop:normalform}}{Proposition 1}}\label{app:normalform}

\emph{Idea.} Every model and estimator in \cref{tab:zoo} is built from covariance, precision or smoothing operators that are functions of $L$, and functions of $L$ act diagonally on the eigencoordinates. What remains is bookkeeping.

\emph{Step 1 (functional calculus).} For a bounded Borel function $h$ on the spectrum, $h(L)u=\sum_jh(\lambda_j)\inner u{\varphi_j}\varphi_j$, and $h(L)^{-1}=(1/h)(L)$ when $h>0$. A centred Gaussian field with covariance operator $h(L)$ therefore has independent coefficients with variances $h(\lambda_j)$, and a Gaussian field with precision operator $p(L)$, $p>0$, has covariance operator $(1/p)(L)$.

\emph{Step 2 (upper block of \cref{tab:zoo}).} \emph{ICAR}: the density of the intrinsic CAR model \citep{Besag1974,BYM1991} is proportional to $\exp\{-\tfrac\tau2u^\top(D-W)u\}$, so the precision is $\tau L$ with $L=D-W$; on a connected graph the null space of $L$ is spanned by $\mathbf 1$, and on $\mathbf1^\perp$ the covariance has eigenvalues $(\tau\lambda_j)^{-1}$. \emph{Leroux CAR}: the precision is $\tau\{(1-\rho)I+\rho L\}$ \citep{Leroux2000}, with covariance eigenvalues $\tau^{-1}(1-\rho+\rho\lambda_j)^{-1}$. \emph{Proper CAR}: the precision is $\tau(D-\rho W)$; on a $k$-regular graph $D=kI$ and $W=kI-L$, so $D-\rho W=k(1-\rho)I+\rho L$. On an irregular graph $D$ and $W$ need not commute, $D-\rho W$ is then not a function of $L$, and the row is only approximate. \emph{SAR}: $(I-\rho W)Y=\varepsilon$ with $\varepsilon\sim N(0,\sigma^2I)$ and symmetric $W$ gives $\Cov(Y)=\sigma^2(I-\rho W)^{-2}$, and on a $k$-regular graph $I-\rho W=(1-\rho k)I+\rho L$. \emph{Whittle--Mat\'ern}: with $a=\nu+d_s/2$, $(\kappa^2+L)^{a/2}Z=\tau\mathcal W$ gives $z_j=\tau(\kappa^2+\lambda_j)^{-a/2}w_j$ with independent standard normal $w_j=\inner{\mathcal W}{\varphi_j}$, hence the spectrum $\tau^2(\kappa^2+\lambda_j)^{-a}$; kriging with this covariance is GLS and best linear unbiased prediction under this prior. \emph{Graph Mat\'ern GP}: the covariance $\sigma^2(2\nu/\kappa^2+L)^{-\nu}$ is the definition of \citet{Borovitskiy2021}, up to normalization. \emph{Moran eigenvector filtering}: on a $k$-regular graph $W\mathbf1=k\mathbf1$, so $W$ commutes with $M=I-\mathbf1\mathbf1^\top/n$; the eigenvectors of $MWM$ are $\mathbf 1$ (eigenvalue $0$) and the eigenvectors of $W$ orthogonal to $\mathbf1$, which are those of $L=kI-W$, with Moran eigenvalues $k-\lambda_j$. Fixed effects on the $K$ eigenvectors with the largest Moran eigenvalues, the $K$ smoothest modes, amount to an improper flat prior on their coefficients. \emph{Random-effects ESF and (M)GWR-type SVC models}: the stated priors on eigenvector coefficients \citep{MurakamiGriffith2015,Murakami2017} are spectral priors by definition. \emph{Linear GNN}: the simplified graph convolution computes $p(L)X$ with $p(\lambda)=(1-\lambda)^K$, $L$ being the normalized Laplacian of the graph with self-loops \citep{Wu2019SGC}.

\emph{Step 3 (the identity \eqref{eq:normalform}).} Substituting $y_j=\gamma_jx_j+v_j=(\beta+\alpha_j)x_j+v_j$,
\[
\sum_j\omega_jx_jy_j=\beta\sum_j\omega_jx_j^2+\sum_j\omega_j\alpha_jx_j^2+\sum_j\omega_jx_jv_j .
\]
Dividing by $\sum_k\tilde\omega_kx_k^2$ gives $\hat\beta=\beta\sum_j\pi_j+\sum_j\pi_j\alpha_j+\sum_j\omega_jx_jv_j/\sum_k\tilde\omega_kx_k^2$, which is \eqref{eq:normalform}.

\emph{Step 4 (lower block of \cref{tab:zoo}).} \emph{OLS}: by Parseval, $\inner XY/\inner XX=\sum_jx_jy_j/\sum_jx_j^2$. \emph{GLS} with covariance $\Sigma=h(L)+\sigma^2I$: $\hat\beta=\inner X{\Sigma^{-1}Y}/\inner X{\Sigma^{-1}X}$ and $\Sigma^{-1}=(h+\sigma^2)^{-1}(L)$, giving the weights $1/(h(\lambda_j)+\sigma^2)$; the same expression is the posterior mean of $\beta$ under a flat prior in the Gaussian mixed model. \emph{RSR}: the restricted model is $Y=X\beta+P^\perp\theta+\varepsilon$ with $P^\perp=I-X(X^\top X)^{-1}X^\top$, $\theta\sim N(0,\Sigma_\theta)$ and $\varepsilon\sim N(0,\sigma^2I)$, so $Y\sim N(X\beta,\Sigma_R)$ with $\Sigma_R=\sigma^2I+P^\perp\Sigma_\theta P^\perp$. Since $P^\perp X=0$, $\Sigma_RX=\sigma^2X$, hence $\Sigma_R^{-1}X=\sigma^{-2}X$, and the GLS estimate (the posterior mean under a flat prior) is $X^\top\Sigma_R^{-1}Y/X^\top\Sigma_R^{-1}X=X^\top Y/X^\top X$, the OLS estimate \citep{Hanks2015,KhanCalder2022}. \emph{Eigenvector filtering}: by the Frisch--Waugh--Lovell theorem the coefficient of $X$ in the regression of $Y$ on $(X,\varphi_1,\dots,\varphi_K)$ equals the OLS coefficient of $(I-P_K)Y$ on $(I-P_K)X$, where $P_K$ projects onto $\operatorname{span}\{\varphi_1,\dots,\varphi_K\}$; in eigencoordinates $I-P_K$ sets the first $K$ coefficients to zero, so $\hat\beta=\sum_{j>K}x_jy_j/\sum_{j>K}x_j^2$. \emph{Spatial+}: the first stage gives $r=(I-S)X$ for a symmetric smoother $S=s(L)$. The second stage minimizes $\norm{Y-\beta r-g}^2+\lambda g^\top Pg$ over $(\beta,g)$, with a penalty $P=p(L)\succeq0$. For fixed $\beta$ the minimizer is $g=T(Y-\beta r)$ with $T=(I+\lambda P)^{-1}=t(L)$. Since $I-T=\lambda PT$ and $P$ and $T$ commute, the profiled objective is $(Y-\beta r)^\top\{(I-T)^2+(I-T)T\}(Y-\beta r)=(Y-\beta r)^\top(I-T)(Y-\beta r)$, whose minimizer is $\hat\beta=r^\top(I-T)Y/r^\top(I-T)r$; in eigencoordinates, $r_j=(1-s_j)x_j$ gives the weights $\omega_j=(1-t_j)(1-s_j)$ and $\tilde\omega_j=(1-t_j)(1-s_j)^2$. \emph{gSEM}: regressing the residual $(I-T)Y$ of the response on the residual $r=(I-S)X$ of the covariate, with smoothers $S=s(L)$ and $T=t(L)$, gives $\hat\beta=r^\top(I-T)Y/r^\top r$, that is, the weights $\omega_j=(1-s_j)(1-t_j)$ and $\tilde\omega_j=(1-s_j)^2$. \emph{High-pass GLS}: GLS on the modes $j>K$, the idealized spectral adjustment of \citet{Guan2023}, has weights $\ind\{j>K\}/(h(\lambda_j)+\sigma^2)$. \emph{GNN features}: regressing $Y$ on $p(L)X$ gives $\inner{p(L)X}{Y}/\norm{p(L)X}^2=\sum_jp_jx_jy_j/\sum_jp_j^2x_j^2$ with $p_j=p(\lambda_j)$.

\emph{Step 5 (part (c)).} If $\omega=\tilde\omega\ge0$, then $\pi_j\ge0$ and $\sum_j\pi_j=1$, so the second term of \eqref{eq:normalform} vanishes and $\hat\beta=\sum_j\pi_j\gamma_j+\sum_j\omega_jx_jv_j/\sum_k\omega_kx_k^2$. \qed

\subsection{Proof of \texorpdfstring{\cref{thm:frontier}}{Theorem 1}}\label{app:frontier}

\emph{Idea.} The complete record is a product over scales. Changing $\beta$ shifts the conditional mean of $y_j$ by an amount proportional to $x_j$, and the Hellinger affinity between the two laws at scale $j$ is explicit. Kakutani's theorem turns the summability of the scale-wise discrepancies, $\cI<\infty$, into equivalence of the laws, which no estimator can overcome; quantitative versions follow from the product formula for affinities and from a Bayesian Cram\'er--Rao bound. Conversely, when $\cI=\infty$, GLS is a weighted average of the scale-specific associations plus a noise term whose variance grows no faster than the total weight, and both pieces behave as required.

\emph{Step 0 (the law of the complete record).} The pairs $(x_j,u_j)$ are jointly Gaussian and uncorrelated across $j$, hence independent across $j$. By \eqref{eq:decomp}, $v_j=u_j-\alpha_jx_j$ is uncorrelated with, hence independent of, $x_j$. Therefore, under $P_{\beta,\alpha}$, the pairs $(x_j,y_j)$ are independent across $j$, with $x_j\sim N(0,f_j)$ and $y_j\mid x_j\sim N(\gamma_jx_j,g_j)$, and $P_{\beta,\alpha}=\bigotimes_jP_j$ with $P_j$ the law of $(x_j,y_j)$. All factors have positive densities on $\R^2$, so the factors of two parameter values are pairwise equivalent.

\emph{Step 1 (scale-wise affinities).} For two parameters with the same $(f,g)$, associations $\gamma,\gamma'$ and $\delta_j=\gamma'_j-\gamma_j$, \cref{lem:gaussaff}(a) and then \cref{lem:gaussaff}(b) with $s=\delta_j^2/(8g_j)$ give
\begin{equation}\label{eq:affinity}
\aff_j=\int\phi_{f_j}(x)\int\sqrt{\phi_{g_j}(y-\gamma_jx)\phi_{g_j}(y-\gamma'_jx)}\,dy\,dx=\int\phi_{f_j}(x)e^{-\delta_j^2x^2/(8g_j)}dx=\Big(1+\tfrac14\delta_j^2q_j\Big)^{-1/2}.
\end{equation}

\emph{Step 2 (summability).} Let $t_j=\delta_j^2q_j/4$, so $H^2_j=1-(1+t_j)^{-1/2}$. If $\sum_jt_j<\infty$, then $H^2_j\le1-e^{-t_j/2}\le t_j/2$, because $(1+t)^{-1/2}=e^{-\frac12\log(1+t)}\ge e^{-t/2}$. Conversely, if $\sum_jH^2_j<\infty$, then $H^2_j\to0$, so $t_j\to0$ and $t_j\le1$ for large $j$; for $t\in[0,1]$, \cref{lem:gaussaff}(c) gives $(1+t)^{-1/2}\le e^{-t/4}$ and hence $H^2_j\ge1-e^{-t_j/4}\ge t_j/8$, because $1-e^{-x}\ge x/2$ for $0\le x\le1$. Thus $\sum_jH^2(P_j,P'_j)<\infty$ if and only if $\sum_j\delta_j^2q_j<\infty$, that is, if and only if $\delta\in\cN_q$.

\emph{Step 3 (part (i): equivalence).} Take $\alpha'=\alpha$ and $\beta'-\beta=\Delta\ne0$, so $\delta_j\equiv\Delta$ and $\sum_j\delta_j^2q_j=\Delta^2\cI<\infty$. By Step 2 and (C1), $P_{\beta,\alpha}\equivm P_{\beta',\alpha}$, and \cref{lem:obstruction} shows that no estimator is consistent at both parameters. Since $\alpha$ is common to both, this holds also when $\alpha=0$ is known.

\emph{Step 4 (part (i): two-point bound).} By (C1), \eqref{eq:affinity} and $\log(1+t)\le t$,
\[
\aff(P_{\beta_0,\alpha},P_{\beta_1,\alpha})=\prod_j\Big(1+\tfrac14\Delta^2q_j\Big)^{-1/2}\ge\prod_je^{-\Delta^2q_j/8}=e^{-\Delta^2\cI/8}.
\]
By \cref{lem:distances}(c) and (a), $P_0(\abs{T-\beta_0}\ge\Delta/2)+P_1(\abs{T-\beta_1}\ge\Delta/2)\ge1-\TV\ge\aff^2/2\ge\frac12e^{-\Delta^2\cI/4}$, and the larger of the two probabilities is at least half the sum. This is \eqref{eq:lecam}.

\emph{Step 5 (part (i): Bayesian bound).} Apply \cref{lem:vantrees} with $\mathsf x=(x_j)_j$, whose law does not depend on $\beta$, $Q=\sum_jx_j^2/g_j$ and $S=\sum_jx_j(y_j-\alpha_jx_j)/g_j$. By \cref{lem:chisq}(a), $0<Q<\infty$ a.s.\ when $\cI<\infty$; on the $P^{\mathsf x}$-null set where $Q\notin(0,\infty)$ we set $Q=1$ and $S=0$, which changes nothing almost surely and makes $Q$ take values in $(0,\infty)$, as \cref{lem:vantrees} requires. Given $\mathsf x$, the $y_j$ are independent $N(\gamma_jx_j,g_j)$. The density of the conditional law under $\beta$ with respect to that under $\beta=0$ (with the same $\alpha$) is, by (C1), the a.s.\ limit of $\prod_{j\le N}\exp\{(\beta x_j(y_j-\alpha_jx_j)-\beta^2x_j^2/2)/g_j\}$, which is $\exp(\beta S-\beta^2Q/2)$, the series for $S$ converging a.s.\ by (C2). Under $\beta$, $S-\beta Q=\sum_jx_jv_j/g_j$ is, given $\mathsf x$, $N(0,Q)$. Since $\int Q\,dP^{\mathsf x}=\sum_jf_j/g_j=\cI$, \cref{lem:vantrees} gives $\int\E_{\beta,\alpha}(T-\beta)^2\pi(d\beta)\ge(\cI+J(\pi))^{-1}$ for every prior with an absolutely continuous density and finite Fisher information.

\emph{Step 6 (part (ii)).} Let $w_j=1/h_j$ and $0<c\le g_j/h_j\le C$. Then $\hat\beta_N-\beta=B_N/D_N+M_N/D_N$ with
\[
D_N=\sum_{j\le N}w_jx_j^2,\qquad B_N=\sum_{j\le N}w_j\alpha_jx_j^2,\qquad M_N=\sum_{j\le N}w_jx_jv_j .
\]
(6a) $D_N=\sum_{j\le N}w_jf_j\xi_j^2$ and $\sum_jw_jf_j=\sum_j(g_j/h_j)q_j\ge c\,\cI=\infty$, so $D_N\to\infty$ a.s.\ by \cref{lem:chisq}(a). (6b) On this event, \cref{lem:toeplitz}(a) with the weights $w_jx_j^2$ gives $B_N/D_N\to\alpha_\infty$. (6c) Fix a realization of $X$ with $D_N\to\infty$ and $x_1\ne0$ (a probability-one event). Conditionally on $X$, the terms $\varepsilon_j=w_jx_jv_j$ are independent with mean zero and variance $w_j^2x_j^2g_j\le C\,w_jx_j^2$. \Cref{lem:kolmogorov} with $c_j=w_jx_j^2$ and $b_N=D_N$ gives $M_N/D_N\to0$ almost surely under the conditional law. Since $X$ and $V$ are independent, Fubini's theorem shows that $M_N/D_N\to0$ almost surely. Hence $\hat\beta_N\to\beta+\alpha_\infty$ a.s. Only the bounds $c\le g_j/h_j\le C$ were used, so one estimator works for every $g$ in the tail-equivalence class of $h$; for $\alpha\in\cA_0$ the limit is $\beta$.

\emph{Step 7 (part (iii)).} By Steps 1--2 and (C1), two parameters $(\beta,\alpha)$ and $(\beta',\alpha')$ with the same $(f,g)$ induce equivalent laws if $\delta=(\beta'-\beta)+(\alpha'-\alpha)\in\cN_q$ and singular laws otherwise. If $c+\alpha-\alpha'\in\cN_q$ for some $c\ne0$ and $\alpha,\alpha'\in\cA$, the parameters $(\beta,\alpha')$ and $(\beta+c,\alpha)$ have $\delta=c+\alpha-\alpha'\in\cN_q$, so their laws are equivalent while their values of $\beta$ differ; \cref{lem:obstruction} applies. If no such triple exists, then for any two parameters with $\beta\ne\beta'$ the sequence $\delta=(\beta'-\beta)+\alpha'-\alpha$ is not in $\cN_q$ (take $c=\beta'-\beta$ and exchange the roles of $\alpha$ and $\alpha'$), so the laws are singular. For $\cA=\cA_0$: if $\cI<\infty$, any $\alpha=\alpha'\in\cA_0$ and $c\ne0$ give $\sum_jc^2q_j=c^2\cI<\infty$. Conversely, if $c+\alpha-\alpha'\in\cN_q$ with $\alpha,\alpha'\in\cA_0$ and $c\ne0$, then $\delta_j=c+\alpha_j-\alpha'_j\to c$, so $\delta_j^2\ge c^2/4$ for $j\ge j_0$ and $\sum_{j\ge j_0}q_j\le(4/c^2)\sum_j\delta_j^2q_j<\infty$, that is, $\cI<\infty$. When $\cI=\infty$, part (ii) supplies a consistent estimator. \qed

\subsection{Proof of \texorpdfstring{\cref{cor:matern}}{Corollary 1}}\label{app:matern}

\emph{Idea.} The spectral dimension converts the decay of $q_j$ in $\lambda_j$ into a power of $j$, and the series $\sum_jq_j$ diverges exactly when that power is at least $-1$. At the boundary the partial sums grow like $\log N$, and the variance of GLS is the reciprocal of this growth.

\emph{Step 1 (divergence).} Since $\lambda_j\asymp j^{2/d_s}$ and $\lambda_j\to\infty$, also $1+\lambda_j\asymp j^{2/d_s}$, with constants valid for all $j\ge1$ after adjustment for finitely many $j$. Hence
\[
q_j=\frac{f_j}{g_j}\asymp(1+\lambda_j)^{-(\nu_X+d_s/2)+(\nu_V+d_s/2)}=(1+\lambda_j)^{\nu_V-\nu_X}\asymp j^{-2(\nu_X-\nu_V)/d_s}.
\]
By \cref{lem:powersums}, $\cI=\sum_jq_j=\infty$ if and only if $2(\nu_X-\nu_V)/d_s\le1$, which is \eqref{eq:threshold}. The constants in $\asymp$ absorb all other parameters (ranges, variances), which therefore do not affect the verdict.

\emph{Step 2 (estimability).} If \eqref{eq:threshold} holds, \cref{thm:frontier}(ii) gives a consistent estimator over $\R\times\cA_0$; if it fails, $\cI<\infty$ and \cref{thm:frontier}(i) excludes consistency.

\emph{Step 3 (variance at the boundary, oracle weights).} At $\nu_X=\nu_V+d_s/2$, $q_j\asymp1/j$, so $\cI_N=\sum_{j\le N}q_j\asymp\log N$ and $\sum_jq_j^2<\infty$; consistency follows from \cref{thm:frontier}(ii). With $w=1/g$, the noise term $M_N/D_N$ has, given $X$, the law $N(0,1/D_N)$ with $D_N=\sum_{j\le N}q_j\xi_j^2$. By \cref{lem:chisq}(b), $D_N/\cI_N\to1$ in probability, so the conditional variance is of order $1/\log N$ in probability. For the unconditional variance $\E D_N^{-1}$, Jensen's inequality gives $\E D_N^{-1}\ge1/\E D_N=1/\cI_N$. Conversely, splitting on $\{D_N\ge\cI_N/2\}$ and using the Cauchy--Schwarz and Chebyshev inequalities,
\[
\E D_N^{-1}\le\frac2{\cI_N}+\big(\E D_N^{-2}\big)^{1/2}\Prob\big(\abs{D_N-\cI_N}>\cI_N/2\big)^{1/2}\le\frac2{\cI_N}+\big(\E D_5^{-2}\big)^{1/2}\frac{\big(8\sum_{j\le N}q_j^2\big)^{1/2}}{\cI_N}
\]
for $N\ge5$, where $\E D_N^{-2}\le\E D_5^{-2}\le(\min_{j\le5}q_j)^{-2}\E(\chi^2_5)^{-2}=(\min_{j\le5}q_j)^{-2}/3$. Hence $\E D_N^{-1}\asymp1/\log N$.

\emph{Step 4 (tail-equivalent working spectra).} For $w=1/h$ with $c\le g_j/h_j\le C$, the conditional variance is $V_N/D_N^2$ with $V_N=\sum_{j\le N}w_j^2x_j^2g_j$ and $D_N=\sum_{j\le N}w_jx_j^2$. Since $c\le w_jg_j\le C$, $cD_N\le V_N\le CD_N$; and since $w_jf_j=(g_j/h_j)q_j$, $c\sum_{j\le N}q_j\xi_j^2\le D_N\le C\sum_{j\le N}q_j\xi_j^2$. Hence $V_N/D_N^2\asymp_P1/\log N$, and $\E[V_N/D_N^2]\le C\E D_N^{-1}\le C'\E\big(\sum_{j\le N}q_j\xi_j^2\big)^{-1}=O(1/\log N)$ by Step 3. Conversely, by Jensen's inequality, $\E[V_N/D_N^2]\ge c\,\E D_N^{-1}\ge c/\E D_N\ge c/(C\cI_N)$. Hence $\E[V_N/D_N^2]\asymp1/\log N$. \qed

\subsection{Proof of \texorpdfstring{\cref{prop:whatgls}}{Proposition 2}}\label{app:whatgls}

\emph{Idea.} GLS converges to the limit of the scale-specific associations whenever the Cameron--Martin energy diverges; in the phase-diagram model all spectral quantities are explicit, and the limit is decided by the ratio of two power laws.

\emph{Step 1 (part (a)).} This is \cref{thm:frontier}(ii), whose proof allows any limit $\alpha_\infty\in\R$; the limit $\beta+\alpha_\infty=\lim_j\gamma_j$ is the association at the finest scale.

\emph{Step 2 (part (b), spectral quantities).} Here $Z$ and $H$ are independent, Gaussian and $L$-stationary, so $f_j=f^Z_j+f^H_j$, $c_j=\Cov(z_j+\inner H{\varphi_j},cz_j)=cf^Z_j$, $s^U_j=c^2f^Z_j$,
\[
\alpha_j=\frac{cf^Z_j}{f^Z_j+f^H_j},\qquad g_j=c^2f^Z_j-\frac{c^2(f^Z_j)^2}{f^Z_j+f^H_j}=\frac{c^2f^Z_jf^H_j}{f^Z_j+f^H_j}>0 .
\]
By the inequality between arithmetic and geometric means, $q_j=f_j/g_j=(f^Z_j+f^H_j)^2/(c^2f^Z_jf^H_j)\ge4/c^2$, so $\cI=\infty$ in every regime, and Step 1 applies with any working spectrum tail-equivalent to $g$.

\emph{Step 3 (part (b), the limit).} Since $f^Z_j/f^H_j\sim(m_Z/m_H)\lambda_j^{a_H-a_Z}$ and $\lambda_j\to\infty$, the ratio $f^Z_j/f^H_j$ tends to $\infty$, $m_Z/m_H$ or $0$ according as $a_Z<a_H$, $a_Z=a_H$ or $a_Z>a_H$, and $\alpha_j=c\,(f^Z_j/f^H_j)/(1+f^Z_j/f^H_j)$ tends to $c$, $c\,m_Z/(m_Z+m_H)$ or $0$, which is \eqref{eq:threeregimes}. For Whittle--Mat\'ern spectra $\tau^2(\kappa^2+\lambda)^{-a}=\tau^2\lambda^{-a}(1+\kappa^2/\lambda)^{-a}\sim\tau^2\lambda^{-a}$, so $m=\tau^2$. The classical Mat\'ern spectral density on $\R^d$ with marginal variance $\sigma^2$ is $\sigma^2C_{\nu,d}\,\kappa^{2\nu}(\kappa^2+\abs\omega^2)^{-(\nu+d/2)}$ with $C_{\nu,d}=\Gamma(\nu+d/2)/\{\Gamma(\nu)\pi^{d/2}\}$, so $m\propto\sigma^2\kappa^{2\nu}$.

\emph{Step 4 (part (c)).} With the confounder independent of $X$, the response coefficients are $y_j=\beta s_jx_j+v_j=\{\beta+\beta(s_j-1)\}x_j+v_j$, so $\alpha_j=\beta(s_j-1)\to-\beta$ and, by part (a), GLS converges to $\beta-\beta=0$. If the confounder also carries confounding $\alpha^U$ as in \cref{sec:confounding}, the same argument gives the limit $\lim_j\alpha^U_j$. \qed

\subsection{Proof of \texorpdfstring{\cref{prop:pointsamples}}{Proposition 3}}\label{app:pointsamples}

\emph{Idea.} (a) The point data are a measurable function of the complete fields plus independent noise, so the equivalence of \cref{thm:frontier}(i) passes to them. (b) Quadrature turns point data into estimates of the first $J$ coefficients. Their errors are controlled by the mesh, the H\"older regularity of the paths and the averaged nugget, and the truncation $J_n$ grows slowly enough that the accumulated error vanishes. The sieve estimator is then a small perturbation of the complete-record GLS estimator, which is consistent.

\emph{Step 1 (part (a), point values as functions of the fields).} For $s\in M$ and $x\in L^2(M)$ let $\Psi_s(x)=\limsup_{r\downarrow0,\,r\in\mathbb Q}\mu(B(s,r))^{-1}\int_{B(s,r)}x\,d\mu$. Since $\operatorname{supp}\mu=M$, $\mu(B(s,r))>0$, and $\Psi_s$ is a Borel function on $L^2(M)$, being a countable $\limsup$ of continuous linear functionals. If $x$ is continuous, $\Psi_s(x)=x(s)$. The fields $X$ and $U$ have continuous paths, hence so does $Y=\beta X+U$, and $(X(s_i),Y(s_i)+\varepsilon_i)_{i\ge1}=(\Psi_{s_i}(X),\Psi_{s_i}(Y)+\varepsilon_i)_{i\ge1}$ almost surely, under every parameter value.

\emph{Step 2 (part (a), conclusion).} The map from coefficient sequences to $(X,Y)\in L^2(M)^2$ is the same measurable map under both parameters. By \cref{thm:frontier}(i), the laws of the coefficient sequences under $(\beta,\alpha)$ and $(\beta',\alpha)$ are equivalent, hence so are the laws of $(X,Y)$ (\cref{lem:kernels}(b)). The data sequence is obtained from $(X,Y)$ by the kernel ``apply $\Psi_{s_i}$ and add independent noise with the common law of $(\varepsilon_i)$'', so \cref{lem:kernels}(b) and \cref{lem:obstruction} conclude. No property of the design was used.

\emph{Step 3 (part (b), quadrature errors).} Let $K_X,K_Y$ be the (a.s.\ finite) H\"older constants of the paths of order $\vartheta$, so that $\abs{X(s)-X(s')}\le K_Xd(s,s')^\vartheta$ and similarly for $Y$, and note that $\operatorname{diam}A_i\le\ell_n$. Since $x_j=\sum_i\int_{A_i}X\varphi_j\,d\mu$ and $\hat x_j=\sum_i\int_{A_i}X(s_i)\varphi_j(s_i)\,d\mu$,
\begin{align*}
\abs{\hat x_j-x_j}&\le\sum_i\int_{A_i}\big(\abs{X(s_i)}\abs{\varphi_j(s_i)-\varphi_j(s)}+\abs{\varphi_j(s)}\abs{X(s_i)-X(s)}\big)d\mu(s)\\
&\le\mu(M)\big(\norm X_\infty\Lambda_j\ell_n+B_jK_X\ell_n^\vartheta\big).
\end{align*}
The same bound holds for $Y$, with $K_Y$ and $\norm Y_\infty$, plus the noise term $\zeta_j=\sum_iw_i\varepsilon_i\varphi_j(s_i)$, for which $\E\zeta_j^2=\sigma^2\sum_iw_i^2\varphi_j(s_i)^2\le\sigma^2B_j^2\mu(M)\max_iw_i$, $\sigma^2$ being the nugget variance. Let $e_j=\abs{\hat x_j-x_j}+\abs{\hat y_j-y_j}$ and $E_J=\sum_{j\le J}e_j/h_j$. For $\ell_n\le1$, with $R_1=\mu(M)(\norm X_\infty+\norm Y_\infty+K_X+K_Y)$,
\[
E_J\le R_1a_J\ell_n^\vartheta+\sum_{j\le J}\frac{\abs{\zeta_j}}{h_j},\qquad\E\sum_{j\le J}\frac{\abs{\zeta_j}}{h_j}\le\sigma\mu(M)^{1/2}\big(\max_iw_i\big)^{1/2}a_J .
\]
Since $a_{J_n}\le\log(1/\varpi_n)$ and $\ell_n,\max_iw_i\le\varpi_n\to0$, we get $R_1a_{J_n}\ell_n^\vartheta\le R_1\varpi_n^\vartheta\log(1/\varpi_n)\to0$ a.s.\ for every $\vartheta>0$, and the expectation of the second term is at most $\sigma\mu(M)^{1/2}\varpi_n^{1/2}\log(1/\varpi_n)\to0$. Hence $E_{J_n}\to0$ in probability. Moreover $J_n$ is finite because $a_J\to\infty$, since $B_j\ge\mu(M)^{-1/2}$ and $h_j\to0$; and $J_n\to\infty$, because each $a_J$ is finite and $\log(1/\varpi_n)\to\infty$.

\emph{Step 4 (part (b), perturbation of the ratio).} Let $N_J=\sum_{j\le J}x_jy_j/h_j$ and $D_J=\sum_{j\le J}x_j^2/h_j$ be the numerator and denominator of the complete-record estimator $\hat\beta_J=N_J/D_J$, and $\hat N_J,\hat D_J$ their plug-in versions, so $\tilde\beta_n=\hat N_{J_n}/\hat D_{J_n}$. Let $Q_0=\max(\norm X_{L^2},\norm Y_{L^2})$, which bounds every $\abs{x_j}$ and $\abs{y_j}$, and $\bar h=\max(1,\sup_jh_j)$. Since $\abs{\hat x_j^2-x_j^2}\le2Q_0e_j+e_j^2$, $\abs{\hat x_j\hat y_j-x_jy_j}\le2Q_0e_j+e_j^2$ and $\sum_je_j^2/h_j=\sum_jh_j(e_j/h_j)^2\le\bar h\big(\sum_je_j/h_j\big)^2$,
\[
\abs{\hat D_J-D_J}+\abs{\hat N_J-N_J}\le2(2Q_0+\bar hE_J)E_J .
\]
Now $D_J\ge D_1=x_1^2/h_1>0$ a.s., and $\sup_J\abs{\hat\beta_J}<\infty$ a.s.\ because $\hat\beta_J$ converges by \cref{thm:frontier}(ii). On the event $\{2(2Q_0+\bar hE_{J_n})E_{J_n}<D_1/2\}$, whose probability tends to one by Step 3, $\hat D_{J_n}\ge D_1/2$ and
\[
\abs{\tilde\beta_n-\hat\beta_{J_n}}=\frac{\abs{(\hat N_{J_n}-N_{J_n})-\hat\beta_{J_n}(\hat D_{J_n}-D_{J_n})}}{\hat D_{J_n}}\le\frac{4(2Q_0+\bar hE_{J_n})E_{J_n}\big(1+\sup_J\abs{\hat\beta_J}\big)}{D_1}\longrightarrow0
\]
in probability. Finally $J_n\to\infty$ and $\hat\beta_J\to\beta$ a.s.\ by \cref{thm:frontier}(ii), so $\tilde\beta_n\to\beta$ in probability. \qed
\subsection{Proof of \texorpdfstring{\cref{thm:resolution}}{Theorem 2} and \texorpdfstring{\cref{cor:latticetarget}}{Corollary 2}}\label{app:resolution}

\emph{Idea.} Lattice data are a stationary Gaussian process on $\Z^d$, so their law is fixed by the aliased spectral matrix on $T_\delta$. This matrix sees the continuum model only through three functions of $\omega\in T_\delta$: the aliased covariate spectrum, the alias mean $\mu=\beta+\bar\alpha^\delta$ of the association, and the alias variance $v^\delta$ of the confounding, plus the aliased residual spectrum. To move $\beta$ to $\beta'$ we redefine the confounding below a cutoff frequency $R$ so that the alias mean is unchanged; this costs an alias variance that becomes small as $R$ grows. Under (a) the residual spectrum absorbs the change in alias variance; under (b) a zero-mean perturbation restores the original alias variance. Identification requires that the alias variance be (nearly) zero somewhere and that $g$ and the nugget variance be known exactly, and then the Cauchy--Schwarz inequality pins down $\beta$.

\emph{Step 1 (spectral representation).} By (C8), $X(s)=\int e^{is\cdot\omega}Z_X(d\omega)$ with an orthogonal random measure satisfying $\E\abs{Z_X(d\omega)}^2=f_X(\omega)d\omega$ and $Z_X(-d\omega)=\overline{Z_X(d\omega)}$, and the filtered field is $(\alpha\ast X)(s)=\int e^{is\cdot\omega}\alpha(\omega)Z_X(d\omega)$, which is real because $\alpha$ is Hermitian. Since $V$ is independent of $X$, the isometry of (C8) gives, for a lag $r\in\R^d$,
\[
\Cov\big(Y(s+r),X(s)\big)=\int e^{ir\cdot\omega}\gamma(\omega)f_X(\omega)\,d\omega,\qquad\Cov\big(Y(s+r),Y(s)\big)=\int e^{ir\cdot\omega}\big(\abs{\gamma}^2f_X+g\big)(\omega)\,d\omega .
\]

\emph{Step 2 (periodization).} The map $(\omega,m)\mapsto\omega_m=\omega+2\pi m/\delta$ is a measure-preserving bijection from $T_\delta\times\Z^d$ (Lebesgue times counting measure) onto $\R^d$. Hence, for integrable $F:\R^d\to\mathbb C$, Tonelli's theorem gives $\int_{T_\delta}\sum_m\abs{F(\omega_m)}d\omega=\int_{\R^d}\abs F<\infty$, so $\sum_mF(\omega_m)$ converges absolutely for a.e.\ $\omega$ and defines an element of $L^1(T_\delta)$, and for $k\in\Z^d$, because $e^{i\delta k\cdot\omega_m}=e^{i\delta k\cdot\omega}e^{2\pi ik\cdot m}=e^{i\delta k\cdot\omega}$,
\[
\int_{\R^d}e^{i\delta k\cdot\omega}F(\omega)\,d\omega=\int_{T_\delta}e^{i\delta k\cdot\omega}\sum_mF(\omega_m)\,d\omega .
\]
We apply this to $F=f_X$, $\gamma f_X$ and $\abs\gamma^2f_X+g$, all integrable because $\gamma$ is bounded and $f_X,g$ are integrable.

\emph{Step 3 (lattice laws).} The point samples form a centred stationary Gaussian process on $\Z^d$ with values in $\R^2$, whose finite-dimensional laws, and hence whose law on $(\R^2)^{\Z^d}$, are determined by its covariances. By Steps 1--2 these covariances are the Fourier coefficients of $(F^\delta_X,F^\delta_{XY},F^\delta_Y)$ with $F^\delta_X=\sum_mf_X(\omega_m)$, $F^\delta_{XY}=\sum_m\gamma(\omega_m)f_X(\omega_m)$ and $F^\delta_Y=\sum_m(\abs\gamma^2f_X+g)(\omega_m)$. Independent $N(0,\sigma^2)$ nugget errors add $\sigma^2$ to the lag-zero covariance of the response, that is, the constant $\sigma^2(\delta/2\pi)^d$ to $F^\delta_Y$, because $\int_{T_\delta}e^{i\delta k\cdot\omega}d\omega=(2\pi/\delta)^d\ind\{k=0\}$. An element of $L^1(T_\delta)$ is determined a.e.\ by its Fourier coefficients. Therefore \emph{two parameter values induce the same law of $D_\delta$ if and only if their aliased spectral matrices agree almost everywhere on $T_\delta$.}

\emph{Step 4 (alias mean and alias variance).} For point samples the weights $w_m=f_X(\omega_m)/F^\delta_X$ are positive and sum to one. With $\mu=\sum_mw_m\gamma(\omega_m)=\beta+\bar\alpha^\delta$, we have $F^\delta_{XY}=F^\delta_X\mu$, and expanding $\abs{\gamma(\omega_m)}^2=\abs{\mu+(\gamma(\omega_m)-\mu)}^2$ with $\sum_mw_m(\gamma(\omega_m)-\mu)=0$ gives $\sum_mw_m\abs{\gamma(\omega_m)}^2=\abs\mu^2+v^\delta$, where $v^\delta=\sum_mw_m\abs{\gamma(\omega_m)-\mu}^2=\sum_mw_m\abs{\alpha(\omega_m)-\bar\alpha^\delta}^2$. Hence, with $G^\delta=\sum_mg(\omega_m)$,
\[
F^\delta_{XY}=F^\delta_X\,\mu,\qquad F^\delta_Y=F^\delta_X\big(\abs\mu^2+v^\delta\big)+G^\delta\ \big(+\sigma^2(\delta/2\pi)^d\big).
\]
So the lattice law is determined by $(F^\delta_X,\mu,v^\delta,G^\delta,\sigma^2)$, which proves the first statement of \cref{cor:latticetarget}.

\emph{Step 5 (conventions and constants).} Functions of $\omega\in T_\delta$ are extended $(2\pi/\delta)$-periodically, and such a function $b$ is regarded as a function on $\R^d$ through $b(\omega_m)=b(\omega)$. Expressions such as $\sum_mw_m(\omega)\chi(\omega_m)\alpha(\omega_m)$ are automatically periodic in $\omega\in\R^d$, since a shift of $\omega$ by $2\pi e/\delta$ only relabels $m$. Because $f_X$ is even, $w_m(-\omega)=w_{-m}(\omega)$, and $(-\omega)_m=-\omega_{-m}$; hence periodic functions built from $w_m$, from even real functions and from the Hermitian $\alpha$ satisfy $b(-\omega)=\overline{b(\omega)}$, and $v^\delta$ is even. By (R1), $F^\delta_X$ is continuous and positive, so $w_0=f_X/F^\delta_X$ is continuous and positive on the compact set $\bar T_\delta$ and $c_0=\inf_{\bar T_\delta}w_0>0$. Put $A_\infty=\norm\alpha_\infty$, $\Delta=\beta'-\beta$, $C_\Delta=\abs\Delta+2A_\infty$ and $\epsilon_R=\sup_{\omega\in\bar T_\delta}\sum_{m:\abs{\omega_m}>R}w_m(\omega)$. For $\omega\in\bar T_\delta$, $\abs{\omega_m}>R$ forces $\abs m\ge(R-\sqrt d\,\pi/\delta)\delta/(2\pi)$, so by the uniform convergence of the alias sums in (R1) and $F^\delta_X\ge\inf F^\delta_X>0$, $\epsilon_R\to0$ as $R\to\infty$.

\emph{Step 6 (part (a): the construction).} Choose an even $\tilde\chi\in C_c(\R^d)$ with values in $[0,1]$ and $\tilde\chi=1$ on $\bar T_\delta$, and let $g_\ast>0$ be the essential infimum of $g$ on $\operatorname{supp}\tilde\chi$, positive by (R1). For $R>\sqrt d\,\pi/\delta$, so that $\bar T_\delta\subset B_R$, choose an even $\chi\in C_c(\R^d)$ with values in $[0,1]$ and $\chi=1$ on $B_R$. Writing $\alpha_m=\alpha(\omega_m)$ and $\chi_m=\chi(\omega_m)$, let $W_\chi=\sum_mw_m\chi_m$, which lies in $[c_0,1]$ because $\chi_0=1$, and define
\[
\bar a=\frac{\sum_mw_m\chi_m\alpha_m-\Delta}{W_\chi},\qquad\gamma'(x)=\beta'+\big(1-\chi(x)\big)\alpha(x)+\chi(x)\,\bar a(x),\qquad\alpha'=\gamma'-\beta'.
\]
The sum defining $\bar a$ has finitely many non-zero terms, since $\chi$ has compact support; so $\bar a$ is bounded, periodic and Hermitian, and continuous if $\alpha$ is. Then $\alpha'-\alpha=\chi(\bar a-\alpha)$ is bounded, Hermitian and compactly supported, and $\alpha'=\alpha$ outside $\operatorname{supp}\chi$, so $\alpha'$ satisfies (R2), and $\alpha'$ is continuous when $\alpha$ is.

\emph{Step 7 (part (a): the alias mean is preserved).} Using $\sum_mw_m=1$,
\[
\sum_mw_m\gamma'(\omega_m)=\beta'+\sum_mw_m(1-\chi_m)\alpha_m+\bar aW_\chi=\beta'+\sum_mw_m\alpha_m-\Delta=\beta+\bar\alpha^\delta=\mu .
\]

\emph{Step 8 (part (a): the new alias variance is small).} We have $\gamma'(\omega_m)-\mu=(1-\chi_m)(\Delta+\alpha_m-\bar\alpha^\delta)+\chi_m(\beta'+\bar a-\mu)$, and a direct computation, using $1-W_\chi=\sum_kw_k(1-\chi_k)$, gives
\[
\beta'+\bar a-\mu=\frac{\Delta W_\chi+\sum_kw_k\chi_k\alpha_k-\Delta-W_\chi\bar\alpha^\delta}{W_\chi}=\frac1{W_\chi}\sum_kw_k(1-\chi_k)\big(\bar\alpha^\delta-\Delta-\alpha_k\big).
\]
Since $\chi=1$ on $B_R$, $\sum_kw_k(1-\chi_k)\le\epsilon_R$. The first bracket is at most $C_\Delta$ in modulus and the second at most $\epsilon_RC_\Delta/c_0$, so, by $\abs{a+b}^2\le2\abs a^2+2\abs b^2$,
\[
v'=\sum_mw_m\abs{\gamma'(\omega_m)-\mu}^2\le2C_\Delta^2\Big(\sum_mw_m(1-\chi_m)^2+\frac{\epsilon_R^2}{c_0^2}\Big)\le2C_\Delta^2\big(\epsilon_R+\epsilon_R^2/c_0^2\big)=:\eta_R\longrightarrow0 .
\]

\emph{Step 9 (part (a): the residual spectrum absorbs the change).} Let $E=F^\delta_X(v'-v^\delta)$, $K=\sum_m\tilde\chi(\omega_m)\ge1$ and $g'=g-\tilde\chi E/K$, with $E$ and $K$ periodic. Then $\sum_mg'(\omega_m)=G^\delta-E\sum_m\tilde\chi(\omega_m)/K=G^\delta-E$, hence
\[
F^\delta_X\big(\abs\mu^2+v'\big)+\sum_mg'(\omega_m)=F^\delta_X\big(\abs\mu^2+v^\delta\big)+G^\delta=F^\delta_Y,
\]
and the aliased cross spectrum is $F^\delta_X\mu$ under both parameters (Step 7). By Step 3 the two parameter values induce the same law of $D_\delta$, with any common nugget. The function $g'$ is even (as $\tilde\chi$, $F^\delta_X$, $v'$, $v^\delta$ and $K$ are) and differs from $g$ only on $\operatorname{supp}\tilde\chi$, by at most $\sup F^\delta_X(\eta_R+4A_\infty^2)$, and it is lower than $g$ by at most $\sup F^\delta_X\,\eta_R$. Choosing $R$ with $\sup F^\delta_X\,\eta_R<g_\ast/2$ gives $g'\ge g_\ast/2$ on $\operatorname{supp}\tilde\chi$ and $g'=g$ elsewhere, so $g'$ satisfies (R1): only finitely many aliases of a compact set meet $\operatorname{supp}\tilde\chi$, so the alias sums of $g'$ converge essentially uniformly on compact sets. If $\alpha$ and $g$ are continuous, then $v'$ and $v^\delta$ are uniformly convergent series of continuous functions (each term is bounded by $(C_\Delta+\abs{\bar a})^2w_m$ and $\sum_mw_m$ converges uniformly), so $E$, and hence $g'$, is continuous.

\emph{Step 10 (part (b): a zero-mean perturbation).} Let $v_t$ be a target alias variance with $0<\operatorname{ess\,inf}v_t\le\sup v_t<\infty$; we construct $\alpha'$ with alias mean $\mu$ and alias variance $v_t$. Take $R$ in Step 6 with $\eta_R<\operatorname{ess\,inf}v_t$, and call the resulting function $\gamma^0$; its alias variance $v^0\le\eta_R$ is smaller than $v_t$ almost everywhere. Let $\chi_1\in C_c(\R^d)$ be even, with values in $[0,1]$, $\chi_1=1$ on $\bar T_\delta$ and support in the cube $\{x:\abs{x_i}<5\pi/(4\delta)\ \forall i\}$, and let $\chi_2(x)=\chi_1(x-2\pi e_1/\delta)+\chi_1(x+2\pi e_1/\delta)$, which is even. For point samples all $w_m$ are positive, and
\[
\Psi=\frac{\sum_mw_m\chi_1(\omega_m)}{\sum_mw_m\chi_2(\omega_m)}
\]
is continuous, periodic, even and positive: the sums are finite, the numerator is at least $w_0>0$, and the denominator is at least $w_{e_1}\chi_2(\omega_{e_1})+w_{-e_1}\chi_2(\omega_{-e_1})=w_{e_1}+w_{-e_1}>0$, because $\chi_2(\omega_{\pm e_1})=\chi_1(\omega)+\chi_1(\omega_{\pm2e_1})=1$ for $\omega\in\bar T_\delta$ (the first coordinate of $\omega_{\pm2e_1}$ has modulus at least $3\pi/\delta$). The function $\zeta=\chi_1-\Psi\chi_2$ is real, even and compactly supported, and $\sum_mw_m\zeta(\omega_m)=0$ by the choice of $\Psi$.

\emph{Step 11 (part (b): $\zeta$ has positive alias variance).} Let $\Upsilon=\sum_mw_m\zeta(\omega_m)^2$, a finite sum of continuous functions. Suppose $\Upsilon(\omega)=0$ for some $\omega\in\bar T_\delta$. Then $\zeta(\omega_m)=0$ for every $m$, and in particular $\zeta(\omega_{\pm e_1})=\chi_1(\omega_{\pm e_1})-\Psi(\omega)=0$, so $\chi_1(\omega_{e_1})=\chi_1(\omega_{-e_1})=\Psi(\omega)>0$. But $\chi_1(\omega_{e_1})\ne0$ forces $\abs{\omega_1+2\pi/\delta}<5\pi/(4\delta)$, that is $\omega_1<-3\pi/(4\delta)$, and $\chi_1(\omega_{-e_1})\ne0$ forces $\omega_1>3\pi/(4\delta)$, a contradiction. Hence $\Upsilon>0$ on the compact set $\bar T_\delta$, and $\Upsilon\ge c_1>0$.

\emph{Step 12 (part (b): solving for the target).} For $\Lambda_0>0$ and $t\in\R$, the function $\gamma^0+t\Lambda_0\zeta$ has alias mean $\mu$ and alias variance
\[
Q_\omega(t)=v^0+2t\Lambda_0\Gamma+t^2\Lambda_0^2\Upsilon,\qquad \Gamma=\operatorname{Re}\sum_mw_m\overline{\big(\gamma^0(\omega_m)-\mu\big)}\,\zeta(\omega_m),
\]
with $\abs\Gamma\le\sqrt{v^0\Upsilon}$ by the Cauchy--Schwarz inequality. Choose $\Lambda_0\ge1+(\sqrt{\sup v_t}+\sqrt{\eta_R})/\sqrt{c_1}$; then $Q_\omega(1)\ge(\Lambda_0\sqrt\Upsilon-\sqrt{v^0})^2\ge(\Lambda_0\sqrt{c_1}-\sqrt{\eta_R})^2>\sup v_t$. Wherever $Q_\omega(0)=v^0<v_t(\omega)$, which holds almost everywhere, the quadratic $t\mapsto Q_\omega(t)-v_t(\omega)$ has positive leading coefficient and a negative value at $0$, hence exactly one positive root, and this root lies in $(0,1)$ because $Q_\omega(1)>v_t(\omega)$:
\[
t^\ast(\omega)=\frac{-\Lambda_0\Gamma+\sqrt{\Lambda_0^2\Gamma^2+\Lambda_0^2\Upsilon\,(v_t-v^0)}}{\Lambda_0^2\Upsilon};
\]
on the remaining null set put $t^\ast=0$. Since $\Gamma,\Upsilon,v^0$ and $v_t$ are measurable, periodic and even, so is $t^\ast$, and it is continuous when they are. Set $\gamma'=\gamma^0+t^\ast\Lambda_0\zeta$ and $\alpha'=\gamma'-\beta'$. With $\alpha^0=\gamma^0-\beta'$, the function $\alpha'-\alpha=(\alpha^0-\alpha)+t^\ast\Lambda_0\zeta$ is bounded, Hermitian (the second term is real and even) and compactly supported, so $\alpha'$ satisfies (R2), and $\gamma'$ has alias mean $\mu$ and alias variance $v_t$ a.e.

\emph{Step 13 (part (b): non-identification).} If $\operatorname{ess\,inf}v^\delta>0$, take $v_t=v^\delta$, which is admissible since $v^\delta\le4A_\infty^2$. Then $(\beta',\alpha',f_X,g)$ has the aliased spectral matrix of $(\beta,\alpha,f_X,g)$, and Step 3 gives equal laws, with any common nugget. If the data carry a nugget of unknown variance $\sigma^2>0$, take any $\sigma'^2\in[0,\sigma^2)$ and $v_t=v^\delta+(\sigma^2-\sigma'^2)(\delta/2\pi)^d/F^\delta_X$. This target is bounded above, and bounded below by $(\sigma^2-\sigma'^2)(\delta/2\pi)^d/\sup F^\delta_X>0$ for every $\alpha$ satisfying (R2). Then $F^\delta_Xv_t+\sigma'^2(\delta/2\pi)^d=F^\delta_Xv^\delta+\sigma^2(\delta/2\pi)^d$, so $(\beta',\alpha',f_X,g,\sigma'^2)$ has the aliased spectral matrix of $(\beta,\alpha,f_X,g,\sigma^2)$.

\emph{Step 14 (part (b): identification).} Let the nugget variance be known (or zero), let $\operatorname{ess\,inf}v^\delta=0$, and suppose that $(\beta',\alpha',f_X,g)$, with $\alpha'$ satisfying (R2), induces the same law of $D_\delta$. By Step 3 the aliased spectral matrices agree a.e.; since $F^\delta_X>0$ and the terms involving $g$ and the nugget coincide, $\mu'=\mu$ and $v'=v^\delta$ a.e., where $\mu'$ and $v'$ are the alias mean and the alias variance of $\gamma'=\beta'+\alpha'$. For $R>0$ let $S_R(\omega)=\{m:\abs{\omega_m}>R\}$ and $t_R=\inf_{\omega\in\bar T_\delta}\sum_{m\in S_R(\omega)}w_m(\omega)$; then $t_R>0$, because for $M_0$ large every $\omega\in\bar T_\delta$ has $M_0e_1\in S_R(\omega)$, and $w_{M_0e_1}$ is continuous and positive on $\bar T_\delta$. Let $\eta_\alpha(R)=\operatorname{ess\,sup}_{\abs x>R}\abs{\alpha(x)}$. By the measure-preserving bijection of Step 2, the set of $\omega\in T_\delta$ for which some alias $\omega_m$ with $\abs{\omega_m}>R$ has $\abs{\alpha(\omega_m)}>\eta_\alpha(R)$ is null; taking $R\in\N$ we may discard one null set for all $R$ at once. For the remaining $\omega$, averaging $\abs{\bar\alpha^\delta(\omega)}\le\abs{\alpha(\omega_m)}+\abs{\alpha(\omega_m)-\bar\alpha^\delta(\omega)}$ over $m\in S_R(\omega)$ with weights $w_m/\sum_{k\in S_R}w_k$ and using the Cauchy--Schwarz inequality gives
\[
\abs{\bar\alpha^\delta(\omega)}\le\eta_\alpha(R)+\Big(\frac{\sum_{m\in S_R}w_m\abs{\alpha(\omega_m)-\bar\alpha^\delta}^2}{\sum_{m\in S_R}w_m}\Big)^{1/2}\le\eta_\alpha(R)+\big(v^\delta(\omega)/t_R\big)^{1/2}.
\]
The same bound holds for $\alpha'$, whose alias variance is also $v^\delta$. Since $\beta+\bar\alpha^\delta=\mu=\mu'=\beta'+\bar\alpha'^{\,\delta}$,
\[
\abs{\beta-\beta'}\le\eta_\alpha(R)+\eta_{\alpha'}(R)+2\big(v^\delta(\omega)/t_R\big)^{1/2}\qquad\text{for almost every }\omega .
\]
Because $\operatorname{ess\,inf}v^\delta=0$, for each $\epsilon>0$ the set $\{v^\delta<\epsilon\}$ has positive measure, so we may choose $\omega$ in it, and $\abs{\beta-\beta'}\le\eta_\alpha(R)+\eta_{\alpha'}(R)+2(\epsilon/t_R)^{1/2}$. Letting $\epsilon\to0$ and then $R\to\infty$, and using (R2) for $\alpha$ and $\alpha'$, gives $\beta'=\beta$.

\emph{Step 15 (block means: reduction).} The block mean over $Q_k$ is $\bar X_k=\delta^{-d}\int_{\delta k+[0,\delta)^d}X(u)\,du=(X\ast\psi_\delta)(\delta k)$ with $\psi_\delta=\delta^{-d}\ind_{[0,\delta)^d}(-\,\cdot)$, and similarly for $Y$. By (C8), filtering both fields with the integrable kernel $\psi_\delta$ multiplies the spectral measure $Z_X(d\omega)$ by $\hat\psi_\delta(\omega)$, so every spectral and cross-spectral density is multiplied by $\abs{\hat\psi_\delta}^2$, with $\abs{\hat\psi_\delta(\omega)}=\prod_i\abs{\operatorname{sinc}(\delta\omega_i/2)}$ and $\operatorname{sinc}(x)=\sin x/x$. Steps 1--4 therefore hold for $B_\delta$ with $f_X\abs{\hat\psi_\delta}^2$ in place of $f_X$ (in $w_m$ and in $F^\delta_X$) and $\sum_mg(\omega_m)\abs{\hat\psi_\delta(\omega_m)}^2$ in place of $G^\delta$. Since $\operatorname{sinc}(x)\ge2/\pi$ for $\abs x\le\pi/2$, $\abs{\hat\psi_\delta}\ge(2/\pi)^d$ on $\bar T_\delta$.

\emph{Step 16 (block means: part (a)).} Write $F^\delta_{X,\psi}=\sum_mf_X(\omega_m)\abs{\hat\psi_\delta(\omega_m)}^2\le F^\delta_X$ for the aliased covariate spectrum of the block means. Steps 5--9 carry over verbatim, with the weights of Step 15, with $F^\delta_{X,\psi}$ in place of $F^\delta_X$ and with $K_\psi=\sum_m\tilde\chi(\omega_m)\abs{\hat\psi_\delta(\omega_m)}^2\ge(2/\pi)^{2d}$ in place of $K$. Indeed $c_0\ge(2/\pi)^{2d}\inf_{\bar T_\delta}f_X/\sup F^\delta_{X,\psi}>0$, and $g'=g-\tilde\chi E/K_\psi$ gives $\sum_mg'(\omega_m)\abs{\hat\psi_\delta(\omega_m)}^2=\sum_mg(\omega_m)\abs{\hat\psi_\delta(\omega_m)}^2-E$. The reduction of $g$ is at most $(\pi/2)^{2d}\sup F^\delta_{X,\psi}\,\eta_R$, and $R$ is chosen to make it smaller than $g_\ast/2$.

\emph{Step 17 (block means: the nugget is identified).} For a non-principal alias $m\ne0$, some $m_i\ne0$, and $\abs{\operatorname{sinc}(\delta(\omega_i+2\pi m_i/\delta)/2)}=\abs{\sin(\delta\omega_i/2)}/\abs{\delta\omega_i/2+\pi m_i}\to0$ as $\omega\to0$. Given $\epsilon>0$, choose $M_1$ with $\sum_{\abs m>M_1}f_X(\omega_m)<\epsilon$ on $\bar T_\delta$; the finitely many terms with $0<\abs m\le M_1$ tend to zero as $\omega\to0$. Since the denominator of $w_0$ is at least $(2/\pi)^{2d}f_X(\omega)$, we get $1-w_0(\omega)=\sum_{m\ne0}w_m(\omega)\to0$ as $\omega\to0$. For every bounded $\alpha$, $v^\delta\le4\norm\alpha_\infty^2(1-w_0)$: the alias mean minimizes $c\mapsto\sum_mw_m\abs{\alpha_m-c}^2$, so taking $c=\alpha_0$ gives $v^\delta\le\sum_{m\ne0}w_m\abs{\alpha_m-\alpha_0}^2\le4\norm\alpha_\infty^2(1-w_0)$. The function $H=F^\delta_Y-\sum_mg(\omega_m)\abs{\hat\psi_\delta(\omega_m)}^2-F^\delta_X\abs\mu^2=F^\delta_Xv^\delta+\sigma^2(\delta/2\pi)^d$ is determined a.e.\ by the law of $B_\delta$ (through $F^\delta_X$, $F^\delta_{XY}=F^\delta_X\mu$ and $F^\delta_Y$) and by $g$, and $0\le H-\sigma^2(\delta/2\pi)^d\le4\sup F^\delta_X\norm\alpha_\infty^2(1-w_0)$. Hence $\sigma^2(\delta/2\pi)^d=\lim_{r\downarrow0}\operatorname{ess\,inf}_{\abs\omega<r}H(\omega)$ is determined by the law and $g$, whatever the bounded $\alpha$. Two parameter values with the same law and the same $g$ therefore have the same nugget variance.

\emph{Step 18 (block means: part (c)).} Let $g$ be known and suppose $(\beta',\alpha',f_X,g)$ induces the same law of $B_\delta$. By Step 17 the nugget variances agree, and as in Step 14, $\mu'=\mu$ and $v'=v^\delta$ a.e. Let $S\subseteq T_\delta$ have positive measure with $\alpha(\omega_m)=0$ for all $m$ and $\omega\in S$. On $S$, $\mu=\beta$ and $v^\delta=0$, hence $v'=0$ a.e.\ on $S$. A block-mean weight $w_m(\omega)$ vanishes only if $\delta(\omega_i+2\pi m_i/\delta)/2\in\pi\Z\setminus\{0\}$ for some $i$, which for $\omega\in T_\delta$ requires $\omega_i=0$; so for a.e.\ $\omega\in S$ all weights are positive, and $\sum_mw_m\abs{\gamma'(\omega_m)-\beta}^2=0$ forces $\gamma'(\omega_m)=\beta$, that is $\alpha'(\omega_m)=\beta-\beta'$, for all $m$. As in Step 14, (R2) for $\alpha'$ gives $\alpha'(\omega_m)\to0$ as $\abs m\to\infty$ for a.e.\ $\omega$. Choosing $\omega\in S$ outside the exceptional null sets, $\beta-\beta'=\lim_m\alpha'(\omega_m)=0$.

\emph{Step 19 (block means: part (b) is void).} By Step 17, $v^\delta\to0$ as $\omega\to0$ for every bounded $\alpha$, so the hypothesis $\operatorname{ess\,inf}v^\delta>0$ of part (b) never holds for block means.

\emph{Step 20 (the final statement of \cref{thm:resolution}).} By part (a), for every $\beta'\in\R$, of either sign, there is an alternative satisfying \cref{ass:resolution} with the same law of the whole lattice, with or without a nugget (the nugget is common to both models). Data from any finite window are the image of the lattice data under a coordinate projection, so their laws also agree (\cref{lem:kernels}(b)). Hence no functional of $\beta$, including its sign, is identified.

\emph{Step 21 (\cref{cor:latticetarget}).} The first statement is Step 4, with $\mu=\gamma^\delta$ by \eqref{eq:aliasavg}. When $\alpha$ is continuous, $\mu$ is continuous (the series converges uniformly), so its value at $\omega^\ast$ is defined. An analysis that sets the coherence to zero at $\omega^\ast$ equates $\beta$ with $\gamma^\delta(\omega^\ast)=\beta+\sum_mw_m(\omega^\ast)\alpha(\omega^\ast_m)$. At $\omega^\ast=(\pi/\delta,\dots,\pi/\delta)$ no coordinate vanishes, so all weights $w_m(\omega^\ast)$ are positive for point samples and for block means; if $\alpha\ge0$ and $\alpha(\omega^\ast_m)>0$ for some $m$, the target exceeds $\beta$. \qed

\subsection{Proof of \texorpdfstring{\cref{thm:window}}{Theorem 3} and \texorpdfstring{\cref{prop:estimators}}{Proposition 4}}\label{app:estimators}

\emph{Idea.} Every estimator considered is signal plus noise. The signal is a weighted average of the $\gamma_j$, and it forgets the confounded low frequencies exactly when the total weight diverges. The noise is conditionally Gaussian given $X$, so it vanishes if and only if its conditional variance does (\cref{lem:condgauss}). Under the power laws \eqref{eq:powerlaw} both the total weight and the conditional variance are power sums, and comparing exponents gives each verdict.

\emph{Step 1 (decomposition and orders).} Write $w_j=1/h_j$, $D_N=\sum_{j\le N}w_jx_j^2$, $B_N=\sum_{j\le N}w_j\alpha_jx_j^2$, $M_N=\sum_{j\le N}w_jx_jv_j$ and $V_N=\sum_{j\le N}w_j^2x_j^2g_j$. Then $\hat\beta_N-\beta=B_N/D_N+M_N/D_N$, and given $X$, $M_N\sim N(0,V_N)$. Under \eqref{eq:powerlaw}, $\E D_N=\sum_{j\le N}w_jf_j$ with $w_jf_j\asymp j^{\varkappa-a}=j^u$, and $\E V_N=\sum_{j\le N}w_j^2f_jg_j$ with $w_j^2f_jg_j\asymp j^{2\varkappa-a-b}=j^{2u+\Delta'}$. Also $\Var D_N=2\sum_{j\le N}(w_jf_j)^2$ and $\Var V_N=2\sum_{j\le N}(w_j^2f_jg_j)^2$. By \cref{lem:powersums} and \cref{lem:chisq}(b): if $u\ge-1$, then $D_N/\E D_N\to1$ in probability; if $2u+\Delta'\ge-1$, then $V_N/\E V_N\to1$ in probability; if $2u+\Delta'<-1$, then $V_N$ converges a.s.\ to a finite limit (\cref{lem:chisq}(a)); if $u<-1$, then $D_N$ converges a.s.\ to a finite limit $D_\infty$.

\emph{Step 2 (\cref{thm:window}(i)).} Let $u<-1$ and take $\alpha=0$, so $\hat\beta_N-\beta=M_N/D_N$ and $D_N\uparrow D_\infty<\infty$. If $2u+\Delta'<-1$, then $V_N\to V_\infty<\infty$, and given $X$ the independent increments $w_jx_jv_j$ have summable variances, so $M_N\to M_\infty$ a.s.\ by (C2), where $M_\infty$ is, given $X$, $N(0,V_\infty)$ with $V_\infty\ge w_1^2x_1^2g_1>0$. Then $\hat\beta_N-\beta\to M_\infty/D_\infty$ a.s., a variable whose conditional law given $X$ is non-degenerate normal, so $\Prob(M_\infty/D_\infty=0)=0$; if $\hat\beta_N$ converged in probability to $\beta$, a subsequence would converge a.s.\ to $\beta$, which is impossible. If $2u+\Delta'\ge-1$, then $V_N\to\infty$ in probability and, given $X$, $\Prob(\abs{M_N}\le KD_\infty\mid X)\le2KD_\infty/\sqrt{2\pi V_N}$, so $\abs{M_N}/D_N\ge\abs{M_N}/D_\infty\to\infty$ in probability.

\emph{Step 3 (\cref{thm:window}(ii)--(iii): reduction).} Let $u\ge-1$. Then $\E D_N\to\infty$ and $D_N\to\infty$ a.s.\ (\cref{lem:chisq}(a)), and $B_N/D_N\to0$ for every $\alpha\in\cA_0$ by \cref{lem:toeplitz}(a). Since $\alpha=0\in\cA_0$, consistency for every $\alpha\in\cA_0$ is equivalent to $M_N/D_N\to0$ in probability, hence, by \cref{lem:condgauss} with $\mathcal G=\sigma(X)$, to $V_N/D_N^2\to0$ in probability.

\emph{Step 4 (\cref{thm:window}(ii)).} Let $u>-1$, so $(\E D_N)^2\asymp N^{2u+2}$. If $2u+\Delta'>-1$, then $V_N\asymp_P\E V_N\asymp N^{2u+\Delta'+1}$, so $V_N/D_N^2\asymp_PN^{\Delta'-1}$, which tends to zero if and only if $\Delta'<1$; when $\Delta'\ge1$ it is bounded below in probability, and \cref{lem:condgauss} shows that the noise does not vanish. If $2u+\Delta'=-1$, then $V_N\asymp_P\log N$ and $V_N/D_N^2\to0$, while $\Delta'=-1-2u<1$. If $2u+\Delta'<-1$, then $V_N$ is bounded and $V_N/D_N^2\to0$, while $\Delta'<-1-2u<1$. In all cases, consistency holds if and only if $\Delta'<1$.

\emph{Step 5 (\cref{thm:window}(iii)).} Let $u=-1$, so $(\E D_N)^2\asymp(\log N)^2$ and $2u+\Delta'=\Delta'-2$. If $\Delta'<1$, then $\Delta'-2<-1$, $V_N$ is bounded and $V_N/D_N^2\to0$. If $\Delta'=1$, then $V_N\asymp_P\log N$ and $V_N/D_N^2\asymp_P1/\log N\to0$. If $\Delta'>1$, then $V_N\asymp_PN^{\Delta'-1}$ and $V_N/D_N^2\asymp_PN^{\Delta'-1}/(\log N)^2\to\infty$. Consistency holds if and only if $\Delta'\le1$.

\emph{Step 6 (Whittle--Mat\'ern form).} With $a=1+2\nu_X/d_s$, $b=1+2\nu_V/d_s$ and $\varkappa=1+2\nu_h/d_s$, $u=2(\nu_h-\nu_X)/d_s$ and $\Delta'=2(\nu_X-\nu_V)/d_s$. Thus $u>-1$ if and only if $\nu_h>\nu_X-d_s/2$, $u=-1$ if and only if $\nu_h=\nu_X-d_s/2$, and $\Delta'<1$ (respectively $\le1$) if and only if $\nu_X<\nu_V+d_s/2$ (respectively $\le$). Working ranges and variances enter only through the constants in \eqref{eq:powerlaw}. Finally, if $\alpha_j\to\alpha_\infty\in\R$, then $y_j=(\beta+\alpha_\infty)x_j+(\alpha_j-\alpha_\infty)x_j+v_j$: the data and the residual coefficients $v_j$ are exactly those of the parameter $(\beta+\alpha_\infty,\alpha-\alpha_\infty)$, whose confounding $\alpha-\alpha_\infty$ lies in $\cA_0$ (it belongs to the jointly Gaussian, $L$-stationary pair $(X,U-\alpha_\infty X)$). Under the conditions of (ii) or (iii), Steps 3--5 therefore give $\hat\beta_N\to\beta+\alpha_\infty$ in probability.

\emph{Step 7 (bounded weights; \cref{prop:estimators}(a)).} We first prove a general statement. Let a normal-form estimator have weights $\omega_j,\tilde\omega_j$, which may depend on the stage $N$, with $0\le\tilde\omega_j\le C$, $\abs{\omega_j}\le C$, and $\abs{\omega_{j_0}}\ge c>0$ and $\tilde\omega_{j_0}\ge c$ for some $j_0$ and all $N\ge j_0$, so that $\tilde D_N\ge c\,x_{j_0}^2>0$ a.s. Take $\alpha=0$. By \eqref{eq:normalform}, $\hat\beta_N-\beta=S_N+\tilde M_N/\tilde D_N$ with $S_N=\beta(\sum_{j\le N}\omega_jx_j^2/\tilde D_N-1)$, $\tilde D_N=\sum_{j\le N}\tilde\omega_jx_j^2$ and $\tilde M_N=\sum_{j\le N}\omega_jx_jv_j$; given $X$, $S_N$ and $\tilde D_N$ are fixed and $\tilde M_N\sim N(0,\tilde V_N)$ with $\tilde V_N=\sum_{j\le N}\omega_j^2x_j^2g_j$. Since $\tilde D_N\le C\norm X^2_{L^2}<\infty$ and $\tilde V_N\ge c^2x_{j_0}^2g_{j_0}$, we have $\tilde V_N/\tilde D_N^2\ge R_0=c^2x_{j_0}^2g_{j_0}/(C^2\norm X^4_{L^2})$, which is positive a.s.\ and does not depend on $N$, and \cref{lem:condgauss}(b) shows that the estimator is not consistent. If the weights do not depend on $N$, then $\tilde D_N$ and $\sum_{j\le N}\omega_jx_j^2$ converge a.s.\ to finite limits, the latter absolutely, and given $X$, $\tilde M_N$ converges a.s.\ by (C2), because $\sum_j\omega_j^2x_j^2g_j\le C^2\sum_jx_j^2g_j<\infty$; so the estimator converges a.s.\ to a limit whose conditional law given $X$ is normal with variance at least $R_0$. This covers GLS with a working nugget, where $\omega_j=\tilde\omega_j=1/(h_j+\sigma^2)\in(0,\sigma^{-2}]$, and the estimators of \cref{prop:estimators}(a), (b) and (e) with bounded weights. For OLS, $w\equiv1$, and $D_N\to\sum_jx_j^2=\norm X^2_{L^2}<\infty$ a.s., because $a>1$. The series $B_\infty=\sum_j\alpha_jx_j^2$ converges absolutely since $\alpha$ is bounded, and, given $X$, $M_N=\sum_{j\le N}x_jv_j$ converges a.s.\ by (C2), since $\E\sum_jx_j^2g_j=\sum_jf_jg_j<\infty$. The limit $\beta+B_\infty/\norm X^2+M_\infty/\norm X^2$ has, given $X$, a normal component with variance $\sum_jx_j^2g_j/\norm X^4>0$, so it is not a.s.\ equal to $\beta$, and OLS is inconsistent for every $\alpha$, including $\alpha=0$. Given $X$, its mean is $\beta+\sum_j\alpha_jx_j^2/\sum_jx_j^2$. By \cref{prop:normalform}, the RSR point estimate is the OLS estimate.

\emph{Step 8 (\cref{prop:estimators}(b)).} The same computation applies to the modes $j>K$, whose sum of squares $\sum_{j>K}x_j^2$ is finite and a.s.\ positive.

\emph{Step 9 (\cref{prop:estimators}(c)).} With weights $\ind\{K_N<j\le N\}$, the signal term is a weighted average of the $\alpha_j$ over $K_N<j\le N$, so its modulus is at most $\sup_{j>K_N}\abs{\alpha_j}\to0$. Write $K=K_N$, $S_1=\sum_{K<j\le N}x_j^2$ and $S_2=\sum_{K<j\le N}x_j^2g_j$. Given $X$, the noise is normal with variance $S_2/S_1^2$. For large $N$, $N\ge2K_N$, and \cref{lem:powersums} gives, since $a>1$ and $a+b>1$,
\begin{gather*}
\E S_1\asymp K^{1-a},\qquad\Var S_1\asymp K^{1-2a}=o\big(K^{2-2a}\big),\\
\E S_2\asymp K^{1-a-b},\qquad\Var S_2\asymp K^{1-2a-2b}=o\big(K^{2-2a-2b}\big).
\end{gather*}
By \cref{lem:chisq}(b), $S_2/S_1^2\asymp_PK^{1-a-b}/K^{2-2a}=K^{\Delta'-1}$. \Cref{lem:condgauss} shows that the estimator is consistent if and only if $\Delta'<1$. Moreover $\E[S_2/S_1^2]\asymp K^{\Delta'-1}$, so the unconditional variance of the noise term has the same order. For the lower bound, $S_2/S_1^2\asymp_PK^{\Delta'-1}$ gives $\E[S_2/S_1^2]\ge cK^{\Delta'-1}\Prob(S_2/S_1^2\ge cK^{\Delta'-1})\ge c'K^{\Delta'-1}$ for a suitable $c>0$ and all large $K$. For the upper bound, $S_1\ge\sum_{K<j\le2K}f_j\xi_j^2\ge c(2K)^{-a}\chi^2_K$ with $\chi^2_K=\sum_{K<j\le2K}\xi_j^2$, and, as in the proof of \cref{lem:inverse}, $\E(\chi^2_K)^{-4}\le(2/K)^4$ for $K\ge16$; hence $\E S_1^{-4}\le CK^{4(a-1)}$. Also $\E S_2^2=(\E S_2)^2+2\sum_{K<j\le N}(f_jg_j)^2\le3(\E S_2)^2\le CK^{2(1-a-b)}$. By the Cauchy--Schwarz inequality, $\E[S_2/S_1^2]\le(\E S_2^2)^{1/2}(\E S_1^{-4})^{1/2}\le CK^{1-a-b+2a-2}=CK^{\Delta'-1}$.

\emph{Step 10 (\cref{prop:estimators}(d)).} With weights $\ind\{K_N<j\le N\}/h_j$, the signal term is bounded by $\sup_{j>K_N}\abs{\alpha_j}\to0$. Since $c\le w_jg_j\le C$, the conditional variance of the noise lies between $c/S_N$ and $C/S_N$, where $S_N=\sum_{K_N<j\le N}w_jx_j^2$, and $c\,T_N\le S_N\le C\,T_N$ with $T_N=\sum_{K_N<j\le N}q_j\xi_j^2$. If $\sum_{K_N<j\le N}q_j\to\infty$, then $T_N\to\infty$ in probability by \cref{lem:chisq}(b), the conditional variance tends to zero in probability, and the estimator is consistent by \cref{lem:condgauss}. If $\sum_{K_N<j\le N}q_j$ stays bounded, $T_N$ is bounded in probability, the conditional variance is bounded below in probability, and the estimator is inconsistent when $\alpha=0$. If $\cI=\infty$, let $K_N=\max\{K:\sum_{j\le K}q_j\le\frac12\sum_{j\le N}q_j\}$. For each fixed $K$, eventually $\frac12\sum_{j\le N}q_j\ge\sum_{j\le K}q_j$, so $K_N\to\infty$; and $\sum_{K_N<j\le N}q_j\ge\frac12\sum_{j\le N}q_j\to\infty$. The same holds for any $K'_N\le K_N$ with $K'_N\to\infty$.

\emph{Step 11 (\cref{prop:estimators}(e)).} With a first-stage smoother $s_j\to0$ and second-stage GLS with working spectrum $h$, \cref{prop:normalform} gives the weights $\omega_j=(1-s_j)w_j$ and $\tilde\omega_j=(1-s_j)^2w_j$. For $j\ge j_0$, $1-s_j\ge1/2$, and the finitely many modes $j<j_0$ contribute a.s.\ finite amounts to every sum below. Signal: $\sum_j\omega_j\gamma_jx_j^2/\sum_j\tilde\omega_jx_j^2=\sum_j\tilde\omega_jx_j^2(\omega_j/\tilde\omega_j)\gamma_j/\sum_j\tilde\omega_jx_j^2$, and $(\omega_j/\tilde\omega_j)\gamma_j=\gamma_j/(1-s_j)\to\beta$ for $\alpha\in\cA_0$. If $u\ge-1$, then $\sum_j\tilde\omega_jx_j^2\ge\frac14\sum_{j\ge j_0}w_jx_j^2=\infty$ a.s., and \cref{lem:toeplitz}(a) gives the limit $\beta$. Noise: the conditional variance is $\tilde V_N/\tilde D_N^2$ with $\tilde V_N=\sum_{j\le N}\omega_j^2x_j^2g_j$ and $\tilde D_N=\sum_{j\le N}\tilde\omega_jx_j^2$; for $j\ge j_0$, $\omega_j^2\in[\frac14w_j^2,w_j^2]$ and $\tilde\omega_j\in[\frac14w_j,w_j]$. Hence $\tilde V_N/\tilde D_N^2$ tends to zero in probability if and only if $V_N/D_N^2$ does: when $D_N\to\infty$, the bounded contributions of $j<j_0$ are negligible in the denominators, and in the numerators they matter only when $V_N$ stays bounded, in which case both ratios tend to zero. If $u<-1$, take $\alpha=0$. Since $s_j\in[0,1]$, $\tilde\omega_j\le w_j$ and $\abs{\omega_j}\le w_j$, and $\sum_jw_jx_j^2<\infty$ a.s.; so $\tilde D_N\uparrow\tilde D_\infty<\infty$, the signal term $\beta\sum_{j\le N}\omega_jx_j^2/\tilde D_N$ converges a.s.\ to an $X$-measurable limit, and $\tilde V_N\ge\omega_{j_0}^2x_{j_0}^2g_{j_0}>0$, because $\omega_{j_0}\ge w_{j_0}/2>0$. \Cref{lem:condgauss}(b), with $R_0=\omega_{j_0}^2x_{j_0}^2g_{j_0}/\tilde D_\infty^2$, shows that the estimator is not consistent. So spatial+ of this form obeys the window of \cref{thm:window}. For the spatial+ of \cref{tab:zoo}, $s_j,t_j\in[0,1]$, so $0\le\tilde\omega_j=(1-t_j)(1-s_j)^2\le1$ and $\abs{\omega_j}\le1$; the same holds for gSEM, whose weights $(1-s_j)(1-t_j)$ and $(1-s_j)^2$ also lie in $[0,1]$. Whenever some $\omega_j\ne0$, Step 7 applies to both: numerator and denominator converge a.s., and the limit has a non-degenerate conditional normal component. \qed

\subsection{Proof of \texorpdfstring{\cref{cor:areal}}{Corollary 3}}\label{app:areal}

\emph{Idea.} Bounded weights leave a bounded denominator and a noise term whose conditional variance stays bounded away from zero, as for OLS. Weights that are powers of $j/n$ become powers of $j$ after a rescaling that does not change the estimator. The order computations in the proof of \cref{thm:window} use the weights only through two-sided bounds, and its almost sure arguments, which rely on weights that do not depend on $N$, are replaced by arguments in probability (Steps 3 and 4). The examples reduce to bounds on the eigenvalues of cycles and discrete tori, and part (c) is summation by parts on the graph.

\emph{Step 1 (part (a)).} Finitely many stages affect neither consistency nor limits, so we may assume that the bounds hold for all $n\ge2$; the same applies in part (b). Take $\alpha=0$, and let $D_n=\sum_{2\le j\le n}w^{(n)}_jx_j^2$, $M_n=\sum_{2\le j\le n}w^{(n)}_jx_jv_j$ and $V_n=\sum_{2\le j\le n}(w^{(n)}_j)^2x_j^2g_j$. The estimator is $\beta+M_n/D_n$, and given $X$, $M_n\sim N(0,V_n)$. Since $D_n\le p_+\norm X^2_{L^2}$ and $V_n\ge p_-^2x_2^2g_2$ for $n\ge2$, $V_n/D_n^2\ge p_-^2x_2^2g_2/(p_+^2\norm X^4_{L^2})$, a positive random variable that does not depend on $n$, and \cref{lem:condgauss}(b) shows that the estimator is not consistent. Now let $p_n=p$ be continuous at $0$ with $p(0)>0$, and let $\lambda^{(n)}_j\to0$ for each $j$, so that $w^{(n)}_j\to p(0)$ for each $j\ge2$. Since $\abs{w^{(n)}_jx_jy_j}\le p_+\abs{x_jy_j}$ and $\sum_j\abs{x_jy_j}\le\norm X_{L^2}\norm Y_{L^2}<\infty$ a.s., dominated convergence, applied to the terms multiplied by $\ind\{j\le n\}$, gives $\sum_{2\le j\le n}w^{(n)}_jx_jy_j\to p(0)\sum_{j\ge2}x_jy_j$ and $D_n\to p(0)\sum_{j\ge2}x_j^2>0$ a.s., for every $\alpha$; hence the estimator converges a.s.\ to $\sum_{j\ge2}x_jy_j/\sum_{j\ge2}x_j^2$, the limit of OLS with an intercept. For the examples: on a graph of maximal degree $k$ the eigenvalues of $L_n=D-W$ lie in $[0,2k]$ by Gershgorin's theorem, and on such graphs the working precisions of \cref{tab:zoo} are bounded above and below. The Leroux precision $\tau(1-\rho+\rho\lambda)$ lies in $[\tau(1-\rho),\tau(1-\rho+2k\rho)]$ for $0\le\rho<1$, and the graph Mat\'ern precision $\sigma^{-2}(2\nu/\kappa^2+\lambda)^\nu$ in $[\sigma^{-2}(2\nu/\kappa^2)^\nu,\sigma^{-2}(2\nu/\kappa^2+2k)^\nu]$. On $k$-regular graphs the proper CAR precision $\tau(k(1-\rho)+\rho\lambda)$ lies in $[\tau k(1-\abs\rho),\tau k(1+\abs\rho)]$ for $\abs\rho<1$, and the SAR precision $\sigma^{-2}(1-\rho k+\rho\lambda)^2$ in $[\sigma^{-2}(1-\abs\rho k)^2,\sigma^{-2}(1+\abs\rho k)^2]$ for $\abs\rho k<1$. All four precisions are continuous, with $p(0)$ equal to $\tau(1-\rho)$, $\sigma^{-2}(2\nu/\kappa^2)^\nu$, $\tau k(1-\rho)$ and $\sigma^{-2}(1-\rho k)^2$, respectively, all positive.

\emph{Step 2 (part (b): rescaling and orders).} Multiplying all weights by $n^{2m/d_s}$ does not change the estimator, and the rescaled weights $\tilde w^{(n)}_j=n^{2m/d_s}w^{(n)}_j$ satisfy $cj^\varkappa\le\tilde w^{(n)}_j\le Cj^\varkappa$ for $2\le j\le n$, with $\varkappa=2m/d_s$. Define $D_n$, $B_n=\sum_{2\le j\le n}\tilde w^{(n)}_j\alpha_jx_j^2$, $M_n$ and $V_n$ as in Step 1 of \cref{app:estimators}, with these weights and with sums over $2\le j\le n$, and let $u=\varkappa-a$. Then $\tilde w^{(n)}_jf_j$ and $(\tilde w^{(n)}_j)^2f_jg_j$ lie between constant multiples of $j^u$ and of $j^{2u+\Delta'}$, uniformly in $n$. \Cref{lem:chisq}(b), which allows weights that depend on $n$, and \cref{lem:powersums} therefore give the conclusions of Step 1 of \cref{app:estimators}: if $u\ge-1$, then $\E D_n\asymp\sum_{j\le n}j^u$ and $D_n/\E D_n\to1$ in probability; if $2u+\Delta'\ge-1$, then $\E V_n\asymp\sum_{j\le n}j^{2u+\Delta'}$ and $V_n/\E V_n\to1$ in probability; if $2u+\Delta'<-1$, then $V_n\le C^2\sum_jj^{2\varkappa}x_j^2g_j<\infty$ a.s.; and if $u<-1$, then $D_n\le\bar D=C\sum_jj^\varkappa x_j^2<\infty$ a.s.

\emph{Step 3 (part (b): $u<-1$).} Take $\alpha=0$. Then $D_n\le\bar D$ and $V_n\ge(\tilde w^{(n)}_2)^2x_2^2g_2\ge c^22^{2\varkappa}x_2^2g_2$, so $V_n/D_n^2$ is bounded below by a positive random variable that does not depend on $n$, and \cref{lem:condgauss}(b) shows that the estimator is not consistent. This is (i) of \cref{thm:window}.

\emph{Step 4 (part (b): $u\ge-1$).} Since $\E D_n\to\infty$ and $D_n/\E D_n\to1$ in probability, $D_n\to\infty$ in probability. For $\alpha\in\cA_0$ and $\epsilon>0$ choose $j_0$ with $\abs{\alpha_j}\le\epsilon$ for $j>j_0$; then $\abs{B_n}\le Cj_0^\varkappa\sup_j\abs{\alpha_j}\sum_{j\le j_0}x_j^2+\epsilon D_n$, so $B_n/D_n\to0$ in probability. As in Step 3 of \cref{app:estimators}, consistency for every $\alpha\in\cA_0$ is therefore equivalent to $V_n/D_n^2\to0$ in probability (\cref{lem:condgauss}(a)). Steps 4 and 5 of \cref{app:estimators} use only the orders established in Step 2, with ``$V_N$ is bounded'' read as ``$V_n$ is bounded by an a.s.\ finite random variable'', and give (ii) and (iii) of \cref{thm:window}.

\emph{Step 5 (part (b): Whittle--Mat\'ern form).} With $\varkappa=2m/d_s$ and $a=1+2\nu_X/d_s$, $u=2(m-\nu_X)/d_s-1$, so that $u>-1$, $u=-1$ or $u<-1$ according as $\nu_X<m$, $\nu_X=m$ or $\nu_X>m$; and $\Delta'=2(\nu_X-\nu_V)/d_s$. Steps 3 and 4 give the stated conditions.

\emph{Step 6 (part (b): the examples).} Let $c(j/n)^{2/d_s}\le\lambda^{(n)}_j\le C(j/n)^{2/d_s}$ for $2\le j\le n$ and $\Lambda=\sup_n\lambda^{(n)}_n<\infty$. Then $\tau\lambda$ satisfies the hypothesis with $m=1$, and so does $\lambda/(\sigma_u^2+\sigma_v^2\lambda)$, which lies between $\lambda/(\sigma_u^2+\sigma_v^2\Lambda)$ and $\lambda/\sigma_u^2$; $\tau\lambda^2$ satisfies it with $m=2$. Since $n^{-2/d_s}\le(j/n)^{2/d_s}\le\lambda^{(n)}_j/c$ for $j\ge2$: if $0\le1-\rho_n\le C'n^{-2/d_s}$ and $n$ is so large that $\rho_n\ge1/2$, the Leroux weights $\tau(1-\rho_n+\rho_n\lambda^{(n)}_j)$ lie between $\tau\lambda^{(n)}_j/2$ and $\tau(C'/c+1)\lambda^{(n)}_j$, and, on $k$-regular graphs, the proper CAR weights $\tau(k(1-\rho_n)+\rho_n\lambda^{(n)}_j)$ between $\tau\lambda^{(n)}_j/2$ and $\tau(kC'/c+1)\lambda^{(n)}_j$, so $m=1$; on $k$-regular graphs, if $0\le1-\rho_nk\le C'n^{-2/d_s}$ and $\rho_n\ge1/(2k)$, then $\lambda^{(n)}_j/(2k)\le1-\rho_nk+\rho_n\lambda^{(n)}_j\le(C'/c+1/k)\lambda^{(n)}_j$, so the SAR weights satisfy the hypothesis with $m=2$.

It remains to verify the eigenvalue bounds. The $n$-cycle has the eigenvalues $4\sin^2(\pi l/n)$, $l=0,\dots,n-1$; since $4\sin^2(\pi l/n)=4\sin^2(\pi(n-l)/n)$ increases in $l\le n/2$, the sorted eigenvalues are $\lambda^{(n)}_j=4\sin^2(\pi\lfloor j/2\rfloor/n)$. For $2\le j\le n$, $\theta=\pi\lfloor j/2\rfloor/n\in[\pi/n,\pi/2]$, and $2\theta/\pi\le\sin\theta\le\theta$ together with $j/3\le\lfloor j/2\rfloor\le j/2$ gives $(16/9)(j/n)^2\le\lambda^{(n)}_j\le\pi^2(j/n)^2$: $d_s=1$. The discrete torus $(\Z/q\Z)^d$, with $n=q^d$ nodes, has the eigenvalues $\mu(l)=\sum_{i\le d}4\sin^2(\pi l_i/q)$, with $l$ in a set of representatives satisfying $\abs{l_i}\le q/2$; as for the cycle, $16\abs l^2/q^2\le\mu(l)\le4\pi^2\abs l^2/q^2$. Let $N(r)$ be the number of representatives with $\abs l\le r$. For $r\ge1$, $N(r)\le(2r+1)^d\le(3r)^d$; for $1\le r\le q/4$ and $q\ge4$, the cube $\{\abs{l_i}\le\lfloor r/\sqrt d\rfloor\}$ consists of representatives and lies in the ball, so $N(r)\ge(r/\sqrt d)^d$. Lower bound: let $j\ge2$ and $r=j^{1/d}/4$. If $\lambda^{(n)}_j<16r^2/q^2$, at least $j$ eigenvalues are below $16r^2/q^2$, and each comes from an $l$ with $\abs l<r$; but there are at most $1<j$ such $l$ if $r<1$, and at most $(3r)^d<j$ if $r\ge1$. Hence $\lambda^{(n)}_j\ge16r^2/q^2=(j/n)^{2/d}$. Upper bound: if $r=\sqrt d\,j^{1/d}\le q/4$, at least $N(r)\ge j$ eigenvalues are at most $4\pi^2r^2/q^2$, so $\lambda^{(n)}_j\le4\pi^2d\,(j/n)^{2/d}$; otherwise $j/n>(4\sqrt d)^{-d}$ and $\lambda^{(n)}_j\le4d\le64d^2(j/n)^{2/d}$. So $d_s=d$.

\emph{Step 7 (part (c)).} For symmetric $W$ with $D_{ii}=\sum_{i'}W_{ii'}$,
\[
x^\top L_ny=\sum_i\sum_{i'}W_{ii'}x_i(y_i-y_{i'})=\tfrac12\sum_{i,i'}W_{ii'}(x_i-x_{i'})(y_i-y_{i'})=\sum_{i<i'}W_{ii'}(x_i-x_{i'})(y_i-y_{i'}),
\]
where the second equality averages the sum with the one obtained by exchanging $i$ and $i'$. The intrinsic CAR model has precision $\tau L_n$, whose null space contains the constants, so its GLS estimator, with or without an intercept, is $X^\top L_nY/X^\top L_nX$, which is the stated ratio; in the eigenbasis of $L_n$ it is the GLS estimator with weights $\lambda^{(n)}_j$. \qed

\subsection{Proof of \texorpdfstring{\cref{thm:minimax}}{Theorem 4}}\label{app:minimax}

\emph{Idea.} The upper bound computes the risk of oracle GLS: the bias is a weighted average of the $\alpha_j$, controlled by their decay, and the variance is an inverse moment of a weighted $\chi^2$ sum. The two lower bounds are two-point arguments. In the first, two parameters have the same association on every resolved mode, so their data laws coincide. In the second, the parameters are as far apart as the Cameron--Martin energy of the resolved modes allows while keeping their affinity bounded below.

\emph{Step 1 (upper bound: decomposition).} For oracle GLS, $w_j=1/g_j$, $w_jx_j^2=q_j\xi_j^2$ and $\hat\beta_N-\beta=B+\mathcal E$, where $B=\sum_{j\le N}q_j\xi_j^2\alpha_j/D_N$ is a function of $X$ and, given $X$, $\mathcal E=\sum_{j\le N}x_jv_j/g_j\big/D_N\sim N(0,1/D_N)$, with $D_N=\sum_{j\le N}q_j\xi_j^2$. Hence $\E\mathcal E^2=\E D_N^{-1}\le CN^{-p}$ by \cref{lem:inverse}, and $\E B\mathcal E=\E[B\,\E(\mathcal E\mid X)]=0$.

\emph{Step 2 (upper bound: bias).} On $\Theta(r,C_\alpha)$, $\abs B\le C_\alpha\sum_{j\le N}q_jj^{-r}\xi_j^2/D_N$. By the Cauchy--Schwarz inequality, Minkowski's inequality in $L^4$ with $\E\xi^8=105$, and \cref{lem:inverse} with $k=4$,
\[
\E B^2\le C_\alpha^2\Big\{\E\Big(\sum_{j\le N}q_jj^{-r}\xi_j^2\Big)^4\Big\}^{1/2}\big(\E D_N^{-4}\big)^{1/2}\le C_\alpha^2\sqrt{105}\Big(\sum_{j\le N}q_jj^{-r}\Big)^2C_4^{1/2}N^{-2p}.
\]
By \cref{lem:powersums}, $\sum_{j\le N}q_jj^{-r}\asymp\sum_{j\le N}j^{p-1-r}$ is of order $N^{p-r}$, $\log N$ or $1$ according as $r<p$, $r=p$ or $r>p$. In the three cases $\E B^2$ is of order $N^{-2r}$, $(\log N)^2N^{-2p}=o(N^{-p})$ and $N^{-2p}\le N^{-p}$. Hence $\E(\hat\beta_N-\beta)^2=\E B^2+\E\mathcal E^2\le C(N^{-2r}+N^{-p})$ for $N\ge N_0$, uniformly over $\Theta(r,C_\alpha)$.

\emph{Step 3 (lower bound: bias part).} Let $\Delta_N=C_\alpha N^{-r}$, $\theta_0=(0,0)$ and $\theta_1=(\Delta_N,\alpha^{(N)})$ with $\alpha^{(N)}_j=-\Delta_N$ for $j\le N$ and $0$ otherwise. Then $\abs{\alpha^{(N)}_j}\le C_\alpha j^{-r}$ for all $j$, so $\theta_1\in\Theta(r,C_\alpha)$, and $\gamma_j=0$ for $j\le N$ under both parameters. By Step 0 of \cref{app:frontier}, the laws of $(x_j,y_j)_{j\le N}$ under $\theta_0$ and $\theta_1$ coincide. For every $T_N$, since $t^2+(t-\Delta)^2\ge\Delta^2/2$,
\[
\E_0T_N^2+\E_1(T_N-\Delta_N)^2=\E_0\big[T_N^2+(T_N-\Delta_N)^2\big]\ge\Delta_N^2/2,
\]
so the maximum risk is at least $\Delta_N^2/4=C_\alpha^2N^{-2r}/4$, for every $N$.

\emph{Step 4 (lower bound: variance part).} Let $\cI_N=\sum_{j\le N}q_j\asymp N^p$, $\theta_0=(0,0)$ and $\theta_1=(\Delta,0)$ with $\Delta=\cI_N^{-1/2}$. By \eqref{eq:affinity} and $\log(1+t)\le t$, the affinity of the laws of the first $N$ modes is at least $\exp(-\Delta^2\cI_N/8)=e^{-1/8}$. By \cref{lem:distances}(c) and (a), $\max_i\Prob_i(\abs{T_N-\beta_i}\ge\Delta/2)\ge\frac12(1-\TV)\ge\frac14e^{-1/4}$, so by Markov's inequality $\max_i\E_i(T_N-\beta_i)^2\ge\frac{\Delta^2}4\cdot\frac14e^{-1/4}\ge cN^{-p}$.

\emph{Step 5 (conclusion).} The maximum of the two lower bounds is at least half their sum, which gives the left inequality in \eqref{eq:minimax}. \qed
\subsection{Proof of \texorpdfstring{\cref{thm:shift}}{Theorem 5} and \texorpdfstring{\cref{prop:shifts}}{Proposition 5}}\label{app:shift}

\emph{Idea.} Given $X$, the latent response is the law of $V$ translated by $\Psi_i(X)$. Two models therefore have equivalent conditional laws exactly when the difference of their translations is an admissible shift, and integrating over $X$ gives (a). The affinity of the joint laws is the average of the conditional affinities, which vanish on inadmissible shifts under the dichotomy property; this gives (b). For \cref{prop:shifts}, the Gaussian case is the Cameron--Martin theorem. Independent coefficients reduce, through Kakutani's theorem, to one-dimensional translation problems controlled by Fisher information, and scale mixtures are handled by conditioning on the mixing variable.

\emph{Step 1 (conditional laws).} For $a\in E$ let $P_V^a=P_V(\cdot-a)$, the law of $V+a$. Addition is continuous on the separable Fr\'echet space $E$, so $(x,v)\mapsto\Psi_i(x)+v$ is measurable, and since $X$ and $V$ are independent, $P_i=P_X\otimes K_i$ with the Markov kernel $K_i(x,\cdot)=P_V^{\Psi_i(x)}$ (measurable in $x$ by Fubini's theorem).

\emph{Step 2 (translations and the group property).} The translation $\tau_b:v\mapsto v-b$ is a Borel bijection of $E$ that maps $P_V^a$ to $P_V^{a-b}$. Hence, by \cref{lem:kernels}(e), $P_V^a\equivm P_V^b$ if and only if $P_V^{a-b}\equivm P_V$, that is, if and only if $a-b\in\cQ(P_V)$, and under the dichotomy property $P_V^a\perpm P_V^b$ otherwise. The set $\cQ(P_V)$ is an additive group: if $a,b\in\cQ$, then $P_V^{a+b}=\tau_{-b}(P_V^a)\equivm\tau_{-b}(P_V)=P_V^b\equivm P_V$, and $P_V^{-a}=\tau_a(P_V)\equivm\tau_a(P_V^a)=P_V$.

\emph{Step 3 (\cref{thm:shift}(a)).} If $\Psi_1(X)-\Psi_2(X)\in\cQ(P_V)$ a.s., then $K_1(x,\cdot)\equivm K_2(x,\cdot)$ for $P_X$-a.e.\ $x$ by Step 2, and \cref{lem:kernels}(a), applied in both directions, gives $P_1\equivm P_2$. An observation scheme is a Markov kernel $R$ from the space of $(X,\eta)$ to the space of data sequences $(D_1,D_2,\dots)$; \cref{lem:kernels}(b) gives $P_1R\equivm P_2R$, and \cref{lem:obstruction} shows that a functional taking different values in the two models is not consistently estimable.

\emph{Step 4 (\cref{thm:shift}(b)).} By \cref{lem:kernels}(c), $\aff(P_1,P_2)=\int\aff(P_V^{\Psi_1(x)},P_V^{\Psi_2(x)})\,P_X(dx)$. The integrand is at most one, and by Step 2 and \cref{lem:distances}(b) it vanishes on $N=\{x:\Psi_1(x)-\Psi_2(x)\notin\cQ(P_V)\}$. A measurable function bounded by $\ind_{N^c}$ has integral at most the outer probability of $N^c$, so $\aff(P_1,P_2)\le1-\Prob_\ast(X\in N)$, with $\Prob_\ast$ the inner probability (the ordinary probability when $N$ is measurable). By \cref{lem:distances}(a), $\TV(P_1,P_2)\ge H^2(P_1,P_2)=1-\aff(P_1,P_2)\ge\Prob_\ast(X\in N)$. If the shift is a.s.\ inadmissible, $\aff=0$ and $P_1\perpm P_2$ by \cref{lem:distances}(b).

\emph{Step 5 (\cref{prop:shifts}(a)).} This is (C5)(a): $\cQ(P_V)=\CM_V$, and every $h\notin\CM_V$ gives a singular translate. (C5) is stated for centred measures; if $P_V$ has mean $m$, the centred law $P_V(\cdot+m)$ has the same admissible and singular shifts, because translation by $m$ is a Borel bijection (\cref{lem:kernels}(e)).

\emph{Step 6 (\cref{prop:shifts}(b): reduction to coefficients).} Let $\ell_j$ be the coordinate functionals, $\ell_j=\inner{\cdot}{\varphi_j}$ when $E$ is continuously embedded in $L^2(M)$, or the $j$-th coordinate when $E$ is a sequence space. They are continuous, $\ell_j(\varphi_k)=\ind\{j=k\}$, and the map $\Pi:E\to\R^{\N}$, $\Pi(v)=(\ell_j(v))_j$, is continuous and injective. Since the series $V=\sum_j\sigma_je_j\varphi_j$ converges a.s.\ in $E$ and each $\ell_j$ is continuous, $\Pi(V)=(\sigma_je_j)_j$ and $\Pi(V+h)=(\sigma_je_j+h_j)_j$ with $h_j=\ell_j(h)$. Thus $\Pi$ maps $P_V$ to $\bigotimes_j\mathcal L(\sigma_je_j)$ and $P_V^h$ to $\bigotimes_j\mathcal L(\sigma_je_j+h_j)$, and by \cref{lem:kernels}(e) equivalence and singularity can be decided for these product measures.

\emph{Step 7 (\cref{prop:shifts}(b): one coordinate).} The $j$-th factors have the densities $\sigma_j^{-1}p_j(x/\sigma_j)$ and $\sigma_j^{-1}p_j((x-h_j)/\sigma_j)$, positive a.e., so they are equivalent. The affinity is invariant under a common rescaling, so their squared Hellinger distance is $H^2(p_j,p_j(\cdot-t_j))$ with $t_j=h_j/\sigma_j$. By \cref{lem:location}(a) and the hypothesis,
\[
c_\ast\min(t_j^2,1)\le H^2\big(p_j,p_j(\cdot-t_j)\big)\le\tfrac18\,t_j^2\sup_kJ(p_k),
\]
where $c_\ast>0$ is the infimum in the hypothesis.

\emph{Step 8 (\cref{prop:shifts}(b): conclusion).} If $\sum_jt_j^2<\infty$, the Hellinger distances are summable, and (C1) gives $P_V^h\equivm P_V$. If $\sum_jt_j^2=\infty$, then $\sum_j\min(t_j^2,1)=\infty$ (otherwise $t_j^2<1$ for large $j$ and $\sum_jt_j^2<\infty$), so the Hellinger distances are not summable and (C1) gives $P_V^h\perpm P_V$. Hence $\cQ(P_V)=\{h\in E:\sum_jh_j^2/\sigma_j^2<\infty\}$, and the dichotomy property holds.

\emph{Step 9 (\cref{prop:shifts}(b): identically distributed coefficients).} If the $e_j$ have a common density $p$ that is positive a.e.\ with $J(p)<\infty$, then $\sup_jJ(p_j)=J(p)$, and the infimum in the hypothesis is positive by \cref{lem:location}(c). The listed densities are positive everywhere, absolutely continuous and have finite Fisher information: the Laplace density $\tfrac12e^{-\abs x}$ has $(p')^2/p=p$ a.e., so $J=1$; the logistic and hyperbolic secant densities have bounded score $p'/p$; a finite Gaussian mixture $p=\sum_k\lambda_kp_k$, with weights $\lambda_k$ and component variances $s_k^2$, has $J(p)\le\sum_k\lambda_kJ(p_k)=\sum_k\lambda_k/s_k^2$, because the Cauchy--Schwarz inequality gives $(\sum_k\lambda_kp_k')^2\le(\sum_k\lambda_kp_k)(\sum_k\lambda_k(p_k')^2/p_k)$ pointwise; and Student's $t_\nu$ density has the bounded score $-(\nu+1)x/(\nu+x^2)$, of modulus at most $(\nu+1)/(2\sqrt\nu)$, so $J<\infty$ (in fact $J=(\nu+1)/(\nu+3)$, and $\nu=1$ gives $J=1/2$ for the Cauchy density).

\emph{Step 10 (\cref{prop:shifts}(c)).} Given $W=w$, $V=\sqrt wG$ has law $N(0,wC_G)$, whose Cameron--Martin space is $\CM_G$ as a set. If $h\in\CM_G$, then $N(h,wC_G)\equivm N(0,wC_G)$ for every $w>0$ by (C5)(a). The laws $P_V$ and $P_V^h$ are the images under $(w,v)\mapsto v$ of $\mathcal L(W)\otimes K$ and $\mathcal L(W)\otimes K_h$, with $K(w,\cdot)=N(0,wC_G)$ and $K_h(w,\cdot)=N(h,wC_G)$, so \cref{lem:kernels}(a) (in both directions) and (b) give $P_V^h\equivm P_V$. If $h\notin\CM_G$, the formula for $\abs h_{\CM}$ in (C5) yields $\ell_k\in E^\ast$ with $\E\ell_k(G)^2\le1$ and $\ell_k(h)\ge k^2$. If some $\ell_k$ has $\E\ell_k(G)^2=0$, then $A=\ker\ell_k$ has $P_V(A)=1$ and $P_V^h(A)=0$. Otherwise rescale so that $\E\ell_k(G)^2=1$, which only increases $\ell_k(h)$, and let $A=\{v:\abs{\ell_k(v)}\le k\text{ for all large }k\}$. Under $N(0,wC_G)$, $\ell_k(v)\sim N(0,w)$, so $\sum_k\Prob(\abs{\ell_k}>k)<\infty$ and $A$ has probability one by the Borel--Cantelli lemma. Under $N(h,wC_G)$, $\ell_k(v)=\ell_k(h)+N(0,w)$ with $\ell_k(h)\ge k^2$, so $\abs{\ell_k(v)}>k$ for all large $k$ almost surely, and $A$ has probability zero. Integrating over $W$, $P_V(A)=1$ and $P_V^h(A)=0$.

\emph{Step 11 (transfer to \cref{thm:frontier} and \cref{cor:matern}).} In the linear model with an $L$-stationary Gaussian $X$ independent of $V$, the models $(\beta,\alpha)$ and $(\beta',\alpha')$ differ by the shift $\delta(L)X$, with coefficients $\delta_jx_j$ and $\delta_j=\beta'-\beta+\alpha'_j-\alpha_j$. In cases (a), (b) and (c) of \cref{prop:shifts} (in case (c) with $G$ $L$-stationary), $\cQ(P_V)=\{h:\sum_jh_j^2/s_j<\infty\}$ for the scales $s_j=g_j$, $\sigma_j^2$ or the spectrum of $G$. By \cref{lem:chisq}(a), the shift is a.s.\ admissible if $\sum_j\delta_j^2f_j/s_j<\infty$ and a.s.\ inadmissible otherwise. \Cref{thm:shift} gives equivalence in the first case and, by the dichotomy property, singularity in the second. With $q_j=f_j/s_j$ this is \cref{thm:frontier}(i) (take $\delta\equiv\Delta$) and (iii), so the frontier $\cI=\infty$ is the same. For the consistency statement \cref{thm:frontier}(ii), Steps~6a--6b of \cref{app:frontier} involve only $X$, and Step~6c uses only that, given $X$, the $v_j$ are independent with mean zero and variance at most a constant multiple of $h_j$. In case (b) this holds when $\E e_j=0$ and $\sup_j\E e_j^2<\infty$, with $g_j=\sigma_j^2$; moreover $\Var e_j\ge1/J(p_j)\ge1/\sup_kJ(p_k)$ by \cref{lem:location}(b), so the coefficient variances are comparable to $\sigma_j^2$ and the conditional variance of the noise term is comparable to that of the Gaussian case, which gives the variance statement of \cref{cor:matern}. In case (c), conditionally on $W$ the $v_j$ are Gaussian with variances $Wg_j$, and the proofs apply conditionally on $W$; the variance statement then holds conditionally on $W$, and unconditionally when $\E W<\infty$, since the conditional variance is $W$ times its Gaussian value.

\emph{Step 12 (the Cauchy example after \cref{prop:shifts}).} Here $V$ is viewed in the sequence space $E=\R^{\N}$, which \cref{prop:shifts}(b) allows; as an $L^2(M)$-valued field, $\sum_j\sigma_je_j\varphi_j$ would converge a.s.\ only if $\sum_j\sigma_j<\infty$, by the three-series theorem \citep[Section~2.5]{Durrett2019}. Let the $e_j$ be independent standard Cauchy, $w_j=\sigma_j^{-2}$ and $\zeta_j=x_j/\sqrt{f_j}$. Given $X$, $\sum_{j\le N}w_jx_j\sigma_je_j=\sum_{j\le N}(x_j/\sigma_j)e_j$, and a linear combination $\sum_ja_je_j$ of independent standard Cauchy variables is Cauchy with scale $\sum_j\abs{a_j}$ (its characteristic function is $e^{-\abs t\sum_j\abs{a_j}}$). So the noise term of GLS is, given $X$, Cauchy with scale $s_N=\sum_{j\le N}\sqrt{q_j}\abs{\zeta_j}/D_N$, where $q_j=f_j/\sigma_j^2$ and $D_N=\sum_{j\le N}q_j\zeta_j^2$. Let $q_j\asymp j^{p-1}$ with $p>0$. Then $\E\sum_{j\le N}\sqrt{q_j}\abs{\zeta_j}=\sqrt{2/\pi}\sum_{j\le N}\sqrt{q_j}\asymp N^{(p+1)/2}$ with variance $(1-2/\pi)\sum_{j\le N}q_j\asymp N^p=o(N^{p+1})$, so this sum is $\asymp_PN^{(p+1)/2}$ by Chebyshev's inequality, and $D_N\asymp_PN^p$ by \cref{lem:chisq}(b) and \cref{lem:powersums}. Hence $s_N\asymp_PN^{(1-p)/2}$. Since $\Prob(\abs{\text{Cauchy}(s)}>\epsilon)=\frac2\pi\arctan(s/\epsilon)$, the noise term tends to zero in probability if $p>1$ and is bounded away from zero in probability if $p\le1$. The laws are singular for every $p\ge0$, because then $\sum_jq_j=\infty$ and the standard Cauchy density is positive with $J=1/2$ (Step 9), so \cref{prop:shifts}(b) applies. In Whittle--Mat\'ern terms $p=1+2(\nu_V-\nu_X)/d_s$, and $p>1$ means $\nu_X<\nu_V$. \qed

\subsection{Regularity of covariate paths}\label{app:regularity}

\begin{lemma}[Composition]\label{lem:composition}
Let $\Omega\subseteq\R^d$ be a bounded open set, $s>d/2$ and $F\in C^\infty(\R)$. If $u\in H^s(\Omega)$, then $F\circ u\in H^s(\Omega)$.
\end{lemma}
\begin{proof}
Write $u=U|_\Omega$ with $U\in H^s(\R^d)$. Since $s>d/2$, $\int\abs{\hat U}\le\big(\int(1+\abs\omega^2)^s\abs{\hat U}^2\big)^{1/2}\big(\int(1+\abs\omega^2)^{-s}\big)^{1/2}<\infty$, so $U$ is bounded and continuous. Let $G=F-F(0)$, so that $G\in C^\infty(\R)$ and $G(0)=0$; by (C6), $G\circ U\in H^s(\R^d)$. Choose $\chi\in C_c^\infty(\R^d)$ with $\chi=1$ on $\Omega$, which is possible because $\Omega$ is bounded. Then $F\circ u=(G\circ U+F(0)\chi)|_\Omega$ with $G\circ U+F(0)\chi\in H^s(\R^d)$, so $F\circ u\in H^s(\Omega)$.
\end{proof}

\begin{lemma}[Exact local regularity]\label{lem:regularity}
Let $G$ be a centred stationary Gaussian field on $\R^d$ with continuous paths and spectral density $f$ satisfying $c_1(1+\abs\omega^2)^{-(\nu+d/2)}\le f(\omega)\le c_2(1+\abs\omega^2)^{-(\nu+d/2)}$ for some $\nu>0$, and let $M\subset\R^d$ be a bounded connected Lipschitz domain. Then $X=G|_M$ satisfies (N2) with $s_X=\nu$, and so does $X=F\circ G$ for every $F\in C^\infty(\R)$ with $F'>0$.
\end{lemma}
\begin{proof}
\emph{Idea.} Localizing the field with a smooth cutoff produces a random function whose Fourier transform has variance comparable to $(1+\abs\omega^2)^{-(\nu+d/2)}$, so its $H^s$ norm has finite mean exactly when $s<\nu$. For $s<\nu$ this gives membership in $H^s$; for $s\ge\nu$, the zero--one law and Fernique's theorem turn an infinite mean into almost sure non-membership.

\emph{Step 1 (spectral bounds for localized fields).} Let $\chi\in C_c^\infty(\R^d)$, $\chi\not\equiv0$, and $\hat\chi(\omega)=\int\chi(s)e^{-is\cdot\omega}ds$. By (C8), $G(s)=\int e^{is\cdot\xi}Z_G(d\xi)$, and the isometry of (C8) and Fubini's theorem give $\widehat{\chi G}(\omega)=\int\hat\chi(\omega-\xi)Z_G(d\xi)$ for the Fourier transform of $\chi G$, so $\E\abs{\widehat{\chi G}(\omega)}^2=(\abs{\hat\chi}^2\ast f)(\omega)$. Since $\abs\omega^2\le2\abs{\omega-\xi}^2+2\abs\xi^2$, Peetre's inequality $1+\abs\omega^2\le2(1+\abs{\omega-\xi}^2)(1+\abs\xi^2)$ holds, and with the roles of $\omega$ and $\xi$ exchanged it gives, for $\sigma=\nu+d/2$,
\[
2^{-\sigma}(1+\abs\omega^2)^{-\sigma}(1+\abs{\omega-\xi}^2)^{-\sigma}\le(1+\abs\xi^2)^{-\sigma}\le2^{\sigma}(1+\abs\omega^2)^{-\sigma}(1+\abs{\omega-\xi}^2)^{\sigma}.
\]
Multiplying by $\abs{\hat\chi(\omega-\xi)}^2$, integrating in $\xi$ and using the bounds on $f$,
\[
c_12^{-\sigma}\kappa_-\,(1+\abs\omega^2)^{-\sigma}\le(\abs{\hat\chi}^2\ast f)(\omega)\le c_22^{\sigma}\kappa_+\,(1+\abs\omega^2)^{-\sigma},\qquad\kappa_\pm=\int\abs{\hat\chi(\eta)}^2(1+\abs\eta^2)^{\pm\sigma}d\eta,
\]
where $0<\kappa_\pm<\infty$ because $\hat\chi$ is a non-zero Schwartz function.

\emph{Step 2 (second moments).} By Tonelli's theorem, $\E\norm{\chi G}^2_{H^s}=\int(1+\abs\omega^2)^s\E\abs{\widehat{\chi G}(\omega)}^2d\omega$, which by Step 1 is finite if and only if $\int(1+\abs\omega^2)^{s-\nu-d/2}d\omega<\infty$, that is, if and only if $s<\nu$.

\emph{Step 3 ((ER$+$)).} Choose $\chi\in C_c^\infty(\R^d)$ with $\chi=1$ on $M$. For $s<\nu$, Step 2 gives $\chi G\in H^s(\R^d)$ a.s., hence $G|_M=(\chi G)|_M\in H^s(M)$ a.s. Taking $s$ along a sequence increasing to $\nu$ and using $H^s(M)\subseteq H^t(M)$ for $s\ge t$ gives one probability-one event on which $G|_M\in H^s(M)$ for all $s<\nu$.

\emph{Step 4 ((ER$-$) for one ball).} Fix a ball $B$ with $\bar B\subset M$ and $\chi_B\in C_c^\infty(B)$, $\chi_B\not\equiv0$. If $G|_B\in H^\nu(B)$, write $G|_B=U|_B$ with $U\in H^\nu(\R^d)$; since $\chi_B$ vanishes outside $B$, $\chi_BG=\chi_BU\in H^\nu(\R^d)$ by (C6). Now $\chi_BG$ is a centred Gaussian random element of the separable Hilbert space $L^2(\R^d)$, and $L_\nu=\{u\in L^2:\int(1+\abs\omega^2)^\nu\abs{\hat u}^2<\infty\}$ is a linear subspace that is Borel in $L^2$, since $u\mapsto\int(1+\abs\omega^2)^\nu\abs{\hat u}^2$ is the supremum of the continuous truncations $\int_{\abs\omega\le k}$. By (C5)(b), $\Prob(\chi_BG\in L_\nu)\in\{0,1\}$. Suppose this probability were one. The inclusion $H^\nu(\R^d)\to L^2(\R^d)$ is continuous and injective, so by (C3) the Borel sets of $H^\nu$ are Borel in $L^2$, and $\chi_BG$ is a random element of $H^\nu$. It is centred Gaussian there: each continuous linear functional $u\mapsto\int(1+\abs\omega^2)^\nu\hat u\,\overline{\hat h}$ on $H^\nu$ is the a.s.\ limit of its truncations to $\abs\omega\le k$, which are $L^2$-continuous and hence Gaussian, and a.s.\ limits of centred Gaussian variables are centred Gaussian. By (C5)(c), $\E\norm{\chi_BG}^2_{H^\nu}<\infty$, which contradicts Step 2 with $s=\nu$. Hence $G|_B\notin H^\nu(B)$ a.s., and therefore $G|_B\notin H^s(B)$ for every $s\ge\nu$.

\emph{Step 5 ((ER$-$) for all balls).} Let $\mathfrak B$ be the countable family of balls with rational centres and radii whose closures lie in $M$. Intersecting the events of Step 4 over $\mathfrak B$ gives a probability-one event on which $G|_{B'}\notin H^s(B')$ for all $B'\in\mathfrak B$ and $s\ge\nu$. Every open ball $B\subseteq M$ contains some $B'\in\mathfrak B$, and $G|_B\in H^s(B)$ would imply $G|_{B'}\in H^s(B')$ by restriction. So (ER$-$) holds.

\emph{Step 6 (the range).} $M$ is connected and $X$ is continuous, so $X(M)$ is an interval. For fixed $s\ne t$ in $M$, $X(s)-X(t)$ is centred Gaussian with variance $\int\abs{e^{is\cdot\omega}-e^{it\cdot\omega}}^2f(\omega)\,d\omega>0$, so $X(s)\ne X(t)$ a.s., and $X(M)$ is non-degenerate.

\emph{Step 7 (transforms).} Let $F\in C^\infty(\R)$ with $F'>0$. (ER$+$) for $F\circ G$ and $s\in(d/2,\nu)$ follows from Step 3 and \cref{lem:composition}. Suppose that $(F\circ G)|_B\in H^s(B)$ for some ball $B\subseteq M$ and some $s\ge\nu$ with $s>d/2$, and take $B'\in\mathfrak B$ with $B'\subseteq B$. The set $(F\circ G)(\bar B')$ is a compact subset of the open interval $F(\R)$, on which $F^{-1}$ is smooth; choose $\tilde F\in C^\infty(\R)$ equal to $F^{-1}$ on a neighbourhood of this compact set. Then $G|_{B'}=\tilde F\circ(F\circ G)|_{B'}\in H^s(B')$ by \cref{lem:composition}, contradicting Step 5. Finally $(F\circ G)(M)=F(G(M))$ is a non-degenerate interval.
\end{proof}

\subsection{Proof of \texorpdfstring{\cref{thm:nonlinear}}{Theorem 6}}\label{app:nonlinear}

\emph{Idea.} Two models differ by the shift $(m_1-m_2)\circ X$ of the latent response, up to Cameron--Martin-negligible terms. Where the difference $\delta=m_1-m_2$ has non-zero derivative, it is locally invertible, so on some small ball the covariate is a smooth function of $\delta\circ X$; if $\delta\circ X$ were an admissible shift, and hence $H^t$-smooth, the covariate would be $H^t$-smooth there, which exact local regularity excludes when $s_X\le t$. If instead the covariate is smoother than order $t$, every smooth function of it is an admissible shift. In part (c) the same local-inversion argument shows that the information matrix diverges in every direction.

\emph{Step 1 (part (a): the shift).} Let $\delta=m_1-m_2$. The shift between the two models is $\delta\circ X+(\mathsf A_1-\mathsf A_2)(X)$. Since $\mathsf A_1(X),\mathsf A_2(X)\in\cQ(P_V)$ a.s.\ and $\cQ(P_V)$ is a group (Step 2 of \cref{app:shift}), the shift is admissible if and only if $\delta\circ X\in\cQ(P_V)$.

\emph{Step 2 (part (a): a ball on which $\delta$ is invertible).} Let $\Omega_0$ be the probability-one event on which $X$ is continuous, $X(M)$ is a non-degenerate interval and (ER$-$) holds for every ball of the family $\mathfrak B$ of Step 5 of \cref{app:regularity}. Fix an outcome in $\Omega_0\cap\{\delta\text{ is non-constant on }X(M)\}$. If $\delta'$ vanished on the interior of $X(M)$, then $\delta$ would be constant on its closure; so there is $u$ in the interior of $X(M)$ with $\delta'(u)\ne0$, and an open interval $I\ni u$ with $\delta'\ne0$ on $\bar I$. Since $u$ is attained by $X$ and $X$ is continuous, $X^{-1}(I)$ is open and non-empty, so it contains a ball $B\in\mathfrak B$ with $X(\bar B)\subset I$. On $\bar I$, $\delta$ is a smooth bijection onto its image with smooth inverse; let $\tilde\delta\in C^\infty(\R)$ extend this inverse. Then $X|_B=\tilde\delta\circ(\delta\circ X)|_B$.

\emph{Step 3 (part (a): the shift is inadmissible).} If $\delta\circ X\in\cQ(P_V)\subseteq H^t_{\mathrm{loc}}(M)$, then $(\delta\circ X)|_B\in H^t(B)$, and \cref{lem:composition} with $t>d/2$ gives $X|_B\in H^t(B)$. As $t\ge s_X$ and $t>d/2$, this contradicts (ER$-$). So on $\Omega_0\cap\{\delta\text{ non-constant on }X(M)\}$, an event that is measurable because it can be written through the values of $X$ on a countable dense subset of $M$, the shift is inadmissible, and \cref{thm:shift}(b) gives $\TV\ge\Prob(\delta\text{ is non-constant on }X(M))$.

\emph{Step 4 (part (a): analytic differences).} A non-constant real-analytic $\delta$ is non-constant on every non-degenerate interval, so the event of Step 3 has probability one, $\TV=1$ and the laws are singular (\cref{lem:distances}(b)).

\emph{Step 5 (part (b)).} Choose $s\in(t,s_X)$; then $s>d/2$, and $X\in H^s(M)$ a.s.\ by (ER$+$). By \cref{lem:composition}, $\delta\circ X\in H^s(M)\subseteq H^t(M)\subseteq\cQ(P_V)$ for every $\delta\in C^\infty(\R)$. With a common $\mathsf A$ the shift is $\delta\circ X$, which is admissible a.s., and \cref{thm:shift}(a) applies. (By Step 1 the same holds for two different mechanisms in $\cA_{\mathrm{CM}}$.)

\emph{Step 6 (part (c): decomposition).} Write $a_j=\inner{\mathsf A(X)}{\varphi_j}$ and $z_j=(\inner{b_k(X)}{\varphi_j})_{k\le K}$, so that $y_j=z_j^\top\theta+a_j+v_j$ and, whenever $I_N$ is invertible,
\[
\hat\theta_N-\theta=I_N^{-1}\sum_{j\le N}\frac{z_ja_j}{g_j}+I_N^{-1}\sum_{j\le N}\frac{z_jv_j}{g_j}.
\]
\emph{Bias.} Let $b=\sum_{j\le N}z_ja_j/g_j$. For a unit vector $e$, the Cauchy--Schwarz inequality gives
\[
e^\top I_N^{-1/2}b=\sum_{j\le N}\frac{(e^\top I_N^{-1/2}z_j)\,a_j}{g_j}\le\Big(\sum_{j\le N}\frac{(e^\top I_N^{-1/2}z_j)^2}{g_j}\Big)^{1/2}\Big(\sum_{j\le N}\frac{a_j^2}{g_j}\Big)^{1/2}\le\norm{\mathsf A(X)}_{\CM_V},
\]
because the first factor is $(e^\top I_N^{-1/2}I_NI_N^{-1/2}e)^{1/2}=1$. Hence $\abs{I_N^{-1}b}\le\lambda_{\min}(I_N)^{-1/2}\norm{\mathsf A(X)}_{\CM_V}$, and $\norm{\mathsf A(X)}_{\CM_V}<\infty$ a.s.\ because $\mathsf A\in\cA_{\mathrm{CM}}$. \emph{Noise.} Given $X$, $\sum_{j\le N}z_jv_j/g_j\sim N(0,I_N)$, so the noise term is $N(0,I_N^{-1})$. It therefore suffices to show that $\lambda_{\min}(I_N)\to\infty$ a.s.

\emph{Step 7 (part (c): the information diverges in every direction).} For $c\in\R^K\setminus\{0\}$ let $\delta_c=\sum_kc_kb_k$. Then $c^\top I_Nc=\sum_{j\le N}\inner{\delta_c\circ X}{\varphi_j}^2/g_j$, which increases to $\norm{\delta_c\circ X}^2_{\CM_V}$ (infinite if $\delta_c\circ X\notin\CM_V$). On $\Omega_0$, every $\delta_c$ with $c\ne0$ is non-constant on the non-degenerate interval $X(M)$, by the assumed linear independence modulo constants; and $\CM_V\subseteq H^t_{\mathrm{loc}}(M)$ with $t\ge s_X$, so Steps 2--3 show $\delta_c\circ X\notin\CM_V$. Thus, on the single event $\Omega_0$, $c^\top I_Nc\to\infty$ for \emph{all} $c\ne0$ simultaneously. For $R>0$ the sets $U_N=\{c\in S^{K-1}:c^\top I_Nc>R\}$ are open, increase with $N$ (as $I_N$ increases in the positive semidefinite order) and cover the compact sphere $S^{K-1}$, so some $U_N$ equals $S^{K-1}$, that is, $\lambda_{\min}(I_N)>R$. Hence $\lambda_{\min}(I_N)\to\infty$ on $\Omega_0$, the bias tends to zero a.s., the noise tends to zero in probability, and $\hat\theta_N\to\theta$ in probability. \qed

\subsection{Proof of \texorpdfstring{\cref{thm:channels}}{Theorem 7} and \texorpdfstring{\cref{cor:kriging}}{Corollary 4}}\label{app:channels}

\emph{Idea.} (\cref{thm:channels}) Data obtained through a Markov kernel from equivalent laws have equivalent laws. If the coefficient estimates are accurate enough relative to the working spectrum, the sieve GLS estimator is a small perturbation of the complete-record estimator; the examples verify the required accuracy by moment bounds for averages of noise and of transformed counts. (\cref{cor:kriging}) A kriging predictor from a fixed network is a finite combination of kernel sections, each of which lies in the Cameron--Martin space of the residual when $2\nu_\ast>\nu_V$. For the structural effect, the response given the station values is a Gaussian element whose mean is a finite combination of fixed functions; that mean is an admissible shift of the conditional noise, so reversing the sign of the effect, or rescaling it together with the covariate variance, leaves the law unchanged up to equivalence.

\begin{lemma}[Conditioning on point values]\label{lem:conditioning}
In the Gaussian model of \cref{sec:confounding}, let the eigenfunctions be continuous, let $m_1,\dots,m_K\in M$ satisfy $\sum_jf_j\varphi_j(m_i)^2<\infty$, put $k_i=\sum_jf_j\varphi_j(m_i)\varphi_j$ and $Z_i=X(m_i)=\sum_jx_j\varphi_j(m_i)$ (an $L^2(\Prob)$-convergent series), and suppose that the covariance matrix $K_Z$ of $Z=(Z_1,\dots,Z_K)$ is nonsingular.
\begin{enumerate}[label=(\alph*),itemsep=1pt]
\item $k_i\in\CM_X$, and for $h=\sum_jh_j\varphi_j\in\CM_X$ the series $h(m_i)=\sum_jh_j\varphi_j(m_i)$ converges absolutely, with $h(m_i)=\inner h{k_i}_{\CM_X}$.
\item $m_Z=\E[X\mid Z]=\sum_i(K_Z^{-1}Z)_ik_i$, and $R=X-m_Z$ is a centred Gaussian element of $L^2(M)$ independent of $Z$.
\item For every $h\in\CM_X$ with $h(m_1)=\dots=h(m_K)=0$, the laws of $R+h$ and $R$ are equivalent.
\item If $X$ has continuous paths and $\operatorname{supp}\mu=M$, then $Z_i$ equals the point value $X(m_i)$ of the continuous version almost surely.
\end{enumerate}
\end{lemma}
\begin{proof}
(a) $\norm{k_i}^2_{\CM_X}=\sum_j(f_j\varphi_j(m_i))^2/f_j=\sum_jf_j\varphi_j(m_i)^2<\infty$. For $h\in\CM_X$, $\sum_j\abs{h_j\varphi_j(m_i)}\le(\sum_jh_j^2/f_j)^{1/2}(\sum_jf_j\varphi_j(m_i)^2)^{1/2}$ and $\inner h{k_i}_{\CM_X}=\sum_jh_jf_j\varphi_j(m_i)/f_j=h(m_i)$.

(b) For $h\in\CM_X$ let $\mathbb W(h)=\sum_jh_j\xi_j/\sqrt{f_j}$, an $L^2(\Prob)$-convergent series with law $N(0,\norm h^2_{\CM_X})$; $h\mapsto\mathbb W(h)$ is a linear isometry from $\CM_X$ onto a centred Gaussian subspace of $L^2(\Prob)$. Then $Z_i=\mathbb W(k_i)$ and $\Cov(Z_i,Z_l)=\inner{k_i}{k_l}_{\CM_X}$, so $K_Z$ is the Gram matrix of $k_1,\dots,k_K$, which are therefore linearly independent. Let $e_1,\dots,e_K$ be an orthonormal basis of $\operatorname{span}\{k_i\}$ and $(e_k)_{k>K}$ an orthonormal basis of its orthogonal complement $S$ in the separable Hilbert space $\CM_X$. The variables $\zeta_k=\mathbb W(e_k)$ are independent standard normal. Since $\inner h{f_j\varphi_j}_{\CM_X}=h_j$, the coefficient $x_j=\sqrt{f_j}\,\xi_j=\mathbb W(f_j\varphi_j)$ equals $\sum_k(e_k)_j\zeta_k$, with $(e_k)_j=\inner{e_k}{\varphi_j}=\inner{f_j\varphi_j}{e_k}_{\CM_X}$; the series converges in $L^2(\Prob)$ and, by (C2), almost surely. Let $m_Z$ have the coefficients $\sum_{k\le K}(e_k)_j\zeta_k$, that is $m_Z=\sum_{k\le K}\zeta_ke_k$, and $R=X-m_Z$, with coefficients $r_j=\sum_{k>K}(e_k)_j\zeta_k$. Since $\operatorname{span}\{\zeta_1,\dots,\zeta_K\}=\mathbb W(\operatorname{span}\{k_i\})=\operatorname{span}\{Z_1,\dots,Z_K\}$ and $(\zeta_k)_{k>K}$ is independent of $(\zeta_k)_{k\le K}$, $R$ is independent of $Z$ and $m_Z$ is a function of $Z$, so $m_Z=\E[X\mid Z]$. As $m_Z\in\operatorname{span}\{k_i\}$ and, by (a), $m_Z(m_l)=\sum_{k\le K}\zeta_k\inner{e_k}{k_l}_{\CM_X}=\mathbb W(k_l)=Z_l$, we have $m_Z=\sum_i(K_Z^{-1}Z)_ik_i$.

(c) Let $h\in\CM_X$ vanish at $m_1,\dots,m_K$. By (a), $h\perp k_i$ in $\CM_X$ for all $i$, so $h=\sum_{k>K}h_ke_k$ with $h_k=\inner h{e_k}_{\CM_X}$ and $\sum_kh_k^2<\infty$. Its coefficients are $h_j=\sum_{k>K}h_k(e_k)_j$, absolutely convergent because $\sum_k(e_k)_j^2\le\norm{f_j\varphi_j}^2_{\CM_X}=f_j$. Let $\Phi$ map a real sequence $(c_k)_{k>K}$ to the sequence of limits $\lim_n\sum_{K<k\le n}(e_k)_jc_k$, $j\ge1$, where they exist, and to $0$ otherwise; $\Phi$ is Borel. Then the coefficient sequence of $R$ is $\Phi((\zeta_k)_{k>K})$ a.s., and that of $R+h$ is $\Phi((\zeta_k+h_k)_{k>K})$ a.s. By \cref{lem:gaussaff}(a), $H^2(N(0,1),N(t,1))=1-e^{-t^2/8}\le t^2/8$, so (C1) and $\sum_kh_k^2<\infty$ give $\mathcal L((\zeta_k+h_k)_k)\equivm\mathcal L((\zeta_k)_k)$. By \cref{lem:kernels}(b) the coefficient sequences of $R+h$ and $R$ have equivalent laws, and \cref{lem:kernels}(e), applied to the continuous injective coefficient map $L^2(M)\to\R^{\N}$, gives the claim.

(d) Let $X$ have continuous paths and $C(s,t)=\E X(s)X(t)$. By (C5)(c) in $E=C(M)$, $\E\sup_M\abs X^2<\infty$, so $C$ is continuous by dominated convergence, and Fubini's theorem gives $\E[X(m)x_j]=\int C(m,t)\varphi_j(t)\,\mu(dt)$. The covariance operator of $X$ on $L^2(M)$ has the eigenpairs $(f_j,\varphi_j)$, so this continuous function of $m$ equals $f_j\varphi_j(m)$ for $\mu$-a.e.\ $m$, hence for every $m\in\operatorname{supp}\mu=M$, the eigenfunctions being continuous. By (C7), $C(m,m)=\sum_jf_j\varphi_j(m)^2$; the expansion of (C7) does not depend on the choice of orthonormal bases within the eigenspaces, and the terms are non-negative. Therefore $\E\abs{X(m)-\sum_{j\le n}x_j\varphi_j(m)}^2=C(m,m)-\sum_{j\le n}f_j\varphi_j(m)^2\to0$, and the $L^2(\Prob)$-limit $Z_i$ equals $X(m_i)$ almost surely.
\end{proof}

\emph{Step 1 (\cref{thm:channels}(a)).} Each cited non-identification statement establishes the equivalence, or the equality, of the laws of one set of observables under two parameter values with different values of the functional of interest: of $(X,\eta)$ for \cref{thm:frontier}(i),(iii), \cref{thm:shift}(a) and \cref{thm:nonlinear}(b); of $(\hat X,Y)$ for \cref{cor:kriging}(a); of $(Z,Y)$ for \cref{cor:kriging}(b); and of $D_\delta$ or $B_\delta$ for \cref{thm:resolution}(a),(b). Data generated through a common Markov kernel inherit equivalence by \cref{lem:kernels}(b) (and equality trivially), and \cref{lem:obstruction} applies.

\emph{Step 2 (\cref{thm:channels}(b): the general criterion).} Let $e_j=\abs{\hat x_j-x_j}+\abs{\hat\eta_j-\eta_j}$ and $E_J=\sum_{j\le J}e_j/h_j$. By hypothesis $E_J/(\epsilon_na_J)$ is tight uniformly in $n$ and $J$, so $E_{J_n}=O_P(\epsilon_na_{J_n})=O_P(\epsilon_n^{1/2})\to0$ in probability, because $a_{J_n}\le\epsilon_n^{-1/2}$. Moreover $J_n\to\infty$, since each $a_J$ is finite and $\epsilon_n^{-1/2}\to\infty$. Step~4 of \cref{app:pointsamples} used only $E_{J_n}\to0$ in probability, $J_n\to\infty$ and the almost sure properties of the complete-record GLS estimator, so it applies verbatim and gives consistency.

\emph{Step 3 (examples (i) and (ii)).} For (i), Step~3 of \cref{app:pointsamples} gives $E_J\le a_J\{R_1\ell_n^\vartheta+\Xi_{n,J}\}$, where $\Xi_{n,J}=a_J^{-1}\sum_{j\le J}\abs{\zeta_j}/h_j$ satisfies $\E\,\Xi_{n,J}\le\sigma\mu(M)^{1/2}(\max_iw_i)^{1/2}$ uniformly in $J$. With $\epsilon_n=1/\log(1/\varpi_n)$, $\ell_n^\vartheta/\epsilon_n\le\varpi_n^\vartheta\log(1/\varpi_n)$ and $(\max_iw_i)^{1/2}/\epsilon_n\le\varpi_n^{1/2}\log(1/\varpi_n)$ are bounded, so $E_J/(\epsilon_na_J)$ is tight uniformly in $n$ and $J$ by Markov's inequality. For (ii), classical measurement error $X(s_i)+\varepsilon^X_i$, with independent errors of mean zero and variance $\tau^2$, adds to $\hat x_j$ the term $\sum_iw_i\varepsilon^X_i\varphi_j(s_i)$, whose second moment is at most $\tau^2B_j^2\mu(M)\max_iw_i$; it is handled exactly like $\zeta_j$.

\emph{Step 4 (example (iii): Poisson counts).} Given $\eta$, the counts $N_i$, $i\le n$, are independent Poisson with means $\mu_i=\Lambda_n\int_{A_i}\rho e^\eta\,d\mu$, where $\rho$ is known with $0<\rho_-\le\rho\le\rho_+<\infty$, and $\eta$ has H\"older paths of order $\vartheta$ with constant $K_\eta$. Let $\hat\eta_i=\log\{(N_i+\frac12)/(\Lambda_n\int_{A_i}\rho\,d\mu)\}$ and $\hat\eta_j=\sum_iw_i\hat\eta_i\varphi_j(s_i)$. Write $\hat\eta_i-\eta(s_i)=Z_i+b_i$ with $Z_i=\log\{(N_i+\frac12)/\mu_i\}$ and $b_i=\log\{\mu_i/(\Lambda_n\int_{A_i}\rho\,d\mu)\}-\eta(s_i)$. The ratio $\mu_i/(\Lambda_n\int_{A_i}\rho)$ is a $\rho$-weighted average of $e^\eta$ over $A_i$, so it lies between $e^{\min_{A_i}\eta}$ and $e^{\max_{A_i}\eta}$, and $\abs{b_i}\le\operatorname{osc}_{A_i}\eta\le K_\eta\ell_n^\vartheta$. Moreover $\mu_i\ge c_\eta\Lambda_nw_i$ with $c_\eta=\rho_-e^{-\norm\eta_\infty}$, so the event $\{\min_i\mu_i\ge1\}$ has probability tending to one because $\Lambda_n\min_iw_i\to\infty$. On it, \cref{lem:counts} gives $\abs{\E[Z_i\mid\eta]}\le C/\mu_i$ and $\E[Z_i^2\mid\eta]\le C/\mu_i$, the $Z_i$ being conditionally independent. Hence
\[
\abs{\hat\eta_j-\eta_j}\le\Big|\sum_iw_i\eta(s_i)\varphi_j(s_i)-\eta_j\Big|+B_j\mu(M)K_\eta\ell_n^\vartheta+\Big|\sum_iw_iZ_i\varphi_j(s_i)\Big|,
\]
where the first term is the quadrature error of Step~3 of \cref{app:pointsamples}, and the last term has conditional expectation at most
\[
B_j\sum_iw_i\frac C{\mu_i}+B_j\Big(\sum_iw_i^2\frac C{\mu_i}\Big)^{1/2}\le CB_j\Big\{\frac{\mu(M)}{c_\eta\Lambda_n\min_iw_i}+\Big(\frac{\mu(M)}{c_\eta\Lambda_n}\Big)^{1/2}\Big\}.
\]
With $\epsilon_n=1/\log(1/\varpi_n)+(\Lambda_n\min_iw_i)^{-1}+\Lambda_n^{-1/2}$, and the a.s.\ finite random constants $\norm\eta_\infty$, $K_\eta$ and $c_\eta^{-1}$, the criterion of Step 2 holds, by Markov's inequality applied conditionally on $\eta$. Note that $\Lambda_n\min_iw_i\to\infty$ implies $\Lambda_n\to\infty$.

\emph{Step 5 (example (iv): binomial counts).} Given $\eta$, let $N_i$ be independent binomial with $m_i$ trials and success probabilities $p_i=1/(1+e^{-\eta(s_i)})$, with $\min_im_i\to\infty$, and let $\hat\eta_i=\log\{(N_i+\frac12)/(m_i-N_i+\frac12)\}$. Since $\eta(s_i)=\log(m_ip_i)-\log(m_i(1-p_i))$, $\hat\eta_i-\eta(s_i)=Z_i-Z_i'$ with $Z_i=\log\{(N_i+\frac12)/(m_ip_i)\}$ and $Z'_i=\log\{(m_i-N_i+\frac12)/(m_i(1-p_i))\}$. Both $N_i$ and $m_i-N_i$ are binomial, with means at least $c_\eta m_i$, $c_\eta=(1+e^{\norm\eta_\infty})^{-1}$. By \cref{lem:counts}, once $c_\eta\min_im_i\ge1$, $\abs{\E[\hat\eta_i-\eta(s_i)\mid\eta]}\le2C/(c_\eta m_i)$ and $\E[(\hat\eta_i-\eta(s_i))^2\mid\eta]\le4C/(c_\eta m_i)$. The argument of Step 4, with $b_i=0$, gives the criterion with $\epsilon_n=1/\log(1/\varpi_n)+(\min_im_i)^{-1/2}$.

\emph{Step 6 (\cref{cor:kriging}(a): regression on a predicted surface).} If $\hat X\in\cQ(P_V)$ a.s., then $(\beta-\beta')\hat X\in\cQ(P_V)$ a.s.\ because $\cQ(P_V)$ is a linear space. The working models $(\beta,\mathsf A)$ and $(\beta',\mathsf A)$ for $Y=\beta\hat X+\mathsf A(\hat X)+V$ differ by this shift, and \cref{thm:shift}(a), with $\hat X$ in the role of $X$, applies.

\emph{Step 7 (\cref{cor:kriging}(a): kernel sections lie in $\CM_V$).} For Gaussian $V$, $\cQ(P_V)=\CM_V$ is linear, and it suffices that each kernel section $k_\ast(\cdot,m)=\sum_jf^\ast_j\varphi_j(m)\varphi_j$ lies in $\CM_V$, where $f^\ast_j=\tau_\ast^2(\kappa_\ast^2+\lambda_j)^{-(\nu_\ast+d_s/2)}$; then $\hat X=\sum_{i\le K}c_ik_\ast(\cdot,m_i)\in\CM_V$ for every value of the weights. Let $N_m(\Lambda)=\sum_{\lambda_j\le\Lambda}\varphi_j(m)^2$. Since $p_\tau(m,m)=\sum_je^{-\tau\lambda_j}\varphi_j(m)^2\ge e^{-1}N_m(1/\tau)$, the heat-kernel bound gives $N_m(\Lambda)\le e\,p_{1/\Lambda}(m,m)\le C\Lambda^{d_s/2}$ for $\Lambda\ge1$. Since $g_j\ge c(1+\lambda_j)^{-(\nu_V+d_s/2)}$, $(f^\ast_j)^2/g_j\le C\tilde F(\lambda_j)$ with the smooth majorant $\tilde F(\lambda)=(\kappa_\ast^2+\lambda)^{-2\nu_\ast-d_s}(1+\lambda)^{\nu_V+d_s/2}$. For $\lambda\ge1$, $\tilde F(\lambda)\le C\lambda^{\nu_V-2\nu_\ast-d_s/2}$ and, from $\tilde F'=\tilde F\,\{-(2\nu_\ast+d_s)/(\kappa_\ast^2+\lambda)+(\nu_V+d_s/2)/(1+\lambda)\}$, $\abs{\tilde F'(\lambda)}\le C\lambda^{\nu_V-2\nu_\ast-d_s/2-1}$. Summation by parts for the Stieltjes integral against the non-decreasing step function $N_m$ gives
\begin{align*}
\sum_{1\le\lambda_j\le\Lambda}\tilde F(\lambda_j)\varphi_j(m)^2&=\tilde F(\Lambda)N_m(\Lambda)-\tilde F(1)N_m(1^-)-\int_1^\Lambda N_m(\lambda)\tilde F'(\lambda)\,d\lambda\\
&\le C\Lambda^{\nu_V-2\nu_\ast}+C\int_1^\infty\lambda^{\nu_V-2\nu_\ast-1}\,d\lambda,
\end{align*}
which is bounded uniformly in $\Lambda\ge1$ whenever $2\nu_\ast>\nu_V$. Only finitely many $\lambda_j$ are below $1$, so $\norm{k_\ast(\cdot,m)}^2_{\CM_V}=\sum_j(f^\ast_j)^2\varphi_j(m)^2/g_j<\infty$.

\emph{Step 8 (\cref{cor:kriging}(b): the conditional law of the response).} Let $\gamma(L)=\beta+\alpha(L)$, a bounded operator because $\alpha$ is bounded, so that $Y=\gamma(L)X+V$. With $m_Z$ and $R$ from \cref{lem:conditioning}, $Y=\gamma(L)m_Z+W$ with $W=\gamma(L)R+V$. The pair $(R,V)$ is independent of $Z$, so the conditional law of $Y$ given $Z=z$ is the law of $\gamma(L)m_z+W$, where $m_z=\sum_i(K_Z^{-1}z)_ik_i$.

\emph{Step 9 (\cref{cor:kriging}(b): admissible shifts of $W$).} For $J$ large, evaluation at $m_1,\dots,m_K$ maps $\operatorname{span}\{\varphi_1,\dots,\varphi_J\}$ onto $\R^K$. Otherwise, since the images increase with $J$, some $a\ne0$ would satisfy $\sum_la_l\varphi_j(m_l)=0$ for all $j$, and then $a^\top K_Za=\Var(\sum_la_lX(m_l))=\sum_jf_j(\sum_la_l\varphi_j(m_l))^2=0$, contradicting the nonsingularity of $K_Z$. So there are $u_i$ in this span with $u_i(m_l)=k_i(m_l)$ for all $l$. Then $k_i-u_i\in\CM_X$ vanishes at the stations. For $z\in\R^K$ and $t\in\R$ write $t\gamma(L)m_z=\gamma(L)h+v_s$ with $h=t\sum_ic_i(k_i-u_i)$ and $v_s=t\sum_ic_i\gamma(L)u_i$, $c=K_Z^{-1}z$. By \cref{lem:conditioning}(c), $\mathcal L(R+h)\equivm\mathcal L(R)$. The function $v_s$ is a finite combination of $\varphi_1,\dots,\varphi_J$, so it shifts finitely many independent Gaussian coefficients of $V$, and $\mathcal L(V+v_s)\equivm\mathcal L(V)$ by (C1) and \cref{lem:kernels}(e). By independence and \cref{lem:kernels}(d), $\mathcal L(R+h,V+v_s)\equivm\mathcal L(R,V)$, and applying the continuous map $(r,v)\mapsto\gamma(L)r+v$ (\cref{lem:kernels}(b)) gives $\mathcal L(W+t\gamma(L)m_z)\equivm\mathcal L(W)$, for every $z$ and $t$.

\emph{Step 10 (\cref{cor:kriging}(b): the sign).} Under $(-\beta,-\alpha)$, $\gamma(L)$ becomes $-\gamma(L)$, and the law of $X$, hence of $(Z,R)$, is unchanged. The new noise $-\gamma(L)R+V$ has the law of $W$, because $(-R,V)$ has the law of $(R,V)$. So the two conditional laws of $Y$ given $Z=z$ are $\mathcal L(W-\gamma(L)m_z)$ and $\mathcal L(W+\gamma(L)m_z)$, both equivalent to $\mathcal L(W)$ by Step 9 with $t=\mp1$. \Cref{lem:kernels}(a), applied in both directions with the common law of $Z$, gives equivalence of the laws of $(Z,Y)$; \cref{lem:kernels}(b) extends this to any Markov kernel acting on the response, and \cref{lem:obstruction} shows that the sign of $\beta\ne0$ is not consistently estimable over any confounding class closed under $\alpha\mapsto-\alpha$.

\emph{Step 11 (\cref{cor:kriging}(b): the scale).} Under $(\beta',\alpha\beta'/\beta,\varsigma f)$ with $\varsigma=(\beta/\beta')^2$, $\gamma(L)$ becomes $(\beta'/\beta)\gamma(L)$ and $X$ has the law of $\sqrt\varsigma\,X$. Then $Z\sim N(0,\varsigma K_Z)$, which is equivalent to $N(0,K_Z)$ since both are nonsingular Gaussian laws on $\R^K$. The kernels $k_i$ and $K_Z$ are multiplied by $\varsigma$, so $m_z=\sum_i(K_Z^{-1}z)_ik_i$ is unchanged as a function of $z$, while $R$ is replaced by a copy of $\sqrt\varsigma\,R$. The new noise $(\beta'/\beta)\sqrt\varsigma\,\gamma(L)R+V=\pm\gamma(L)R+V$ has the law of $W$. So the conditional law of $Y$ given $Z=z$ is $\mathcal L(W+(\beta'/\beta)\gamma(L)m_z)$, equivalent to $\mathcal L(W)$ and hence to $\mathcal L(W+\gamma(L)m_z)$ by Step 9. With equivalent marginal laws of $Z$ and equivalent conditional laws, \cref{lem:kernels}(a) gives equivalent joint laws, and \cref{lem:kernels}(b) any observation of the response. Finally, $\alpha\beta'/\beta\in\cA_0$ whenever $\alpha\in\cA_0$. None of these steps used the residual spectrum or the sampling scheme of the response. \qed
\section{Simulation details}\label{app:sims}

All experiments use $\beta=1$ unless stated otherwise and are fully reproducible from the supplementary scripts (fixed random seeds). Convergence exponents are least-squares slopes of $\log$ RMSE on $\log N$ (or $\log n$) over the upper half of the logarithmic range.

\paragraph{E1 (spectral model).} Eigenvalues $\lambda_j=j^{2/d_s}$, $j\le N_{\max}=2^{17}$ ($2^{18}$ in panel c); $f_j=(1+\lambda_j)^{-(\nu_X+d_s/2)}$, $g_j=(1+\lambda_j)^{-(\nu_V+d_s/2)}$ with $\nu_V=1/2$; confounding $\alpha_j=0.5(1+\lambda_j)^{-\varrho}$ with $\varrho=1$ in panels (a)--(b) and $\varrho=0.75$ in \cref{tab:estimators}, which uses $d_s=1$, $\nu_X=0.75$ (inside) or $1.125$ (outside), and working smoothness $0.05$ (undersmoothed) or $\nu_V+1$ (oversmoothed); in panel (c) $\alpha_j=(1+\lambda_j)^{-\varrho}$ with $\varrho=r d_s/2$, $d_s=1$, $\nu_X=\nu_V$. Seventeen values of $\Delta\nu/(d_s/2)$ in $[-0.5,1.5]$; 200 replicates; RMSE evaluated at 25 logarithmically spaced values of $N$.

\paragraph{E3 (finest-scale association).} $d_s=2$, $N=2^{18}$, $f^Z=\sigma_Z^2\kappa_Z^{2\nu_Z}(\kappa_Z^2+\lambda)^{-(\nu_Z+1)}$ and likewise for $H$; $c=1$, $\sigma_Z=\sigma_H=1$, $\kappa_Z=1$. Range sweep: $\nu_Z=\nu_H=1$, eleven values of $\kappa_H/\kappa_Z\in[1/4,4]$, 60 replicates. Regimes: $(\nu_Z,\nu_H)\in\{(1.5,1),(1,1),(1,1.5)\}$ with $\kappa_H=3$, 100 replicates; oracle GLS weights $1/g$.

\paragraph{E4 (resolution).} $d=1$, $\delta=1$; $f_X(\omega)=(1+\omega^2)^{-1}$ ($\nu_X=1/2$), $g(\omega)=0.5(1+\omega^2)^{-1}$, $\alpha(\omega)=0.5(1+(\omega/2)^2)^{-1}$, $\beta=1$. The constructions follow the proof of \cref{thm:resolution}: $\chi=1$ on $\abs\omega\le R$ with a $\cos^2$ taper to zero at $R+\pi/\delta$; $\tilde\chi=1$ on the principal cell with a taper to $1.25\pi/\delta$; and, for part (b), $\chi_1=\tilde\chi$ and $\chi_2$ its translates by $\pm2\pi/\delta$. In part (a) the cutoff $R$ is the smallest integer multiple of $\pi/\delta$ for which $g'>0$, found by bisection. In part (b) it is the smallest multiple of $100\pi/\delta$ for which $\max v^0/v^\delta<0.8$ on a grid of 2001 points of the closed cell, including the boundary, where $v^\delta$ is smallest. This gives $R=800\pi/\delta$. An independent computation with exact alias tails gives $\max v^0/v^\delta=0.87<1$ at this cutoff. Alias sums are truncated at $\abs m\le3000$ (identity checks) or $1500$ (cutoff searches). The independent check integrates the continuum spectral differences over $\R$ by the trapezoidal rule (step $\approx10^{-2}$ on the support of the construction, explicit tail formulas beyond it) for lags $0,\dots,10$. Block means use $\abs{\hat\psi_\delta}^2=\operatorname{sinc}^2(\delta\omega/2)$. The resolution floor (panel c) uses 30 values of $\delta\in[0.01,4]$ and alias sums with $\abs m\le2000$.

\paragraph{E7 (non-Gaussian residuals).} As E1 with $d_s=1$, $\nu_V=1/2$, $\varrho=1$, and $v_j=\sqrt{g_j}\,e_j$ with $e_j$ standard Gaussian, Laplace scaled to unit variance, or Student-$t_5$ scaled to unit variance; 200 replicates. The heavy-tailed comparison uses $\alpha=0$, oracle weights, $\nu_X=\nu_V+(1-p)/2$ for $p\in\{1.5,1,0.5\}$, and standard Cauchy $e_j$; it reports median absolute errors, since the RMSE is undefined.

\paragraph{E2 (2-torus with nugget).} Whittle--Mat\'ern fields on $[0,1)^2$ simulated on a $512\times512$ master grid by FFT, normalized to unit marginal variance, constant mode removed (equivalent to an intercept), and subsampled to $m\times m$ grids with $m\in\{16,\dots,512\}$. $\nu_V=1$, $\nu_X=\nu_V+\Delta\nu$; $\kappa=2\pi\cdot3$ (solid) or $2\pi\cdot8$ (dotted) for both fields; $\alpha(k)=0.5/(1+\lambda_k/(8\pi)^2)$; nugget variance $0.05$. GLS uses the residual spectrum at the principal frequencies (without aliasing correction) plus the nugget, computed in the DFT domain; this is a working GLS with a tail-equivalent spectrum, not exact oracle GLS for the subsampled data; 100 replicates.

\paragraph{E5 (kriged covariate).} $256\times256$ master grid; $\nu_X=\nu_V=1$, $\kappa=8\pi$; 36 stations on a regular $6\times6$ sub-lattice; simple kriging with the true covariance. The response is $Y=X+\alpha\ast X+V$ plus nugget $0.05$, with $\alpha\equiv0$ or $\alpha(k)=0.5/(1+\lambda_k/(4\pi)^2)$. Regression on the kriged surface uses the favourable working spectrum $g+f_X$ plus nugget; 150 replicates.

\paragraph{E6 (nonlinear covariate effect).} $d=1$, master grid $2^{16}$, subsamples $2^9,\dots,2^{16}$; $\kappa=8\pi$; $\nu_V=1/2$ and $\nu_X\in\{0.75,1.5\}$ (inside and outside, the frontier being $\nu_X=1$). Covariate $X=G$ (Gaussian) or $X$ the standardized $e^{G/2}$ (log-normal). Response $X+0.5(X^2-\overline{X^2})+0.5S[X]+0.3S[X^2-\overline{X^2}]+0.5V$ with $S=(1+\lambda/(4\pi)^2)^{-3}$ and nugget $0.01$; GLS on the features $(X,X^2-\overline{X^2})$ with the constant mode removed; 200 replicates. Because $S$ decays like $\lambda^{-3}$, this mechanism is Cameron--Martin negligible.

\paragraph{E8 (Poisson counts).} Here $\beta=0.5$; $d=1$, master grid $2^{15}$, $m\in\{2^7,\dots,2^{15}\}$ cells; latent $\eta=0.5X+\alpha\ast X+0.3V$ with $X$ of standard deviation $0.5$, $\alpha(k)=0.3/(1+\lambda_k/(4\pi)^2)$, $\nu_V=1/2$, $\nu_X\in\{0.5,1.5\}$. Counts are Poisson with means $e_m e^{\eta}$, with $\eta$ evaluated at the cell's grid point rather than integrated over the cell, and $e_m=4(m/128)^{1/2}$, so that the total expected count grows like $m^{3/2}$. The estimator is GLS on $\log\{(N+1/2)/e_m\}$ with working spectrum $0.09g$ and nugget $1/e_m$; 200 replicates.

\paragraph{E9 (nugget, point sampling).} $d=1$, master grid $2^{17}$, subsamples $n\in\{2^9,2^{11},\dots,2^{17}\}$; $\kappa=8\pi$; $\nu_V=0.5$, $\nu_X=0.9$; both fields normalized to unit variance with the constant mode removed; $\alpha(\lambda)=0.5/(1+\lambda/(4\pi)^2)$; nugget variance $0.01$. Working covariance: a unit-variance Whittle--Mat\'ern field with smoothness $\nu_h$ and range parameter $\kappa_h=\kappa$ (or $8\kappa$) plus the true nugget, through DFT-domain weights $1/(n^2h(\lambda)+n\cdot0.01)$; 300 replicates.

\paragraph{E10 (areal models).} Cycle graphs with $n$ nodes, whose Laplacian eigenvalues are $2-2\cos(2\pi k/n)$, $k=0,\dots,n-1$. Physical design: master grid of $2^{16}$ points on the unit 1-torus, $n\in\{2^7,\dots,2^{15}\}$, $\kappa=8\pi$, both fields normalized to unit variance with the constant mode removed, $(\nu_X,\nu_V)=(0.75,0.5)$ or $(1.5,1.25)$, $\alpha(\lambda)=0.5/(1+\lambda/(4\pi)^2)$ or $\alpha\equiv0$, no nugget. GLS is computed in the DFT domain, without the constant mode, with the weights $\lambda$ (ICAR), $\lambda/(1+\lambda)$ (BYM with $\sigma_u^2=\sigma_v^2=1$), $1-\rho+\rho\lambda$ (Leroux, with $\rho=0.9$ or $\rho_n=1/(1+\kappa^2/n^2)$), $(1-\rho+\rho\lambda/2)^2$ (SAR error with row-standardized $W$ and $\rho=0.4$) and $1/g$ at the principal frequencies (oracle); 100 replicates. Without confounding the RMSE at $n=2^{15}$ is $0.024$ (ICAR), $0.027$ (BYM), $0.224$ (Leroux, $\rho=0.9$) and $0.224$ (OLS) for $\nu_X=0.75$, and $0.166$, $0.167$, $0.318$ and $0.318$ for $\nu_X=1.5$. Spectral design: the model \eqref{eq:model} with $\lambda_j=j^2$, spectra $(1+\lambda_j)^{-(\nu+1/2)}$, $\alpha_j=0.5(1+\lambda_j)^{-1}$ or $0$, the first $n\in\{2^7,\dots,2^{17}\}$ modes and the sorted cycle eigenvalues $4\sin^2(\pi\lfloor j/2\rfloor/n)$ as graph eigenvalues; 200 replicates. With confounding and $\nu_X=0.75$ the exponents are $-0.262$ (ICAR), $-0.257$ (BYM), $-0.326$ (rescaled Leroux), $+0.001$ (Leroux, $\rho=0.9$), $0.000$ (SAR error and OLS) and $-0.262$ (oracle), against $-0.25$, $-0.25$, $-0.25$, $0$, $0$ and $-0.25$; with $\nu_X=1.5$ they lie between $-0.001$ and $0.000$ for ICAR, BYM and rescaled Leroux, and equal $-0.262$ for oracle GLS.

\end{document}